\documentstyle[amscd,amssymb,verbatim]{amsart}
\newcommand{\lrar}[1]{\begin{picture}(50,10)(-25,-5)                          
\put(-25,0){\vector(1,0){50}}
\put(0,5){\makebox(0,0)[b]{\mbox{$#1$}}}
\end{picture}}

\newcommand{\ldar}[1]{\begin{picture}(10,50)(-5,-25)
\put(0,25){\vector(0,-1){50}}
\put(5,0){\mbox{$#1$}}
\end{picture}}

\newcommand{\fin}{\operatorname{fin}}
\newcommand{\mA}{{\mathbb A}}
\newcommand{\res}{\operatorname{res}}
\newcommand{\smooth}{\operatorname{smooth}}
\newcommand{\card}{\operatorname{card}}
\newcommand{\GL}{\operatorname{GL}}

\newcommand{\vol}{\operatorname{vol}}

\renewcommand{\div}{\operatorname{div}}
\newcommand{\SO}{\operatorname{SO}}

\newcommand{\Sym}{\operatorname{Sym}}

\newcommand{\Tr}{\operatorname{Tr}}

\newcommand{\Lines}{\operatorname{Lines}}

\newcommand{\SL}{\operatorname{SL}}
\newcommand{\F}{{\Bbb F}}
\newcommand{\G}{{\Bbb G}}

\newcommand{\frg}{{\frak g}}

\newcommand{\lan}{\langle}
\newcommand{\ran}{\rangle}

\newcommand{\Aut}{\operatorname{Aut}}

\newcommand{\Frob}{\operatorname{Frob}}

\renewcommand{\Sp}{\operatorname{Sp}}

\renewcommand{\P}{{\Bbb P}}

\newcommand{\ga}{\gamma}
\newcommand{\de}{\delta}
\newcommand{\eps}{\epsilon}

\renewcommand{\ker}{\operatorname{ker}}

\numberwithin{equation}{section}

\newtheorem{thm}{Theorem}[subsection]
\newtheorem{prop}[thm]{Proposition}
\newtheorem{lem}[thm]{Lemma}
\newtheorem{cor}[thm]{Corollary}
\newenvironment{defi}{\vspace{3mm}\noindent
{\bf Definition.}}{\vspace{3mm}}
\newenvironment{rem}{\vspace{3mm}\noindent
{\bf Remark.}}{\vspace{3mm}}

\newcommand{\Pf}{\noindent {\it Proof}}
\newcommand{\id}{\operatorname{id}}

\newcommand{\ov}{\overline}

\newcommand{\ra}{\rightarrow}

\newcommand{\FF}{{\cal F}}

\newcommand{\SS}{{\cal S}}
\newcommand{\LL}{{\cal L}}
\newcommand{\MM}{{\cal M}}
\newcommand{\OO}{{\cal O}}

\newcommand{\Res}{\operatorname{Res}}

\newcommand{\norm}{\operatorname{norm}}
\renewcommand{\a}{\alpha}
\renewcommand{\b}{\beta}
\newcommand{\om}{\omega}
\newcommand{\De}{\Delta}
\newcommand{\la}{\lambda}
\newcommand{\th}{\theta}
\newcommand{\C}{{\Bbb C}}
\newcommand{\R}{{\Bbb R}}
\newcommand{\Z}{{\Bbb Z}}
\newcommand{\Q}{{\Bbb Q}}
\newcommand{\La}{\Lambda}
\newcommand{\Ga}{\Gamma}

\newcommand{\wt}{\widetilde}

\newcommand{\sub}{\subset}
\newcommand{\ed}{\qed\vspace{3mm}}
\newcommand{\triv}{\operatorname{triv}}
\newcommand{\grad}{\operatorname{grad}}

\title{Minimal representations: spherical vectors and automorphic functionals}
\author{D. Kazhdan and A. Polishchuk}
\thanks{The work of both authors was partially supported by NSF grants}

\begin{document}
\begin{abstract}
In the first part of this paper we study minimal
representations of simply connected simple split groups $G$ of type $D_k$ or $E_k$ over
local non-archimedian fields. Our main result
is an explicit formula for the spherical vectors in these representations.
In the case of groups over $\R$ and $\C$, such a formula
was obtained recently in \cite{KPW}. We also use our techniques to
study the structure of the space of smooth vectors in
the minimal representation. In the second part we consider groups $G$ as above
defined over a global field $K$. 
In this situation we describe the form of the automorphic functional on the 
minimal representation of the corresponding adelic group. 
\end{abstract}

\maketitle


\centerline{\sc Introduction}

\bigskip

The study of minimal representations of simple groups of types $D_k$ and $E_k$
over local non-archimedian
fields started with papers \cite{K} and \cite{KS}. In the former paper such a representation
was constructed explicitly for (not necessarily split) groups of type $D_4$.
In \cite{KS} an explicit realization of the minimal unitary representation
$\pi:G(F)\ra\Aut(R_2)$
in the split case was given, where $F$ is arbitrary local non-archimedian
field (in the archimedian case the construction is similar).
On the other hand,
minimal representations appear naturally in the theory of automorphic forms
as residues of Eisenstein series (see \cite{K}, \cite{GRS}).
It would be desirable to connect the global theory with the explicit local realizations.
To achieve this it is important to understand these local realizations better, namely,
to understand the structure of the space $R\sub R_2$ of smooth vectors and
to have explicit formulas for spherical vectors. In the case of real and
complex groups such formulas were obtained in the work \cite{KPW} that served as
a motivation for us. In the present paper we compute spherical vectors in
the non-archimedian case. The answer is given in Theorem \ref{sphericalthm} below.

In the case of the type $D_4$ our calculation uses a different realization of
the minimal representation in terms of functions on a quadratic cone and the results
of \cite{BK} that allow to identify the spherical vector in that realization.
The intertwining operator between the two realizations is given by some kind of Fourier
transform and we can easily calculate the image of the spherical vector
(this method can be generalized to the type $D_k$ with any $k\ge 4$). For all other types
our method of computation is roughly as follows.
First we analize the action of some Hecke operators
on Iwahori invariant vectors. Using the results of Savin in \cite{Savin}
and the calculations of
Lusztig in \cite{L1} we obtain the form of the spherical vector with one unknown constant.
Then we determine this constant by direct calculation involving the finite Fourier transform.

In the non-archimedian case we also obtain some information on the structure of the space of
smooth vectors $R$ consisting of all vectors in $R_2$ with open stabilizer.
Note that in the case of exceptional groups these results were obtained by Magaard and Savin in
\cite{MS}. Recall that the minimal representation is realized 
in the space of $L_2$-functions $f(y,x_0,\ldots,x_n)$ where
$x_0,\ldots,x_n\in F$
are variables of the Schr\"odinger representation of certain Heisenberg group $H$ contained in $G$
(see section \ref{constrsec}). This subgroup $H$ appears as a unipotent radical of some maximal 
parabolic subgroup $P\sub G$. The space of smooth vectors $R$ in this representation
contains as a subspace the space of all locally constant functions with compact support contained in the
set $y\neq 0$. To understand the quotient of $R$ by this subspace we define
$$\ov{f}(x_0,\ldots,x_n)=\lim_{y\ra 0} \psi(\frac{I_G(x_1,\ldots,x_n))}{x_0y})f(y,x_0,x)$$
where $I_G$ is certain cubic form that appears in the formulas for operators of $G$ acting
on $R$. We prove that for $x_0\neq 0$ the limit defining $\ov{f}(x_0,\ldots,x_n)$ exists
for all $f\in R$. Then we consider the variety $C$ defined as the graph of the gradient of
the rational function $I_G(x_1,\ldots,x_n)/x_0$. It turns out that for all $f\in R$
the functions $\ov{f}(x_0,\ldots,x_n)$ extend to locally constant sections 
$\res(f)$ of some line bundle $\LL$ over
$C(F)$. Moreover, $\LL$ is equipped with the natural action of the Levi subgroup $L(F)\sub P(F)$ and the obtained
map is $L(F)$-equivariant. We prove that all locally constant sections with compact support are obtained
in this way. The final piece of $R$ is the space of asymptotics of functions $\ov{f}$ near
zero. We show that this space coincides with the Jacquet functor $J(R)$ of $R$, 
i.e., the space of $H(F)$-coinvariants on $R$. It is known that 
$J(R)$ is the direct sum of a one-dimensional representation and of the space of smooth vectors in the minimal
representation of $L(F)$ (twisted by certain character).

Some of these results have analogues in the case $F=\R$ or $\C$.
Namely, let us denote by $\frg$ the Lie algebra of $G(F)$.
If $\pi:G(F)\ra \Aut(R_2)$ is a unitary
representation of G then the space of smooth vectors in it
is defined to be the subspace $R\sub R_2$ consisting of vectors $f\in  R_2$
which belong to the domain of the definition
of $\pi(u)$ for all $u\in U(\frg)$. It is well-known that $R$ is dense in $R_2$
(see \cite{Gard}). For every
$u\in U(\frg)$ we consider the norm on $R$ given by 
$\| f\|_u := \|\pi(u)f \|$ and introduce the topology on $R$ given by the
system of norms $(\|\cdot \|_u)$ .
Then $R$ is a complete locally convex topological space.
Let $\SS(G(F))$ be the  ring of Schwartz 
functions on $G(F)$.
We have an action $\pi:\SS(G(F))\otimes R\ra R: (\varphi,r)\mapsto \varphi*r$
and $\pi$ is surjective. Similar action can be defined for subgroups
of $G$. We define Jacquet functor
$J(R)=J_H(R)$ with respect to $H$ as the quotient of $R$
by the closure of the
span of vectors of the form $\varphi*f-(\int_{H(F)}\varphi(h)dh)f$, where
$f\in R$, $\varphi\in\SS(H(F))$. On the other hand,
we call a linear functional
$f^*$ on $R$ {\it tempered} (resp. $G'$-{\it tempered}, where $G'\sub G$
is a subgroup) if for every $r\in R$ the
functional  $\varphi\ra f^*(\varphi*r)$ on $\SS(G(F))$ (resp. $\SS(G'(F))$)
is continuous. The two definitions
are compatible in the following way: if
a functional $f^*$ is $H$-tempered and $H(F)$-invariant then
$f^*$ factors through $J(R)$. We refer to \ref{archsec} for
more details on the archimedian picture (which is still partially
conjectural).

In the second part of the paper
we use our formula for the spherical vectors to determine the form of an
automorphic (i.e., $G(K)$-invariant tempered)
functional on the minimal representation of the
adelic group $G(\mA)$, where
$K$ is a global field, $\mA=\mA_K$ is the ring of adeles for it.
We prove that every $G(K)$-invariant tempered functional
on the adelic minimal representation is proportional to the
functional of the form
$$\th(f)=\sum_{y\in K^*,x_0,\ldots,x_n\in K}f(y,x_0,\ldots,x_n)+|D|^{\frac{n+3}{2}}\cdot
\sum_{z\in C(K)}\res(f)+\a$$
where the second sum is defined using a natural trivialization of
the line bundle $\LL$ over $K$-points of $C$, $\a$ is the functional
that factors through the Jacquet functor $J(R)$, $|D|$ is the absolute value of the discriminant
in the number field case, $|D|=q^{2g-2}$ in the functional field case, where $g$ is the genus.
We study in more details the case when $G$ is of type $E_k$ 
and $K$ is the field of functions on a curve $X$ of genus $g$ over $\F_q$. In this case it is
more convenient to work with the twisted realization of the minimal representation of $G(\mA)$, where
all variables are (adelic) $1$-forms on $X$. The line bundle $\LL$ in this case is trivial but
has a non-trivial $L$-equivariant structure given by some character of $L$. Because of this
the second term of $\th$ can be rewritten in terms of the function
$\ov{f}^{\norm}(x_0,\ldots,x_n)=|x_0|^{s+1}\ov{f}(x_0,\ldots,x_n)$ that extends to a locally
constant function on $C$. We prove that in this realization the $G(K)$-invariant tempered functional
takes form
$$\th(f)=\sum_{y\in \om_K\setminus\{0\},(x_0,x)\in \om_K^{n+1}}f(y,x_0,x)+
q^{(2g-2)(s+1)}\cdot\sum_{z\in C(\om_K)}\ov{f}^{\norm}(z)+\a_1(f)+\a_2(f),$$
where $s=1,2,4$ for $G=E_6, E_7, E_8$,
$\om_K$ is the set of rational $1$-forms on our curve, 
the functional $\a_1$
factors through the trivial component of the Jacquet functor of $R$, $\a_2$ factors
through the minimal representation of the Levi subgroup of $P$.
The functionals $\a_1$ and $\a_2$ are uniquely determined by the constants
$\a_1(f_0)$ and $\a_2(f_0)$, where $f_0$ is the product of spherical vectors over all places.
We compute these constants in the cases $g=0$ and $g=1$. Basing on this computation
we conjecture that the following relation holds:
$$\a_1(f_0)=q^{(2g-2)(s+1)}L(X,q^s),$$
where $L(X,t)$ is the $L$-function of the curve $X$.
In a sequel to this paper we will show that there is a formula for $\a_1(f_0)$ as a product
of local factors (defined by certain integrals). It seems that these factors are computable and
in this way we plan to prove the above formula for $\a_1(f_0)$. This method should also provide
a formula for $\a_1(f_0)$ in the number field case. 

To complete our description of the automorphic
functional on the minimal representation in the number field case,
we need to know the Jacquet functor $J(R)$ of local
minimal representations at archimedian places.
We conjecture that the decomposition of
this Jacquet functor into irreducible representations of the
Levi subgroup has the same
structure as in the non-archimedian case. D.~Barbasch informed us that
he checked this for the case of series $D_k$. 
In section \ref{globalEsec} assuming the validity of this conjecture
we show that for $G$ of type $E_k$ the $G(K)$-invariant tempered functional 
can be written as the sum of $4$ terms $\th_1+\th_2+\a_1+\a_2$ as above.
Hence, in this case the calculation of this functional again reduces to
finding the constants
$\a_1(f_0)$ and $\a_2(f_0)$.

The paper is organized as follows: sections \ref{constrsec}-\ref{archsec} constitute part I
devoted to local results, while sections \ref{globalsec} and \ref{globalEsec} constitute part II
describing the global picture. In section \ref{constrsec} we review the standard 
model for the local minimal representation and formulate our main local result: the formula
for the spherical vector. In section \ref{D4sec} we prove this formula for $G=D_4$.
In section \ref{structuresec} we study properties of the map $f\mapsto\ov{f}$ and
describe the form of $\ov{f}$ for $L(\OO)$-invariant vectors $f$. In section \ref{heckesec}
we use the action of Hecke algebra to derive the formula for the spherical vector with
one unknown constant. In section \ref{fouriersec} we find this constant by explicit
computation involving the Fourier transform over the residue field.
In section \ref{revisitsec} we apply our formula for spherical vectors to obtain the
description of the space of smooth vectors in local minimal representations.
In section \ref{archsec} we (partially) describe the archimedian picture.
In section \ref{globalsec} we determine 
the general form of the automorphic functional on the adelic
minimal representation. Finally,
in section \ref{globalEsec} we get more precise results
about this functional in the case of $G=E_k$.

\noindent
{\it Acknowledgment.} We are grateful to D.~Barbasch, A.~Braverman and P.~Etingof 
for useful discussions, and to B.~Gross and J.~Harris for their help with the
proof of Lemma \ref{ellnumlem2}. Special thanks are to G.~Savin for valuable
remarks on the first draft of the paper.
The second author would like to thank Max-Planck-Institut
f\"ur Mathematik and Institut des Hautes \'Etudes Scientifiques  
for the hospitality during the final stages of the preparation of this paper.

\bigskip

\centerline{\bf Part I. Local results}

\section{Formula for the spherical vector}\label{constrsec}

Let $G$ be a simply connected simple split group of type $D_k$ or $E_k$ over a local
non-archimedian field $F$. 
We denote by $\Phi$ the corresponding root system, by
$\Phi^+$ the set of positive roots, by
$\Delta$ the set of simple roots and by $\om$
the highest positive root.
We fix a splitting of $G$ (see \cite{Borel}). It gives
for every root $\a$ an isomorphism 
$t\mapsto e_{\a}(t)$ between $\G_a$ and the 
one-dimensional unipotent subgroup in $G$ corresponding to $\a$. We have
$$[e_{\a}(t_1),e_{\b}(t_2)]=\cases e_{\a+\b}(c(\a,\b)t_1t_2),\ 
& \a+\b\in\Phi,\\
1,\ &\text{otherwise,}\endcases$$
where $c(\a,\b)=\pm 1$ is some sign (see e.g. \cite{KS}, sec.1).
The important properties of this sign function are:
$$c(\a,\b)=c(-\a,\b),$$
$$c(\a,\b)c(\b,\a)=(-1)^{\lan\a,\b\ran},$$
$$c(\a,\b+\ga)=c(\a,\b)c(\a,\ga),$$
$$c(\a,\a)=-1\ \text{for}\ \a\in\Phi.$$
Let us denote by $\OO\sub F$ the valuation ring in $F$.
We choose an additive character $\psi:F\ra\C^*$ of order $0$, i.e., such
that $\psi(\OO)=1$ and $\psi(\pi^{-1}\OO)\neq 1$, where $\pi\in F$
is a uniformizer. The corresponding self-dual
measure $dx$ on $F$ is such that $\vol(\OO)=1$.

\subsection{Heisenberg group and standard realization}\label{standsec}

Let $\b_0$ be the simple root to which the affine root attaches on
the extended Dynkin diagram. The following lemma is easily checked
case by case.

\begin{lem}\label{rootlem1} 
The fundamental weight corresponding to $\b_0$ coincides with $\om$.
\end{lem}

Let $P\sub G$ be the maximal parabolic subgroup
corresponding to $\Delta\setminus\{\b_0\}$. We denote by
$H$ the unipotent radical of $P$. The Lie algebra of $H$ is spanned
by the roots in the set $\{\om\}\cup\Sigma$, where $\Sigma$ is the set
of all roots $\ga$ such that $\om-\ga$ is also a root.
Equivalently, $\Sigma=\{\ga\in\Phi:\ \lan\om,\ga\ran=1\}$.
We have a natural involution $*:\ga\mapsto\om-\ga$ on $\Sigma$.
Let us define the subset $\Pi\subset\Sigma$ by setting
$\Pi=\{\b_0\}\cup\Pi_1$, where
$$\Pi_1=\{\b\in\Sigma:\ \lan\b_0,\b\ran=1\}=
\{\b\in\Sigma:\ \b-\b_0\in\Phi\}.$$
Then $\Sigma$ is the disjoint union of $\Pi$ and $*(\Pi)$.

We have $P=LH$, where $L$ is the Levi subgroup: 
$L=G_1\b_0^{\vee}(\G_m)$, where $G_1=[L,L]$ is
a semisimple group with simple roots $\Delta\setminus\{\b_0\}$.
The center $Z(L)$ of $L$ is generated by the fundamental coweight 
$\om_{\b_0}^{\vee}(\G_m)$ corresponding to $\b_0$. 
On the other hand, $\om$ (considered as a character of the maximal
torus) extends to the character $\om:L\ra\G_m$ such that
$G_1=\ker\om$.

The center of $H$ is the subgroup $Z(H)=e_{\om}(\G_a)\subset G$.
The group $H$ is the central extension of the vector space
$V=H/Z(H)$ by $Z(H)$.
Moreover, the commutator on $H$ induces a symplectic form $E$ on
$V$. By the definition, we have
\begin{equation}\label{symptransfor}
E(gv_1,gv_2)=\om(g)E(v_1,v_2),
\end{equation}
where $g\in L$, $v_1,v_2\in V$.
Let $\La\sub V$ (resp. $\La'\sub V$) be the span of the one-parametric
subgroups corresponding to $\Pi$ (resp. $*\Pi$). 
Then $\La$ and $\La'$ are Lagrangian
subspaces in $V$. We denote the elements of $\Pi_1$ by $\b_1,\ldots,\b_n$,
so that $\Pi=\{\b_0,\b_1,\ldots,\b_n\}$ and we set 
$\ga_i=*(\b_i)$ for $i=0,\ldots,n$. In the case $G=D_k$ with $k\ge 5$
we make the following special choice of $\b_1$. 
If we remove the vertex $\b_0$ from the Dynkin diagram of $D_k$,
we get a disjoint union of $A_1$ and $D_{k-2}$. Let us denote by
$\a_1$ the simple root corresponding to $A_1$ in the obtained diagram.
Then we set $\b_1=\b_0+\a_1$. For every root $\a$ let us set
$e_{\a}=e_{\a}(1)$.
Then $(e_{\b_0},\ldots,e_{\b_n},e_{\ga_0},\ldots,e_{\ga_n})$ 
is a basis of $V$ over $F$ and we denote by 
$(x_0,\ldots,x_n,x'_0,\ldots,x'_n)$ the corresponding
coordinates on $V$. Note that the symplectic form in this basis has
the following form: $E(e_{\b_i},e_{\ga_j})=\delta_{ij}c(\b_i,\ga_j)$.
Let us denote $\eps_i:=E(e_{\b_i},e_{\ga_i})=\pm 1$.
Using the properties of $c(\cdot,\cdot)$
we obtain the following expressions for these signs:
$$\eps_i=c(\b_i,\ga_i)=-c(\b_i,\om)=c(\ga_i,\om).$$
When $\b$ is one of the roots $\b_i$ we will also use the notation
$x_{\b}$ for $x_i$.

The minimal representation of $G(F)$ is realized in the space of
$L^2$-functions $f(y,x_0,\ldots,x_n)$. The action of $H(F)$ is given
as follows:
$$e_{\om}(t)f=\psi(ty)f,$$
$$e_{\ga_i}(t)f=\psi(\eps_i t x_i)f,$$
$$e_{\b_i}(t)f(y,x_0,\ldots,x_i,\ldots,x_n)=f(y,x_0,\ldots,x_i+ty,\ldots,x_n),
$$
where $t\in F$, $i=0,\ldots,n$.
The important role is played by the following two operators:
$$Sf(y,x_0,\ldots,x_n)=\int_{F^{n+1}}f(y,\wt{x}_0,\ldots,\wt{x}_n)
\psi(y^{-1}(\wt{x}_0x_0+\ldots+\wt{x}_nx_n))
|y|^{-\frac{n+1}{2}}d\wt{x}_0\ldots d\wt{x}_n,$$
\begin{equation}\label{operatorA}
Af(y,x_0,\ldots,x_n)=\psi(-\frac{I(x_1,\ldots,x_n)}{x_0y})
f(-x_0,y,x_1,\ldots,x_n),
\end{equation}
where $I=I_G$ is the cubic form
$$I(x_1,\ldots,x_n)=
\sum_{i<j<k:\b_i+\b_j+\b_k=\b_0+\om}c(\b_i,\b_j)c(\b_i,\b_k)c(\b_j,\b_k)
c(\b_0,\om)x_ix_jx_k.$$
The operator $A$ represents the Weyl generator 
$s_{\b_0}=e_{\b_0}(1)e_{-\b_0}(1)e_{\b_0}(1)$, while
$S$ represents the longest element in the Weyl group of $L$ (or rather
its canonical lifting to $G$).
The action of an element $t$ of the maximal torus $T\sub G$ is given by
$$tf(y,x_0,\ldots,x_n)=|\om(t)\prod_{i=0}^n\ga_i(t)|^{\frac{1}{2}}\cdot
f(\om(t)y,\ga_0(t)x_0,\ldots,\ga_n(t)x_n).$$
In particular, $Z(L)$ acts as follows:
$$\om^{\vee}_{\b_0}(c)f(y,x_0,\ldots,x_n)=
|c|^{\frac{n+3}{2}}f(c^2y,cx_0,\ldots,cx_n),$$
where $c\in F^*$. On the other hand, the coroot corresponding to $\b_0$
acts by
$$\b_0^{\vee}(c)f(y,x_0,x_1,\ldots,x_n)=f(cy,c^{-1}x_0,x_1,\ldots,x_n).$$
Note that the operator $A^2$ corresponds to the action of
$s_{\b_0}^2=\b_0^{\vee}(-1)$.

We denote by $R$ the space of smooth vectors in the minimal representation of
$G(F)$. Later we will discuss the structure of this space in more details.
Now we restrict ourselves to the following simple observations.

\begin{lem}\label{supportlem} 
(a) Every $f\in R$ has bounded support.

\noindent
(b) For every $f\in R$ and every $\eps>0$ the restriction of $f$ to 
the set $\{(y,x_0,\ldots,x_n):|y|\ge\eps\}$ is locally constant.

\noindent
(c) The space $R$ contains all locally constant functions $f$
with compact support contained in $F^*\times F^{n+1}$.
\end{lem}

\Pf . (a) This follows from the invariance of $f$ under the action of
$e_{\om}(t)$ and $e_{\ga_i}(t)$, $i=0,\ldots,n$,
for sufficiently small $t\in F$.

\noindent
(b) This follows from the invariance of $f$ under the action of
$e_{\b_i}(t)$, $i=0,\ldots,n$, for sufficiently small $t\in F$.

\noindent
(c) It is clear that such a function $f$ is invariant with respect
to some open subgroup in $P(F)$ (since it acts essentially through
Weil representation). Assume first that the support of $f$ is contained
in the set $y\neq 0, x_0\neq 0$. Then $Af$ also has compact support in
$F^*\times F^{n+1}$, so we get that $Af$ is also invariant under
the action of some open subgroup in $P(F)$. This implies that $f$
is a smooth vector. The general case is reduced to this one,
using the action of $L(F)$.
\ed

Our main result is the following formula for the spherical vector 
$f_0$ in $R$. Note that because of the invariance of $f_0$ with respect to
$A$, it suffices to compute $f_0(y,x_0,\ldots,x_n)$ for $|y|\le |x_0|$.

\begin{thm}\label{sphericalthm} Let $f_0\in R$ be the spherical vector normalized by
the condition $f(1,0,\ldots,0)=1$. Then one has the following formulas
for $f_0(y,x_0,x_1,\ldots)$, where $|y|\le |x_0|$. 

\noindent
(i) Case $G=E_k$.
$$f_0(y,x_0,x_1,\ldots)=
\cases \psi(-\frac{I}{yx_0})|x_0|^{-s-1}\cdot
\frac{q^s |(x,\grad(\frac{I}{x_0})x)|^{-s}-1}{q^s-1}, 
&|(x,\grad(\frac{I}{x_0})x)|\le 1,\\
0, &\text{otherwise},\endcases
$$
where $x=(x_0,x_1,\ldots)$; $s=1,2,4$ 
for $k=6,7,8$, respectively, $q$ is the number of elements in $\OO/\pi\OO$.

\noindent
(ii) Case $G=D_k$, where $k\ge 5$. 
\begin{align*}
&f_0(y,x_0,x_1,\ldots)=\\
&=\cases\psi(-\frac{I}{yx_0})
|x_0|^{-1}\max(1,|\frac{x_1}{x_0}|)^{k-4}\cdot
\frac{q^{k-4}|(x,\grad(\frac{I}{x_0})x)|^{-k+4}-1}{q^{k-4}-1}, 
&|(x,\grad(\frac{I}{x_0})x)|\le 1,\\
0, &\text{otherwise}.\endcases 
\end{align*}

\noindent
(iii) Case $G=D_4$.
$$f_0(y,x_0,x_1,x_2,x_3)=
\cases \psi(-\frac{x_1x_2x_3}{x_0y})|x_0|^{-1}\cdot
[1+v(x,\grad(\frac{x_1x_2x_3}{x_0})x)], &v(x_0,x_1,\ldots)\ge 0,\\
0, &\text{otherwise},\endcases
$$
where for every vector $y$ we denote by $v(y)\in\Z$ the valuation of $y$
(so that $|y|=q^{-v(y)}$).
\end{thm}

The proof of this theorem will constitute sections \ref{D4sec}--\ref{fouriersec}.
We end this section with several useful technical lemmas.
First, we are going to look at the action of some subgroups of $G$ on $R$.
Let us consider the semisimple subgroup $G_0\sub G_1$
corresponding to the root subsystem in $\Phi$ spanned by the
simple roots not attached to $\b_0$. 
We have the following table:
\begin{equation}
\begin{array}{ccc}
G   & G_1                 & G_0       \nonumber\\
D_k & A_1\times D_{k-2}   & D_{k-3}       \nonumber\\
E_6 & A_5                 & A_2\times A_2 \nonumber\\  
E_7 & D_6                 & A_5           \nonumber\\
E_8 & E_7                 & E_6 \nonumber
\end{array}
\end{equation}
Since, $G_0$ centralizes
$Z(H)$ and the
action of $G_0$ on $V$ preserves the subspaces $\La$ and $\La'$,
the corresponding operators on $R$ are
$$g_0f(y,\la)=f(y,g^{-1}\la).$$
Now let us consider the Borel subgroup $B\sub L$ 
corresponding to positive roots in $L$, and let $U=[B,B]\sub B$ 
be the unipotent radical of $B$.
It is easy to see that the action of $B$ on $V$ preserves the subspace
$\La'\sub V$. Furthermore, the element $e_{\ga_0}\in V$ is invariant under
the action of $U$. Hence, the action of $U$ on $V$ leaves the
coordinate $x_0$ on $V$ invariant.  
Now for every $g\in U(F)$ we have
$g(\la,\la')=(A_g\la,(A_g^*)^{-1}(\la'+B_g\la))$, where $\la\in\La$,
$\la'\in\La'$, $A_g$ is an operator on $\La$, $A_g^*$ is
the dual operator on $\La'$, $B_g$ is a self-dual linear map from 
$\La$ to $\La'$. The corresponding operator on $R$ is
$$gf(y,\la)=\psi(\frac{E(B_gA_g^{-1}\la,A_g^{-1}\la)}{2y})
f(y,A_g^{-1}\la),$$
where $\la\in\La$, $E$ is the symplectic form on $V$.
In the following lemma we consider
$I=I$ as a function on $\La$ by setting 
$I(x_0,x_1,\ldots,x_n)=I(x_1,\ldots,x_n)$.

\begin{lem}\label{cubiclem1} 
For every $g\in U(F)$ one has
$$\frac{E(\la,B_g(\la))}{2}=\frac{I(A_g\la)-I(\la)}{x_0(\la)},
$$
where $\la\in \La$.
\end{lem}

\Pf . Since $B_{g_1g_2}=A^*_{g_2}B_{g_1}A_{g_2}+B_{g_2}$
and the action of $U$ preserves the coordinate $x_0$, it suffices
to check the identity in question for $g=e_{\a}(t)$, where
$\a$ is a simple root $\neq\b_0$. If $\lan\a,\b_0\ran=0$, then
$\a$ belongs to the root system of $G_0$. In this case $B_g=0$, so we
have to prove that the operator $A_g$ preserves $I$. This follows 
easily from the fact that the operator $A$ representing $s_{\b_0}$ 
commutes with $g$. Thus, we are reduced to the case $\lan\a,\b_0\ran=-1$. 
In this case the proof is obtained by writing explicitly the equality
of operators in the minimal representation corresponding to the identity
$$s_{\b_0}e_{\a}(t)s_{\b_0}^{-1}=[e_{\b_0}(1),e_{\a}(t)].$$
\ed

\begin{lem}\label{cubiclem2} 
For every $t$ in the maximal torus of $G$ one has
$$I(\ga_1(t)x_1,\ldots,\ga_n(t)x_n)=\ga_0(t)\om(t)I(x_1,\ldots,x_n).$$
\end{lem}

\Pf . The cubic form $I$ contains only monomials of the form
$x_ix_jx_k$ such that $\b_i+\b_j+\b_k=\b_0+\om$. For such $(i,j,k)$ 
one has also $\ga_i+\ga_j+\ga_k=\ga_0+\om$.
\ed

\begin{lem}\label{rootlem2} Let us denote $\de=\om+\sum_{i=0}^n\ga_i$.
Then in the case $G=E_k$ one has
$$\de=2(s+1)(\om+\ga_0),$$
where $s=1,2,4$ for $k=6,7,8$. On the other hand, for $G=D_k$ ($k\ge 4$) 
one has
$$\de=(k-2)\om+2\ga_0+(k-4)\a_1,$$
where for $k\ge 5$ we denote by $\a_1$ the simple root attached only to
$\b_0$.
\end{lem}

The proof is a straightforward computation.

\subsection{Twisted realizations}\label{twistedsec}

Note that the additive character $\psi$ above was assumed to be of order $0$.
Now we are going to explain how to construct a twisted version of the standard realization
associated with an arbitrary non-trivial additive charater of $F$. 

Let $M$ be a free $\OO$-module of rank $1$ and let 
$\psi:F\otimes_{\OO} M\ra\C^*$ be a character which is trivial
on $M$ and non-trivial on $\pi^{-1}M$. Then one can
realize the minimal representation of $G(F)$ in $L_2$-functions of
$y,x_0,\ldots,x_n\in F\otimes_{\OO} M$ using essentially the same formulas as in
section \ref{standsec}. One just has to notice that the arguments of $\psi$
in these formulas are always expressions that are homogeneous of degree $1$
in the variables, so they can be considered as elements of $F\otimes_{\OO} M$.
Also, the norm on $F\otimes_{\OO} M$ is defined by $|a\otimes m|=|a|$ for
$a\in F$, $m\in M\setminus\pi M$.
Finally, the measure $|y|^{-\frac{1}{2}}dx$ appearing in the formula
for the action of $S$ should be replaced by the self-dual measure on $F\otimes_{\OO} M$
with respect to the biadditive character $(x_1,x_2)\mapsto \psi(x_1x_2/y)$ on this group
(see \cite{Weil-BNT}, VII-2). 

The formula of Theorem \ref{sphericalthm} for 
the spherical vector is still valid in this realization (with the above definition of
norm on $F\otimes_{\OO} M$). Note that this spherical vector is normalized by the 
condition $f_0(1\otimes m, 0,\ldots, 0)=1$ for $m\in M\setminus\pi M$.

Now let $\psi:F\ra\C^*$ be an arbitrary non-trivial character. Then
we can take $M=\pi^{-\nu}\OO\sub F$, where $\nu$ is the order of $\psi$, so
that $\psi$ is trivial on $\pi^{-\nu}\OO$ and non-trivial on $\pi^{-\nu-1}\OO$.
Using the embedding of $M$ in $F$ we can identify $F\otimes_{\OO} M$ with $F$.
Then the realization of the minimal representation associated with $M$ can be
considered as a realization in $L_2$-functions of $y,x_0,\ldots,x_n\in F$.
Let us denote the space of smooth vectors in this realization by $R_{\psi}$.
Of course, it coincides with $R$ but the action of $G(F)$ on $R_{\psi}$
is different. Namely, let $\psi_0$ be a character of order $0$ used to define $R$.
Then there exists $a\in\pi^\nu\OO^*$ such that $\psi(t)=\psi_0(at)$. The map
$f(y,x_0,\ldots,x_n)\mapsto f(ay,ax_0,\ldots,ax_n)$ is a $G(F)$-isomorphism from
$R$ to $R_{\psi}$.
The spherical vector $f_0$ in the realization $R_{\psi}$ is supported on the
$(\pi^{-\nu}\OO)^{n+2}$. We always normalize it by the condition that
$f_0(\pi^{-\nu},0,\ldots,0)=1$. Then $f_0$ is given by the formula of Theorem \ref{sphericalthm}
with $|x|$ replaced by $|\pi^n x|$.

We will need twisted standard realizations for global applications in the second part of the paper.
Until then we will always assume that $\psi$ is of order $0$.

\section{Case of $D_4$}\label{D4sec}

Our general strategy of finding the spherical vector $f_0$ is based
on the study of the action of certain double Iwahori coset on
Iwahori-invariant vectors in $R$.
The case $G=D_4$ differs from all the others
in that the action of this coset is not semisimple.
One can overcome this difficulty, but we prefer to give a completely
different derivation of the formula for $f_0$ in this case.
We are going to use a realization of the minimal representation of $D_4$
on the space of functions on the quadratic cone 
$u_1v_1+u_2v_2+u_3v_3=0$ (similar realization exists for all $D_n$). 
Namely, the subgroup $D_3\subset D_4$ acts
geometrically by linear transformations of variables preserving the
quadratic form. One more Weyl generator acts by the partial Fourier
transform 
\begin{align*}
& w\phi(u_1,u_2,u_3,v_1,v_2,v_3)=
|u_1|^{-2}\int_{u'_2,u'_3,v'_2,v'_3}
\psi(\frac{u'_2v_2+u'_3v_3+u_2v'_2+u_3v'_3}{u_1})\times\\
&\phi(-u_1,u'_2,u'_3,\frac{u'_2v'_2+u'_3v'_3}{u_1},v'_2,v'_3)du'_2du'_3
dv'_2dv'_3.
\end{align*}
It follows from the results of Braverman and Kazhdan in \cite{BK} that 
the spherical vector in this realization is the function
$$\phi_0(x)=\delta_{\OO}(x)\frac{|x|^{-1}-q^{-1}}{1-q^{-1}},$$
where for a vector $x=(u_1,u_2,u_3,v_1,v_2,v_3)$ we denote by
$|x|$ the maximum of norms of the coordinates, $\delta_{\OO}$
is the characteristic function of the set where all coordinates are integers 
(the normalization is chosen in such a way that $\phi_0(x)=1$ for $|x|=1$).
More precisely, we use the fact that
our quadratic cone coincides with the basic affine
space for the group $\SL_3$ and that the generalized Fourier
transforms $\Phi_w$ acting on functions on this space (see \cite{BK}),
are part of the action of the Weyl group of $D_4$. 
Theorem 3.13 of \cite{BK} describes an explicit basis in the space
of $K$-invariants for some compact subgroup $K$ in $D_4$ and the action
of the operators $\Phi_w$ on it. Combining this description with
the fact that the spherical vector is supported on points
with integer coordinates, we obtain the formula above. 

The intertwining operator from this realization to the standard one
is given by the formula
$$\SS\phi(y,x_0,x_1,x_2,x_3)=\int_{yy'+x_0x'_0+x_1x_2=0}
\psi(\frac{x_3x'_0}{y})\phi(y,y',x_0,x'_0,x_1,x_2)\frac{dx'_0}{|y|}.$$
Thus, we can compute the spherical vector $f_0$ in the standard 
realization using the formula $f_0=\SS\phi_0$. In other words, we have
to compute the following integral
\begin{align*}
&f_0(y,x_0,x_1,x_2,x_3)=\\
&\int_{y',x'_0\in\OO:\ yy'+x_0x'_0+x_1x_2=0}
\psi(\frac{x_3x'_0}{y})\frac{|(y,y',x_0,x'_0,x_1,x_2)|^{-1}-q^{-1}}{1-q^{-1}}
\cdot\frac{dx'_0}{|y|}
\end{align*}
for $y,x_0,x_1,x_2\in\OO$, $x_3\in F$.
First let us consider the case $|x_0|\le |y|$.
Then from the condition $y'=-\frac{x_0x'_0+x_1x_2}{y}$
we see that the domain of integration is empty unless $\frac{x_1x_2}{y}\in\OO$.
On the other hand, if $\frac{x_1x_2}{y}\in\OO$ then $y'$ is automatically
an integer. Furthermore, in this case $|(y,x_0)|=|y|$ and
$|(x'_0,y')|=|(x'_0,\frac{x_1x_2}{y})|$, so we can write
$$f_0(y,x_0,x_1,x_2,x_3)=\int_{x'_0\in\OO}
\psi(\frac{x_3x'_0}{y})
\frac{|(y,\frac{x_1x_2}{y},x'_0,x_1,x_2)|^{-1}-q^{-1}}{1-q^{-1}}
\cdot\frac{dx'_0}{|y|}.$$
Let us set $a=|(y,\frac{x_1x_2}{y},x_1,x_2)|$. Note that by our
assumptions $a\le 1$. Now we can split the domain of the above integral
into two pieces: $|x'_0|\le a$ and $a<|x'_0|\le 1$, so that
\begin{align*}
&|y|f_0(y,x_0,x_1,x_2,x_3)=\\
&\frac{a^{-1}-q^{-1}}{1-q^{-1}}
\int_{|x'_0|\le a}\psi(\frac{x_3x'_0}{y})dx'_0+
\int_{a<|x'_0|\le 1}\psi(\frac{x_3x'_0}{y})
\frac{|x'_0|^{-1}-q^{-1}}{1-q^{-1}}dx'_0.
\end{align*}
The first integral can be computed immediately:
$$I_1=\int_{|x'_0|\le a}\psi(\frac{x_3x'_0}{y})dx'_0=\cases
0, & |\frac{x_3}{y}|>a^{-1},\\
a, & |\frac{x_3}{y}|\le a^{-1}.\endcases $$
The second integral an be rewritten as follows:
$$I_2=\int_{a<|x'_0|\le 1}\psi(\frac{x_3x'_0}{y})
\frac{|x'_0|^{-1}-q^{-1}}{1-q^{-1}}dx'_0=
\sum_{b\in q^{\Z}: a<b\le 1}\frac{b^{-1}-q^{-1}}{1-q^{-1}}\cdot
\int_{|x'_0|=b}\psi(\frac{x_3x'_0}{y})dx'_0.$$
Now we use the fact that for $b\in q^{\Z}$ and $t\in F^*$ one has
$$\int_{|x|=b}\psi(tx)dx=\cases 
0, & |t|>qb^{-1},\\
-q^{-1}b, & |t|=qb^{-1},\\
(1-q^{-1})b, & |t|\le b^{-1}.\endcases $$
It follows that $I_2=0$ unless $|\frac{x_3}{y}|\le a^{-1}$.
Now let us consider two cases: (1) $|\frac{x_3}{y}|\le 1$,
(2) $1<|\frac{x_3}{y}|\le a^{-1}$. In the first case, denoting
$a=q^{-n}$ (where $n\ge 0$) we obtain
$$I_2=\sum_{0\le i<n}(1-q^{-i-1})=n-\frac{q^{-1}(1-a)}{1-q^{-1}},$$
and therefore
$$|y|f_0(y,x_0,x_1,x_2,x_3)=1+v(y,\frac{x_1x_2}{y},x_1,x_2).$$
In the second case, setting $m=v(\frac{y}{x_3})$ we get
$$I_2=\sum_{m\le i<n}(1-q^{-i-1})-\frac{q^{-1}}{1-q^{-1}}
(1-|\frac{y}{x_3}|)=n-m-\frac{q^{-1}(1-a)}{1-q^{-1}}$$
provided that $m\le n$ (otherwise, $I_2=0$).
Hence, in this case
$$|y|f_0(y,x_0,x_1,x_2,x_3)=1+v(y,\frac{x_1x_2}{y},x_1,x_2)-v(\frac{y}{x_3})=
1+v(x_3,\frac{x_1x_2x_3}{y^2},\frac{x_1x_3}{y},\frac{x_1x_2}{y})$$
provided that all the components of the last vector are integers
(otherwise, we get zero).
We can unify these answers by writing
$$f_0(y,x_0,x_1,x_2,x_3)=|y|^{-1}\cdot
[1+v(y,x_1,x_2,x_3,\frac{x_1x_2}{y},\frac{x_1x_3}{y},
\frac{x_2x_3}{y},\frac{x_1x_2x_3}{y^2})],$$
if all the components of the argument of $v$ are integers
and is zero otherwise.
The formula in the case $|y|\le|x_0|$ follows from the invariance of
$f_0$ with respect to the operator $A$.
Here is the resulting formula for $f_0$:
\begin{align*}
&f_0(y,x_0,x_1,x_2,x_3)=\\
&\cases |y|^{-1}\cdot 
[1+v(y,x_1,x_2,x_3,\frac{x_1x_2}{y},\frac{x_1x_3}{y},
\frac{x_2x_3}{y},\frac{x_1x_2x_3}{y^2})], & |x_0|\le |y|, 
v(y,x_1,\ldots)\ge 0\\
\psi(-\frac{x_1x_2x_3}{x_0y})|x_0|^{-1}\cdot
[1+v(x_0,x_1,x_2,x_3,\frac{x_1x_2}{x_0},\frac{x_1x_3}{x_0},
\frac{x_2x_3}{x_0}),\frac{x_1x_2x_3}{x_0^2})], & |y|\le |x_0|,
v(x_0,x_1,\ldots)\ge 0,\\
0, & \text{otherwise}.
\endcases
\end{align*}

\section{The space of smooth vectors in the minimal representation}
\label{structuresec}

In this section we take a closer look at the space of
smooth vectors in the minimal representation. Then we apply the
results of this study to derive some information on $L(\OO)$-invariant
smooth vectors. In the case of exceptional groups most of the results proved in
\ref{Levisec} and \ref{conesec} can be found in section 6 of \cite{MS}.

\subsection{Representation of the Levi subgroup and
the Lagrangian subvariety}\label{Levisec}

Let us look at the action of $L$ on $V$. We already know how
the maximal torus acts on $V$, so it remains to determine
the action of the commutator subgroup $G_1\sub L$.

\begin{lem} Assume that $G=E_k$ and let $\a_0$ be the unique simple
root attached to $\b_0$. Then $V$ is the fundamental representation
of $G_1$ corresponding to $\a_0$. In the case $G=D_k$, $k\ge 5$,
we have a decomposition $G_1=\SL_2\times \SO(2k-4)$.
Then $V$ is the tensor product of the standard representation of $\SL_2$
with the standard representation of $\SO(2k-4)$.
\end{lem} 

\Pf . This follows from the fact that $\lan \ga_0,\a\ran=0$ for all
simple roots $\a$ that are not attached to $\b_0$ while
$\lan\ga_0,\a\ran=1$ for simple roots attached to $\b_0$.
\ed


Let us denote by $Q_{\ga_0}\sub L$ the stabilizer subgroup of  
the line spanned by $e_{\ga_0}$. 


\begin{lem}\label{parablem}
The group $P_1=Q_{\ga_0}\cap G_1$ is a parabolic subgroup in
$G_1$ containing $B\cap G_1$. In the case $G=E_k$, it is
the maximal parabolic corresponding to the simple root $\a_0$.
In the case $G=D_k$ we have $G_1=\SL_2\times\SO(2k-4)$ 
and $P_1$ is the product of the Borel subgroup $B\cap\SL_2\sub\SL_2$
with the stabilizer subgroup of the line spanned by the
highest weight vector in the standard representation of $\SO(2k-4)$.
\end{lem}

For every algebraic group $\Ga$
over $F$ we denote
by $X(\Ga)$ the group of algebraic homomorphisms $\Ga\ra\G_m$
(defined over $F$).

\begin{lem}\label{charlem}
Let $T\sub G$ be the maximal torus, $T_0\sub T$ be the maximal torus
of $G_0$. The natural homomorphism
$$T(F)\ra Q_{\ga_0}(F)/[Q_{\ga_0}(F),Q_{\ga_0}(F)]$$
is surjective. Similar statement is true for the algebraic
homomorphism $T\ra Q_{\ga_0}/[Q_{\ga_0},Q_{\ga_0}]$.
The image of the induced embedding
$X(Q_{\ga_0})\ra X(T)$ consists of all characters of $T$
that are trivial on $T_0$. Same statements hold
for $Q_{\b_0}$, the stabilizer subgroup of the line spanned by $e_{\b_0}$.
\end{lem}

\Pf . We have a decomposition of $Q_{\ga_0}(F)$ into
a semi-direct product of the normal subgroup 
$P_1(F)=Q_{\ga_0}\cap G_1(F)$ and the one-dimensional
torus $\b_0^{\vee}(F^*)$. Therefore, we have an exact sequence
$$0\ra P_1(F)/[Q_{\ga_0}(F),P_1(F)]\ra
Q_{\ga_0}(F)/[Q_{\ga_0}(F),Q_{\ga_0}(F)]\stackrel{\om}\ra F^*\ra 0$$
which has a splitting induced by $\b_0^{\vee}$.
It remains to observe that 
$$P_1(F)/[P_1(F),P_1(F)]=L_1(F)/G_0(F)=T_1(F)/T_0(F),$$ 
where $L_1\sub P_1$ is the Levi subgroup,
$T_i=T\cap G_i$. Hence, the natural map
$$T_1(F)/T_0(F)\ra P_1(F)/[Q_{\ga_0}(F),P_1(F)]$$
is an isomorphism, which implies the first statement.
Similar argument works for algebraic characters. Finally,
$Q_{\b_0}$ is conjugate to $Q_{\ga_0}$ by the action of the longest
element of the Weyl group of $G_1$, so we obtain the same results
for $Q_{\b_0}$.
\ed

Note that in the case $G=E_k$ the lattice $\ker(X(T)\ra X(T_0))$ is
spanned by $\om$ and $\b_0$, while in the case $G=D_k$, $k\ge 5$,
it is spanned by $\om$, $\b_0$ and $(\a_1+\om)/2$.

Let us define the subvariety $C\subset V$ by setting $C=Le_{\b_0}$.

\begin{prop}\label{lagrprop} 
(i) $C$ is a Lagrangian subvariety invariant with
respect to the involution $*:V\ra V$. 

\noindent
(ii) $C\cup\{0\}$ coincides with the closure of the
graph of the birational map
$\grad(\frac{I}{x_0}):\La\ra\La'$,
where we identify $\La'$ with the dual of $\La$ by the
map $\la'\mapsto E(\cdot,\la')$ ($E$ is the symplectic form on
$V$).
\end{prop}

\Pf . (i) Since the symplectic form on $V$ is $L$-invariant,
it suffices to check that the tangent space to $C$ at $e_{\b_0}$
is a Lagrangian subspace of $V$. But this tangent space $T_{e_{\b_0}}C$ is
spanned by $e_{\b_0}$ and by the basis vectors of the form $e_{\b_0+\a}$ where
$\a$ belongs to the root system generated by $\De\setminus\{\b_0\}$.
Therefore, $T_{e_{\b_0}}C$ coincides with $\La$.
The involution $*$ corresponds to the action of the longest element of the
Weyl group of $L$, hence $C$ is invariant with respect to it.
                                     
\noindent
(ii) Since both varieties in question are irreducible and have dimension
$n+1$, it suffices to check that the graph $\Ga$ of
$\grad(\frac{I}{x_0})$ is $L$-invariant. The invariance of $\Ga$ with
respect to the action of the Borel subgroup $B\sub L$ follows
from Lemmas \ref{cubiclem1} and \ref{cubiclem2}. The fact that
$*(\Ga)=\Ga$ follows from the equality
$$\grad(I/x_0)\circ\grad(I/x_0)=\id,$$
where we identify $\La$ with $\La'$ using the involution $*:V\ra V$.
This equality follows easily from Lemma
\ref{prehomlem} below.
\ed

\begin{lem}\label{prehomlem} 
The action of the Levi subgroup $L_1\sub P_1$ 
preserves the subspaces $(e_{\b_0})\sub \La\sub V$.
The induced action of $L_1$ 
on variables $\La/(e_{\b_0})$ is
prehomogeneous (i.e., has an open orbit) and the cubic form $I$ is
a relative invariant for this action.
Let us identify $\La/(e_{\b_0})$ with $\La'/(e_{\ga_0})$
using the involution $*$ (these spaces are also naturally dual to each other).
Then one has
$$\grad I(\grad I(x))=I(x)x,$$
$$I(\grad I(x))=I(x)^2$$
for $x\in\La/(e_{\b_0})$.
\end{lem}

\Pf . In fact, these prehomogeneous spaces are well-known.
For $G=E_6$ this is (a subgroup of) $\GL_3\times\GL_3$ acting  
on $3\times 3$-matrices by left multiplication, and
our cubic form coincides with the determinant.
For $G=E_7$ we will have $L_1=\GL_6$ acting on $6\times 6$
skew-symmetric matrices and $I$ is the Pfaffian.
Finally, for $G=E_8$, we will have $L_1=E_6\times\G_m$
acting on the irreducible
$27$-dimensional representation of $E_6$ and $I$
coincides with the standard cubic form on this representation
(see \cite{Spr}).
The identities for $\grad I$ are well-known in the theory
of irreducible prehomogeneous spaces (see \cite{S}).

For $G=D_k$ we have $I=x_1Q(x_2,\ldots,x_n)$,
where $Q$ is a non-degenerate quadratic form, so these identities
are easy to check directly.
\ed 

Let us define the twisted action of $L$ on $V$ by the formula
$$g*v=\om(g)^{-1}gv.$$

\begin{lem}\label{transitivelem}
The twisted action of $L(F)$ on
$C(F)$ is transitive.
\end{lem}

\Pf . The similar statement over $\ov{F}$ is clear, so
it suffices to prove that the kernel (=the preimage of the trivial
class) of the natural map
$H^1(F,S)\ra H^1(F,L)$ is trivial, where $S=S_{\ga_0}$ is the
stabilizer subgroup
of $e_{\ga_0}\in C(F)$ (with respect to the twisted action).
Since $Q_{\ga_0}=S_{\ga_0}\times\om^{\vee}_{\b_0}(\G_m)$,
we have $H^1(F,S_{\ga_0})\simeq H^1(F,Q_{\ga_0})$.
On the other hand, $Q_{\ga_0}$ is a semi-direct product
of the normal subgroup $P_1$ and $\b_0^{\vee}(\G_m)$,
so $H^1(F,Q_{\ga_0})\simeq H^1(F,P_1)\simeq H^1(F,L_1)$.
Furthermore, since $L_1$ is a semi-direct product of
$G_0$ and an $F$-split torus, we have $H^1(F,L_1)\simeq H^1(F,G_0)$.
In the case $G=E_k$ the set $H^1(F,G_0)$ is trivial. On the other hand,
in the case $G=D_k$ the map
$H^1(F,G_0)\ra H^1(F,G_1)\simeq H^1(F,L)$ is an isomorphism.
\ed

Let us denote by $\La^0$ the open subset of $\La$ consisting of elements
with $x_0\neq 0$. By Proposition \ref{lagrprop}(ii) 
the projection $C\cap(\La^0\times\La')\ra\La^0$ is an isomorphism, so
we can identify $\La^0$ with an open subset in $C$.

Let $U\sub B\sub L$ be the unipotent radical of $B$, and let 
$\b_0^{\vee}(\G_m)\sub B$ be the one-dimensional subtorus corresponding
to $\b_0$. It is easy to see that 
$$\La^0=(\b_0^{\vee}(\G_m)U)*e_{\b_0}.$$
We are going to construct a subvariety 
$X\subset \b_0^{\vee}(\G_m)U$, such that
the natural map $X\ra\La^0: x\mapsto x*e_{\b_0}$
is an isomorphism.
Let $S_{\b_0}\sub L$ be the stabilizer
subgroup of $e_{\b_0}\in C(F)$ with respect to the twisted action of $L$:
$S_{\b_0}=\{g\in L: g*e_{\b_0}=e_{\b_0}\}$.
It is easy to see that the intersection
$S_{\b_0}\cap \b_0^{\vee}(\G_m)U$ is contained in $U$.
More precisely, we have
$$S_{\b_0}\cap U\b_0^{\vee}(\G_m)=S_{\b_0}\cap U=Q_{\b_0}\cap U=P_1^-\cap U,$$
where $Q_{\b_0}$ be the stabilizer of the line spanned by $e_{\b_0}$,
$P_1^-\sub G_1$ is the opposite subgroup to the maximal parabolic
$P_1$ (see Lemma \ref{parablem}).
Therefore, it suffices to find a subvariety $X_0\sub U$ such that
the projection $X_0\ra U/P_1^-\cap U$ is an isomorphism (then
we can set $X=\b_0^{\vee}(\G_m)X_0$). But 
$$P_1^-\cap U=L_1\cap U=G_0\cap U,$$
where $G_0\sub G_1$ is the semisimple subgroup corresponding to simple roots
which are not attached to $\b_0$. Let us denote by $\Phi_1\sub\Phi$ the
root system of $G_1$ and let $\Sigma_1\sub\Phi_1$ be the set of
roots in $\Phi_1$ which have a positive coefficient with some 
simple root attached to $\b_0$. Then we can set
$$X_0=\prod_{\a\in\Sigma_1}e_{\a}(F)\sub U.$$
Now the subvariety $X=\b_0^{\vee}(\G_m)X_0\sub L$
projects isomorphically to $\La^0$ by the twisted action of $L$ on 
$e_{\b_0}$. Thus, there is a unique map
\begin{equation}\label{sigmamap}
\sigma:\La^0\ra X
\end{equation}
such that 
$$\sigma(v)*e_{\b_0}=v$$
for all $v\in\La^0$.

\subsection{Minimal representation and the
$L(F)$-equivariant line bundle on $C(F)$}
\label{conesec}

Recall that we denote by $R$ the space of smooth vectors in
the minimal representation of $G(F)$.
Let us denote by $R'\sub R$ the subspace consisting of locally
constant functions with compact support in $F^*\times F^{n+1}$
(see Lemma \ref{supportlem}). It is clear that the subgroup
$Z(H)(F)=e_{\om}(F)\sub G(F)$ acts trivially on the quotient $R/R'$.
Therefore, the space $R/R'$
is a representation of the semi-direct product of 
$V(F)=H(F)/Z(H)(F)$ and $L(F)$. In this section
we construct an $L(F)\ltimes V(F)$-equivariant map from $R/R'$ to the
space of locally constant sections of certain $L(F)$-equivariant
line bundle over $C(F)\sub V(F)$.
Later we will show that in fact
$R/R'$ coincides with coinvariants of $Z(H)(F)$ 
on $R$ and that the constructed map is an embedding
(see Proposition \ref{coinvprop}). 

For every non-zero vector $v\in V$
let us consider the linear map $l_v:V(F)\ra F$ given by
$l_v(x)=E(v,x)$, 
and extend it to a homomorphism $l_v:H(F)\ra F$,
trivial on $Z(H)(F)$. 
To understand the spectrum of the action of $V(F)$ on $R/R'$ we 
have to compute the (twisted) spaces of coinvariants
$$R(v):=R_{H,l_v}=R/\{hf-\psi(l_v(h))f, h\in H(F), f\in R\}$$
for all $v\neq 0$.

Let $f(y,x_0,x)$ be an element of $R$. 
We claim that for fixed $x_0\in F^*$ and $x\in F^n$, the expression
$\psi(\frac{I(x)}{x_0y})f(y,x_0,x)$ stabilizes as
$y$ tends to $0$. Indeed, applying the operator $A$ we see that our claim
is equivalent to the statement that for $y\neq 0$ the expression
$f(y,x_0,x)$ stabilizes as $x_0$ tends to $0$.
But this follows from invariance
of $f$ with respect to $e_{\b_0}(\pi^m\OO)$ for sufficiently
large $m$.

Recall that we denoted by $\La^0$ the open subset of 
$\La$ consisting of elements with $x_0\neq 0$. 
Let us define a map $R\ra \C(\La^0(F)):f\mapsto \ov{f}$ by
setting 
\begin{equation}\label{fbareq}
\ov{f}(x_0,x)=\lim_{y\ra 0} \psi(\frac{I(x)}{x_0y})f(y,x_0,x).
\end{equation}

\begin{lem}\label{simpleidlem}
One has the following identities
$$\ov{e_{\om}(a)f}=\ov{f},$$
$$\ov{e_{\ga_i}(a)f}=\psi(\eps_i a x_i)\ov{f},$$
$$\ov{e_{\b_i}(a)f}=\psi(-a\frac{\partial(I/x_0)}{\partial x_i})\ov{f},$$
where $a\in F$. Also, for $t$ in the maximal torus we have
$$\ov{tf}(x_0,\ldots,x_n)=|\de(t)|^{\frac{1}{2}}
\ov{f}(\ga_0(t)x_0,\ldots,\ga_n(t)x_n),$$
where the character $\de$ was defined in Lemma \ref{rootlem2}.
\end{lem}

\Pf . The first three identities follow easily from the explicit formulas
for the action of the relevant operators on $R$. The last identity 
follows from the formula for the action of $t$ on $R$ and from
Lemma \ref{cubiclem2}.
\ed

Consider the linear functional $f^*_0$ on $R$ defined by
$$f^*_0(f)=\ov{f}(\b_0)=\ov{f}(1,0)=\lim_{y\ra 0}f(y,1,0).$$
We have
$$f^*_0(e_{\om}(t)f)=f^*_0(e_{\b_i}(t)f)=f^*_0(e_{\ga_j}(t)f)=f^*_0(f),$$
where $i=0,\ldots,n$, $j=1,\ldots,n$, while 
$$f^*_0(e_{\ga_0}(t)f)=\psi(\eps_0t)f^*_0(f).$$
In other words, $f^*_0$ descends to a functional on the 
space of coinvariants $R(e_{\b_0})$.
The following Proposition 
implies that $f^*_0$ is characterized by this property up
to a scalar. In the case of exceptional groups it is equivalent to Lemma 6.2 of
\cite{MS} (in the case $G=E_8$ Magaard and Savin impose the additional assumption
that the residual characteristic of $F$ differs from $2$). 

\begin{prop}\label{onedimlem}
(i) The functional $f^*_0$ is not equal to zero.

\noindent
(ii) For every $v\in C(F)$
the space $R(v)$ is one-dimensional. 

\noindent (iii) For $v\neq 0$ such that $v\not\in C(F)$, one has $R(v)=0$. 
\end{prop}

\Pf . (i) Let $f\in R'$ be
the characteristic function of a small compact neighborhood of
$(1,0,\ldots,0)\in F^*\times F^{n+1}$. Then $f^*_0(Af)\neq 0$, hence
$f^*_0\neq 0$.

\noindent
(ii) For every $g\in L(F)$ we have $l_v(Ad(g^{-1})h)=l_{\om(g)^{-1}gv}(h)$.
Hence, the action of $g$ induces an isomorphism
$g:R(v)\wt{\ra} R(\om(g)^{-1}gv)=R(g*v)$.
Since by Lemma \ref{transitivelem} the twisted action of $L(V)$ on $C(F)$
is transitive, it suffices to prove that  
the space $R(e_{\b_0})$ is one-dimensional.
From part (i) we know that $f^*_0$ is a non-zero element of $R(e_{\b_0})^*$. 
Now let $f^*$ be a functional on $R$ that factors through $R(e_{\b_0})$.
Consider the functional 
$$f^*_1(f)=f^*(Af)$$
where $A$ is the operator \ref{operatorA}.
Since the simple reflection $s_{\b_0}$ preserves $\ga_j$ with $j\ge 1$
invariant, we obtain
$$f^*_1(e_{\ga_j}(t)f)=f^*_1(f)$$
for all $j=1,\ldots,n$. On the other hand, since $s_{\b_0}$ switches
$\om$ and $\ga_0$, we have the relations
$$e_{\om}(t)=s_{\b_0}e_{\ga_0}(\eps_0t)s_{\b_0}^{-1}$$ 
$$e_{\ga_0}(t)=s_{\b_0}e_{\om}(-\eps_0t)s_{\b_0}^{-1}.$$
Hence,
$$f^*_1(e_{\ga_0}(t)f)=f^*_1(f),$$ 
$$f^*_1(e_{\om}(t)f)=\psi(-t)f^*_1(f).$$
We claim that the latter condition implies that $f^*_1(f)$ depends only on
the restriction $f'(x_0,\ldots,x_n)=f(-1,x_0,\ldots,x_n)$ of $f$ to $y=-1$.
Indeed assume that $f(-1,x_0,\ldots,x_n)\equiv 0$. We want to prove that
$f^*_1(f)=0$. Let us write
$f=f_1+f_2$, where $f_1=f\de_{|y|<\eps}$,
$f_2=f\de_{|y|\ge\eps}$, $\eps>0$, and  
$\de_C$ denotes the characteristic function of a set $C$. By Lemma
\ref{supportlem}, we have $f_2\in R'$.
Hence, $f_1\in R$. Let us pick $t_0\in F^*$ such that $\psi(-t_0)\neq 1$.
For sufficiently small $\eps$ we have
$\psi(yt_0)f_1=f_1$, hence
$$f_1=\frac{\psi(yt_0)-\psi(-t_0)}{1-\psi(-t_0)}\cdot f_1=
(e_{\om}(t_0)-\psi(-t_0))\frac{f_1}{1-\psi(-t_0)}.$$
This implies that $f^*_1(f_1)=0$. 
On the other hand, since $f_2(-1,x_0,\ldots,x_n)\equiv 0$,
it follows that $f_2$ lies in $\sum_{t\in F^*}(e_{\om}(t)-\psi(-t))(R')$, so
$f^*_1(f_2)=0$. Therefore, $f^*_1(f)=0$ which proves our claim.
Note that $f'=f|_{-1\times F^{n+1}}$ can be an arbitrary locally constant
function with compact support in $(x_0,\ldots,x_n)$ (as follows from
Lemma \ref{supportlem}(c)). Since 
of $e_{\ga_i}(t)f=\psi(\eps_itx_i)f$ we obtain that
$$f^*_1(\psi(tx_i)f')=f^*_1(f')$$
for all $i\ge 0$, which implies that
$f^*_1$ is proportional to $f\mapsto f(-1,0,\ldots,0)$. 
It follows that $f^*$ is proportional to $f^*_0$.

\noindent (ii) Using the twisted action of $L(V)$ as in the proof of
part (i) we can assume that $x_0(v)\neq 0$.
Let $f^*$ be a functional on $R$ that factors through $R(v)$.
The same argument as in part (i) shows that $f^*$ is proportional
to the functional $f\mapsto \ov{f}(x_0(v),\ldots,x_n(v))$.
Now imposing the conditions 
$f^*(e_{\b_i}(t)f)=\psi(t E(v,e_{\b_i}))f^*(f)$ for $i=0,\ldots,n$ 
and using Lemma \ref{simpleidlem} and Proposition \ref{lagrprop} 
we immediately see that $v\not\in C(V)$ should imply
that $f^*=0$.
\ed

Recall that $S_{\b_0}\sub L$ is the stabilizer
subgroup of $e_{\b_0}\in C(F)$ with respect to the twisted action of $L$.

\begin{cor} The functional $f^*_0:R\ra\C$ is an eigenvector for 
the action of $S_{\b_0}(F)$. 
\end{cor}

Let $\chi:S_{\b_0}(F)\ra\C^*$ be the character defined by the condition
\begin{equation}\label{eigencharacter}
f_0^*(g^{-1}f)=\chi(g)f_0^*(f)
\end{equation}
for all $g\in S_{\b_0}(F)$, $f\in R$.
Using this character we can define an 
an $L(F)$-equivariant complex line bundle over $C(F)$,
such that its fiber at the point $v\in C(F)$ is canonically
identified with $v$. 

\begin{defi} We define $\LL$ to
be the $L(F)$-equivariant line bundle over $C(F)$
corresponding to the character $\chi$ of $S_{\b_0}(F)$:
$\LL$ is the quotient $\C\times L(F)/S_{\b_0}(F)$,
where the (right) action of $S_{\b_0}(F)$ on $\C\times L(F)$ is given by
$g_0(z,g)=(\chi(g_0)z,gg_0)$. The projection $\LL\ra C(F)$ is given
by $(z,g)\mapsto g*e_{\b_0}$.
\end{defi}

We have a natural identification 
$R(v)\wt{\ra}\LL|_v$: an element $f\in R(v)$ corresponds
to an orbit of $(f^*_0(g^{-1}f),g)$, where $g\in L(F)$ is such that
$g*e_{\b_0}=v$. The relation (\ref{eigencharacter}) garantees that this
map is well-defined. Under this isomorphism the
action of $L(F)$ on $\LL$ is given by the natural operators
$g:R(v)\ra R(g*v)$, where $g\in L(F)$, $v\in C(F)$
(recall that
the action of $L(F)$ on $\LL$ is compatible with its twisted action on
$C(F)$).

Let us denote by $\pi_v:R\ra R(v)$
the natural projection. Using the above identification,
we get an $L(F)$-equivariant map
$$\res:R\ra \Ga(C(F),\LL):f\mapsto (v\mapsto \pi_v(f)),$$
where $\Ga(C(F),\LL)$ is the space of all global sections of $\LL$.
In fact, this map is $P(F)$-equivariant, where $H(F)$ acts on
$\Ga(C(F),\LL)$ through the quotient $V(F)$ and the action of
$v\in V(F)$ on $s\in\Ga(C(F),\LL)$ is given by
\begin{equation}\label{Vaction}
v\cdot s(x)=\psi(E(x,v))s(x).
\end{equation}
Indeed, this follows immediately from the properties of the functional $f^*_0$
and from the identity $E(g*x,gy)=E(x,y)$ for $x,y\in V$, $g\in L$.

Below we will compute the character $\chi$
explicitly (see Lemma \ref{eigenlem}). In particular, we will
see that it is locally constant. Hence, the transition functions for
the line bundle $\LL$ are also locally constant. Therefore,
it makes sense to talk about locally constant sections of $\LL$.
Since $R$ consists of smooth vectors, the image of the map $\res$ is 
contained in the subspace of smooth vectors in $\Ga(C(F),\LL)$ considered
as a representation of $L(F)$. In particular, the sections $\res(f)$ for
$f\in R$ are locally constant. Also, from Lemma \ref{supportlem},
we obtain that $\res(f)$ has bounded support.

In the next lemma we compute the character $\chi$ of $S_{\b_0}$. 
Recall that in Lemma \ref{rootlem2} we defined a character $\de:T\ra\G_m$.
It is easy to check that $\de|_{T_0}\equiv 1$, where $T_0\sub T$
is the maximal torus of $G_0$. Therefore,
using Lemma \ref{charlem} we can extend $\de$ to a character of $Q_{\b_0}$.

\begin{lem}\label{eigenlem} 
One has
\begin{equation}\label{eigeneq}
\chi(g)=|\de(g)|^{-\frac{1}{2}}
\end{equation}
where $g\in S_{\b_0}(F)$.
\end{lem}

\Pf . By Lemma \ref{charlem}, the root $\ga_0$ (considered as a character
of the maximal torus) extends to the character $\ga_0:Q_{\b_0}\ra\G_m$, 
such that
$S_{\b_0}=\ker\ga_0$. Furthermore, since $\lan\ga_0,\b_0\ran=-1$, the group
$Q_{\b_0}$ is the direct product of $S_{\b_0}$ and 
$Z(L)=\om^{\vee}_{\b_0}(\G_m)$. 

Set $T'=T\cap S_{\b_0}$, so that
$T'$ is the kernel of the character $\ga_0:T\ra\G_m$. 
We claim that every character of $S_{\b_0}(F)$
is determined by its restriction to $T'(F)$. Indeed, composing
a character $\mu:S_{\b_0}(F)\ra\C^*$
with the projector $Q_{\b_0}\ra S_{\b_0}$ along $Z(L)$,
we get a character $\wt{\mu}$ of $Q_{\b_0}(F)$. Assuming
that $\mu|_{T'(F)}\equiv 1$ we get $\wt{\mu}|_{T'(F)}\equiv 1$.
On the other hand, $\wt{\mu}$ is trivial on $Z(L)(F)$. Since
$T$ is the direct product of $T'$ and $Z(L)$, we obtain that
$\wt{\mu}|_{T(F)}\equiv 1$. Now Lemma \ref{charlem} implies that
$\wt{\mu}\equiv 1$, hence $\mu\equiv 1$. This proves our claim.

Therefore, it suffices to check the equality (\ref{eigeneq}) for
$g\in T'(F)$. But in this case it follows 
from the explicit formula for the action of $T(F)$ on $R$.
\ed

\begin{rem}
In the case $G=E_k$ the character $\de$ of $S_{\b_0}$
is the restriction of the character $2(s+1)\om$ of $L$ 
where $s=1,2,4$ for $k=6,7,8$ (see Lemma \ref{rootlem2}). 
It follows that in this
case the line bundle $\LL$ has a canonical trivialization (however,
the $L$-equivariant structure is twisted). 
In the case $G=D_4$ we have a similar situation, since $\de$ coincides
with the restriction of $2\om$.
However, in the case $G=D_k$ with $k\ge 5$, the character $\de$ of
$S_{\b_0}$ does not extend to a character of $L$, so we cannot
trivialize $\LL$.
\end{rem}

\subsection{Trivialization over an open part}\label{opensec}

In this section we consider a trivialization of $\LL$ over an
open subset $\La^0(F)\sub C(F)$,
under which the natural map 
$$\res:R\ra\Ga(C(F),\LL)\ra\Ga(\La^0(F),\LL)$$ 
gets identified with the map $f\mapsto\ov{f}$ given by (\ref{fbareq}). 
Namely, let $\sigma:\La^0\ra X\sub\b_0^{\vee}(\G_m)U$ be the
map (\ref{sigmamap}). Then for every $v\in\La^0$ the action of 
$\sigma(v)$ defines an isomorphism $R(e_{\b_0})\wt{\ra} R(v)$
(recall that $\sigma(v)*e_{\b_0}=v$).
Since $f^*_0$ identifies $R(e_{\b_0})$ with $\C$,
we get a trivialization of $\LL$ over $\La^0$.
We claim that under this trivialization, the section $\res{f}$
corresponds to the function $\ov{f}$ on $\La^0$. Here is the
precise statement.

\begin{prop}\label{trivprop} 
For every $f\in R$ and $v\in\La^0$ one has
$$f^*_0(\sigma(v)^{-1}f)=\ov{f}(v).$$
\end{prop}

\Pf . Let us write $g=\sigma(v)=tu$, where $t\in\b_0^{\vee}(F^*)$, $u\in U$.
By the definition we have $g*e_{\b_0}=\om(t)^{-1}ge_{\b_0}=v$. 
On the other hand,
$$g^{-1}f(y,e_{\b_0})=|\de(t)|^{-\frac{1}{2}}\cdot
\psi(\frac{E(B_{u^{-1}}A_u e_{\b_0},A_u e_{\b_0})}{2y})
f(\om(t)^{-1}y,\om(t)^{-1}tA_u e_{\b_0}).$$
Note that
$$\lan\b_0,\de\ran=\lan\b_0,\om\ran+\sum_{i=0}^n\lan\b_0,\ga_0\ran=0$$
since $\lan\b_0,\om\ran=1$, $\lan\b_0,\ga_0\ran=-1$ and $\lan\b_0,\ga_i\ran
=0$ for $i\ge 1$. Therefore, $\de(t)=1$. Also, by Lemma
\ref{cubiclem1}, we have
$$\frac{E(B_{u^{-1}}A_ue_{\b_0},A_ue_{\b_0})}{2}=
\frac{-I(e_{\b_0})+I(A_ue_{\b_0})}{x_0(A_ue_{\b_0})}=
I(A_ue_{\b_0}).$$
Substituting this into the above expression we get
$$f^*_0(g^{-1}f)=\lim_{y\ra 0}\psi(\frac{I(A_ue_{\b_0})}{y})
f(\om(t)^{-1}y,\om(t)^{-1}tA_u e_{\b_0})=
\lim_{y\ra 0}\psi(\frac{I(A_ue_{\b_0})}{\om(t)y})
f(y,\om(t)^{-1}tA_u e_{\b_0}).$$
Using Lemma \ref{cubiclem2}, it is easy to see that this
limit is equal to $\ov{f}(\om(t)^{-1}ge_{\b_0})=\ov{f}(v)$.
\ed

\begin{cor}\label{trivlinecor} 
Assume that $G=E_k$ (resp. $G=D_4$). 
Then for every $f\in R$ the function 
$|x_0|^{s+1}\ov{f}(x_0,\ldots,x_n)$ (resp.
$|x_0|\ov{f}(x_0,\ldots,x_n)$) 
extends to a locally constant function
on $C(F)$ and the action of $L(F)$ on $f$ is compatible with the 
action of $L(F)$ on $|\om|^{s+1}\otimes\C(C(F))$
(resp. $|\om|\otimes\C(C(F))$), where $|\om|$ is the line bundle
of volume forms on $C(F)$. 
\end{cor}

\subsection{$L(\OO)$-invariant vectors}\label{invvecsec}

In this section we derive a general form of $\ov{f}$,
where $f\in R$ is an $L(\OO)$-invariant vector.
The derivation is based on the study of $L(\OO)$-invariant sections of 
the line bundle $\LL$ over $C(F)$.

Recall that by Lemma \ref{transitivelem}, 
$C(F)$ can be identified with $L(F)/S_{\b_0}(F)$.
Therefore, every global section $\Phi$ of $\LL$
can be described by the function
\begin{equation}\label{PhiPsi}
\Psi(g)=f^*_0(g^{-1}\Phi(g*e_{\b_0}))
\end{equation}
on $L(F)$, where $f^*_0$ is the functional
on $\LL|_0=R(e_{\b_0})$ considered in section \ref{conesec}.
Now we are going to use the Iwasawa decomposition
$$L(F)=L(\OO)Q_{\b_0}(F)=L(\OO)\om^{\vee}_{\b_0}(F^*)S_{\b_0}(F).$$

\begin{lem}\label{invseclem} 
Let $\Phi$ be an $L(\OO)$-invariant global section of $\LL$.
Then there exists a function $\Psi_{\Z}$ on $\Z$ such that
$$\Psi(k\om^{\vee}_{\b_0}(a)g_0)=|\de(g_0)|^{-\frac{1}{2}}\Psi_{\Z}(v(a)),$$
where $\Psi$ is the function on $L(F)$ given by (\ref{PhiPsi}),
$k\in L(\OO)$, $a\in F^*$, $g_0\in S_{\b_0}(F)$.
\end{lem}

\Pf . Lemma \ref{eigenlem} implies that 
$$
\Psi(gg_0)=|\de(g_0)|^{-\frac{1}{2}}\Psi(g)
$$
for $g_0\in S_{\b_0}(F)$. On the other hand,
the condition of $L(\OO)$-invariance of
$\Phi$ is equivalent to the condition $\Phi(kg)=\Phi(g)$ for
$k\in L(\OO)$. Therefore,
$$\Psi(k\om^{\vee}_{\b_0}(a)g_0)=|\de(g_0)|^{-\frac{1}{2}}
\Psi(\om^{\vee}_{\b_0}(a)).$$
It remains to notice that $\Psi(\om^{\vee}_{\b_0}(a))$
depends only on $|a|$ (by $L(\OO)$-invariance of $\Psi$).
\ed

\begin{prop}\label{invariantprop} 
Let $f\in R$ be an $L(\OO)$-invariant vector. Then
there exists a function $f_{\Z}$ on $\Z$ such that
$$\ov{f}(x)=c(x)f_{\Z}(v(x)),$$
where $x\in C(F)$, the function  
$c(x)$ is given by
$$c(x)=\cases |x_0|^{-s-1}, & G=E_k,\\
|x_0|^{-1}\max(1,|\frac{x_1}{x_0}|)^{k-4}, & G=D_k,\endcases$$
where $s=1,2,4$ for $G=E_6,E_7,E_8$.
\end{prop}

\Pf . Let $\Phi=\res(f)$ be the $L(\OO)$-invariant global section of $\LL$
obtained from $f$ by the construction of section \ref{conesec}.
According to Proposition \ref{trivprop}, the function $\ov{f}$
coincides with the function on $\La^0$
obtained from $\Phi$ by the trivialization of $\LL|_{\La^0}$. 
More precisely, for $x=(x_0,\ldots,x_n)\in\La^0$ we have
\begin{equation}\label{ovfPsisigma}
\ov{f}(x)=f^*_0(\sigma(x)^{-1}\Phi(x))=\Psi(\sigma(x)),
\end{equation}
where $\Psi$ is given by (\ref{PhiPsi}).
Let us denote $g=\sigma(x)$, so that $x=g*e_{\b_0}$. 
Applying the Iwasawa decomposition we can write
$g=kt'g_0$, where $k\in L(\OO)$, $t'=\om^{\vee}_{\b_0}(a)$ for some 
$a\in F^*$, $g_0\in S_{\b_0}$.
Moreover, we can assume that $k\in G_1(\OO)$. Indeed, we can write
$k=k'\b_0^{\vee}(c)$, where $k'\in G_1(\OO)$, $c\in F^*$.
But 
$$\b_0^{\vee}(c)=\om^{\vee}_{\b_0}(c^{-1})\cdot 
(\om^{\vee}_{\b_0}+\b^{\vee})(c),$$
and $(\om^{\vee}_{\b_0}+\b^{\vee})(\G_m)\sub S_{\b_0}$.
Applying Lemma \ref{invseclem}, we can rewrite (\ref{ovfPsisigma}) as
$$\ov{f}(x)=|\de(g_0)|^{-\frac{1}{2}}\Psi_{\Z}(v(a)).$$
It remains to find $|a|$ and $|\de(g_0)|$ in terms of the vector $x$.
Since the action of $L(\OO)$ on $V$ preserves the norm,
we have $|x|=|t'*e_{\b_0}|=|a|^{-1}$, so $|a|=|x|^{-1}$. Finding $\de(g_0)$
requires some more work.

Recall that by the definition, $g=\sigma(x)$ has form
$g=ut$, where $t=\b_0^{\vee}(b)$ for some $b\in F^*$, 
$u\in\prod_{\a\in\Sigma_1}e_{\a}(F)$. Therefore, we have
$$x_0=x_0(ut*e_{\b_0})=\ga_0^{-1}(t)=b.$$
On the other hand, applying the character $\om:L(F)\ra F^*$ to
the equality $g=ut=kt'g_0$, we get
$$\om(g_0)=\frac{b}{a^2}=x_0\cdot|x|^2.$$
In the case $G=E_k$ this finishes the proof: since $\de=2(s+1)\om$,
the function $x\mapsto |\de(g_0)|^{-\frac{1}{2}}$ has the required form.
In the case $G=D_k$ we need in addition to compute $|\a_1(g_0)|$,
where $\a_1$ is considered as a character of $Q_{\b_0}$.
Since $\om(t'g_0)=b$, we can write $t'g_0=g't$,
where $g'\in G_1\cap Q_{\b_0}$. Then we have the equality $u=kg'$ in $G_1$.
Let us denote elements of $G_1=\SL_2\times\SO(2k-4)$ as pairs $h=(h_1,h_2)$,
where $h_1\in\SL_2$, $h_2\in\SO(2k-4)$. Then we have
$u_1=k_1g'_1$, $u_2=k_2g'_2$. Recall that $V=V_1\otimes V_2$, where
$V_1$ and $V_2$ are standard representations of $\SL_2$ and $\SO(2k-4)$.
Let $(e_{-1},e_1)$ be the standard basis of $V_1$, and let $\xi$ be
the lowest weight vector in the standard basis of $V_2$. 
Then $e_{\b_0}=e_{-1}\otimes \xi$, $e_{\b_1}=e_1\otimes \xi$.
Together with the equality
$$x_0^{-1} x=ue_{\b_0}=u_1e_{-1}\otimes u_2\xi$$
this implies that 
$$|u_1e_{-1}|=\max(1,|\frac{x_1}{x_0}|).$$
On the other hand, $u_1=k_1g'_1$, so
$$|u_1e_{-1}|=|g'_1e_{-1}|=|\frac{\om-\a_1}{2}(g'_1)|.$$
Therefore,
$$|\a_1(g')|=|\a_1(g'_1)|=\max(1,|\frac{x_1}{x_0}|)^{-2}.$$
It remains to apply $\a_1$ to the equality $t'g_0=g't$ to deduce that
$$|\a_1(g_0)|=|b^{-1}\a_1(g')|=|x_0|^{-1}\cdot
\max(1,|\frac{x_1}{x_0}|)^{-2}.$$
Hence,
$$|\delta(g_0)|^{-\frac{1}{2}}=|\om(g_0)|^{\frac{-k+2}{2}}
|\a_1(g_0)|^{\frac{-k+4}{2}}=|x_0|^{-1}\max(1,|\frac{x_1}{x_0}|)^{k-4}\cdot
|x|^{-k+2}.
$$ 
\ed

%
%
%
%
%
%

\subsection{Twisted case}\label{twistedconesec}

This section will be used only for global considerations in part II.

Let $R_{\psi}$ be a twisted realization of the minimal representation 
associated with a non-trivial character $\psi:F\ra\C^*$ (see \ref{twistedsec}).
Then we still have a unique (up to constant)
$P(F)$-equivariant map 
\begin{equation}\label{psimap}
R_{\psi}\ra \Ga(C(F),\LL)_{\psi}
\end{equation}
where $H(K_v)$ acts on sections of $\LL_v$ by the formula (\ref{Vaction}).

Let $\psi(t)=\psi_0(at)$, where $a\in F^*$ and $\psi_0$ is of order $0$
and let $R$ be the standard realization of the
minimal representation associated with $\psi_0$.
We have a $G(F)$-equivalence
$R\ra R_{\psi}: f\mapsto f(ay,ax_0,\ldots,ax_n)$.
On the other hand, the action of $\om_{\b_0}^{\vee}(a^{-1})\in Z(L)(F)$ gives an
equivalence of $P(F)$-representations
$\Ga(C(F),\LL)_{\psi_0}\ra\Ga(C(F),\LL)_{\psi}$.
We can normalize the map (\ref{psimap}) above by requiring the following diagram
to be commutative (the normalized map is denoted by $\res_{\psi}$):
\begin{equation}
\begin{array}{ccc}
R & \lrar{} & R_{\psi}\\
\ldar{\res} & &\ldar{\res_{\psi}}\\
\Ga(C(F),\LL)_{\psi_0} & \lrar{} & \Ga(C(F),\LL)_{\psi}
\end{array}
\end{equation}
Equivalently, for every $f\in R_v$ we have
$$\res_{\psi}(f)(a^{-1}e_{\b_0})=\lim_{y\to 0}f(y,a^{-1},0)\cdot\om_{\b_0}^{\vee}(a)(1),$$
where the map
$\om_{\b_0}^{\vee}(a):\LL_v|_{e_{\b_0}}\ra\LL_v|_{a^{-1}e_{\b_0}}$ is
given by $L(F)$-equivariant structure on $\LL$, the element $1\in\LL_{e_{\b_0}}$
corresponds to the canonical trivialization of $\LL_{e_{\b_0}}$.

Let us compute $\res_{\psi}(f)(e_{\b_0})$ for $f\in R_v$.
For every $s\in\Ga(C(F),\LL)$ we have
$$\om_{\b_0}^{\vee}(a)(s(e_{\b_0}))=(\om_{\b_0}^{\vee}(a)s)(a^{-1}e_{\b_0}).$$
Applying this to $s=\res_{\psi}(f)$ and using the $L(F)$-equivariance of $\res_{\psi}$
we get
\begin{align*}
&\om_{\b_0}^{\vee}(a)(\res_{\psi}(f)(e_{\b_0}))=\res_{\psi}(\om_{\b_0}^{\vee}(a)f)(a^{-1}e_{\b_0})=\\
&\lim_{y\to 0}(\om_{\b_0}^{\vee}(a)f)(y,a^{-1},0)\cdot\om_{\b_0}^{\vee}(a)(1)=
|a|^{\frac{n+3}{2}}\lim_{y\to 0}f(a^2y,1,0)\cdot \om_{\b_0}^{\vee}(a)(1).
\end{align*}
Hence,
\begin{equation}\label{resnormal}
\res_{\psi}(f)(e_{\b_0})=|a|^{\frac{n+3}{2}}\lim_{y\to 0} f(y,1,0).
\end{equation}

\section{Action of the affine Hecke algebra}\label{heckesec}


\subsection{Lusztig's representation}
\label{Lusztigsec}

It is known (see \cite{Savin}) that the representation of $G(F)$ on the space $R$ is 
equivalent to the spherical representation $J(\chi_{sr})$ with
the Satake parameter $\chi_{sr}$ corresponding to the subregular unipotent
orbit in the Langlands dual group $\sideset{^L}{}G$.

Let $I\sub G(F)$ be the Iwahori subgroup.
Then $R^I$ is a representation of the affine Hecke algebra 
$\wt{H}$ of $I$-biinvariant compactly supported functions on $G(F)$.
It is equivalent to the representation of $\wt{H}$ on $J(\chi_{sr})^I$
which was described by Lusztig in \cite{L1}.
Below we are going to recall some details of this description.

Let $\wt{\De}=\De\cup\{-\om\}$ be the set of vertices 
of the affine Dynkin diagram of $G$ (thought of as roots of $G$).
Recall that $\wt{H}$ is generated by elements $(T_{\a}, \a\in\wt{\De})$ 
with the defining relations
$$T_{\a}T_{\b}=T_{\b}T_{\a}\ \text{if}\ \lan\a,\b\ran=0,$$
$$T_{\a}T_{\b}T_{\a}=T_{\b}T_{\a}T_{\b}\ \text{if}\ \lan\a,\b\ran=-1,$$
$$(T_{\a}-q)(T_{\a}+1)=0.$$
Lusztig has shown in \cite{L1} that
the representation of $\wt{H}$ on $J(\chi_{sr})^I$ is equivalent
to the representation of $\wt{H}$ on the space 
$$\wt{E}=\oplus_{\a\in\wt{\De}}\C e_{\a}$$
given by the formulas
$$T_{\a}e_{\b}=\cases -e_{\b} &\ \text{if}\ \a=\b,\\ 
qe_{\b}+q^{\frac{1}{2}}e_{\a} &\ \text{if}\ \lan\a,\b\ran=-1,\\
qe_{\b} &\ \text{if}\ \lan\a,\b\ran=0.\endcases$$

The important part of the work \cite{L1} is the study of the action of
Bernstein's elements in $\wt{H}$ on $\wt{E}$. Recall that for every simple
root $\a\in\De$ one defines an element $T^{\a}\in\wt{H}$ as the characteristic
function of the double coset $I\om^{\vee}_{\a}(\pi)I\sub G(F)$, where
$\pi\in\OO$ is a uniformizer. It is known 
that these elements commute between themselves. One can also
describe some commutation relations
between the elements $T_{\a}$ and $T^{\b}$. We will only need the following
simple relation:
\begin{equation}\label{comrel}
T_{\a}T^{\b}=T^{\b}T_{\a}\ \text{if}\ \a(\om^{\vee}_{\b})=0.
\end{equation}

Lusztig describes an explicit basis in $\wt{E}$ with respect to
which the action of the elements $T^{\a}$ is almost diagonal. More
precisely, let $\a_b$ be the simple root corresponding to the branch point
of the Dynkin diagram of $G$ (recall that $G$ is assumed to be
either of type $D$ or of type $E$). Then according to Theorem 4.7 of \cite{L1},
the action of $T^{\a}$ on $\wt{E}$ is semisimple for all $\a\neq\a_b$.
To describe the eigenvalues of $T^{\a}$ computed in {\it loc. cit.}
we need to introduce some notation. For every $\a,\b\in\De$ 
one defines an integer $\la(\a,\b)$ using the following procedure.
First we define the homomorphism $\la_b:P\ra\frac{1}{2}\Z$, where
$P$ is the coweight lattice by its values on simple coroots:
$$\la_b(\a^{\vee})=\cases 1, &\ \a\neq\a_b,\\ 0, &\ \a=\a_b.\endcases$$
Next we set
$$\la(\a_b,\b)=\la_b(\om^{\vee}_{\b}).$$
Finally, for every $\a\neq\a_b$ we define
$$\la(\a,\b)=\cases \la(\a_b,\b)-d(\a_b,\b), &\ \text{if}\ \b\
\text{belongs to the geodesic from}\ \a_b\ \text{to}\ \a,\\
\la(\a_b,\b), &\ \text{otherwise},\endcases$$ 
where $d(\a_b,\b)$ is the distance between the vertices $\a_b$ and $\b$
on the Dynkin diagram.
Let $\rho$ denotes the half-sum of all positive roots of $G$.
Then the part of Theorem 4.7 of \cite{L1} that we need can
be summarized as follows. 

\begin{thm}\label{Lusztigthm} 
For every $\b\in\De$ the eigenvalues of
$q^{-\rho(\om^{\vee}_{\b})}T^{\b}$ acting 
on $\wt{E}$ are $(q^{\la(\a,\b)},\a\in\De)$, where
$q^{\la(\a_b,\b)}$ is taken twice.  
For $\b\neq\a_b$ the action of $T^{\b}$ on $\wt{E}$ is semisimple.
\end{thm}

\begin{cor}\label{heckecor} 
For $G=E_k$ 
the eigenvalues of $T^{\b_0}$ acting on $\wt{E}$ are
$q^{2(n+1)-s}$ with multiplicity $2$ and $q^{2(n+1)-2s}$ with
multiplicity $k-1$ (where $n=6s+3$). For $G=D_k$ these eigenvalues are
$q^{2(n+1)}$ and $q^{2(n+1)+4-k}$ (where $n=2k-5$). For $G\neq D_4$
the action of $T^{\b_0}$ on $\wt{E}$ is semisimple.
\end{cor}

\Pf . This follows from the above Theorem and from the equality
$\rho(\om^{\vee}_{\b_0})=n+2$.
\ed

\subsection{The form of the spherical vector}

In this section we are going to combine the above information on $R^I$
with Proposition \ref{invariantprop} 
to derive the specific form of the spherical
vector in the standard representation (with two unknown constants).

First of all, we claim that the spherical vector $f_0\in R$
is uniquely determined by $\ov{f_0}$. Indeed, we have
$$f_0(y,x_0,x)=\psi(-\frac{I(x)}{x_0y})\ov{f_0}(x_0,x,
\grad(\frac{I(x)}{x_0})),$$
provided $x_0\neq 0$ and $|y|\le |x_0|$. Since $f_0$ is also invariant
with respect to the operator $A$, 
we obtain that $f_0$ can be uniquely recovered from 
$\ov{f_0}$.
Therefore, in order to find the spherical vector $f_0$ it suffices
to determine $\ov{f_0}$. Let us normalize $f_0$ by the
condition $\ov{f_0}(1,0,\ldots,0)=1$.

\begin{thm}\label{sphform} 
For $G=E_k$, one has
$$\ov{f_0}(v)=|x_0|^{-s-1}\cdot\frac{1+a|v|^{-s}}{1+a}\cdot\delta_{C(\OO)},$$
for some $a\in\C$. For $G=D_k$ with $k\ge 5$ one has
$$\ov{f_0}(v)=|x_0|^{-1}\max(1,|\frac{x_1}{x_0}|)^{k-4}\cdot
\frac{1+a|v|^{-k+4}}{1+a}\cdot\delta_{C(\OO)}$$
for some $a\in\C$. In both these formulas $v$ is an element of $\La^0(F)$.
\end{thm}

\Pf . Let us consider the action of the double coset 
$T^{\b_0}=IzI$ on $R^I$, where $z=\om^{\vee}_{\b_0}(\pi)$.
We claim that for every $f\in R^I$ one has
$$\ov{T^{\b_0}f}(v)=
q^{2(n+2)-\frac{n+3}{2}}\ov{f}(\pi v)\cdot\delta_{C(\OO)}(v),$$
Indeed, by Lemma \ref{simpleidlem}, we have
$$\ov{zf}(v)=q^{-\frac{n+3}{2}}\ov{f}(\pi v).$$ 
Consider the decomposition $I=H(\OO)B(\OO)U^-(\pi\OO)$,
where $U^-\sub G$ is the subgroup spanned by one-parametric subgroups 
corresponding to negative roots. Since $z$ commutes with $B(\OO)$
and $U^-(\pi\OO)z\sub zU^-(\pi\OO)$, we get $IzI=H(\OO)zI$. More precisely,
for $h\in H(\OO)$ we have $hzI=zI$ if and only if 
$h\in H'\sub H(\OO)$, where
$$H'=e_{\om}(\pi^2\OO)\prod_i e_{\b_i}(\pi\OO)e_{\ga_i}(\pi\OO).$$
When defining the action of $T^{\b_0}$ on $R^I$ we normalize the 
measure on $I$ in such a way that $\vol(I)=1$. Hence,
$$\ov{T^{\b_0}f}=\sum_{h\in H(\OO)/H'}\ov{hzf}.$$
Recalling the explicit formulas for the action of $H(\OO)$ and $z$ on 
$\C(\La^0(F))$ we obtain
$$\ov{T^{\b_0}f}(v)=
q^{2-\frac{n+3}{2}}\sum_{v'\in V(\OO)/V(\pi\OO)}\psi(v\cdot v')\ov{f}(\pi v)$$
(note that since $f$ is $H(\OO)$-invariant, $\ov{f}$ is supported on
$C(\OO)$, so this sum is well-defined). Since sums over $v'$ are
zero unless $v\in C(\OO)$ we get the required formula. 

According to Corollary \ref{heckecor} 
$T^{\b_0}$ acts on $R^I$ semisimply with eigenvalues 
$q^{2(n+1)-s}$ and $q^{2(n+1)-2s}$ in the case $G=E_k$, 
resp. $q^{2(n+1)}$ and $q^{2(n+1)+4-k}$ in the case $G=D_k$, $k\ge 5$.
Furthermore, since 
$\a(\om^{\vee}_{\b_0})=0$ for every simple root $\a\neq\b_0$, by (\ref{comrel})
the action of $T^{\b_0}$ on $R^I$ commutes with the action of the
finite Hecke algebra associated with $G_1$. It follows that 
$T^{\b_0}$ preserves the subspace in $R^I$ consisting of $L(\OO)$-invariant
vectors.

Let us consider the case $G=E_k$. Then we can write the spherical vector
in the form $f_0=f_1+f_2$, where $f_i$ are $L(\OO)$-invariant vectors and
$T^{\b_0}f_i=q^{2(n+1)-is}f_i$ for $i=1,2$.
Therefore, 
$$q^{2(n+2)-\frac{n+3}{2}}\ov{f_i}(\pi v)\cdot\delta_{C(\OO)}(v)=
q^{2(n+1)-is}\ov{f_i},$$
where $i=1,2$.
Now we observe that for every $u\in\C$ there is a unique (up to scalar factor)
function $\phi_u$ on $C(F)$ supported on $C(\OO)$, such that
$\phi_u(v)$ depends only on $|v|$ and 
$$\phi_u(\pi v)\cdot\delta_{C(\OO)}(v)=q^{-u} \phi_a(v),$$
namely, $\phi_u(v)=|v|^u\cdot\delta_{C(\OO)}$.
Applying Proposition \ref{invariantprop} (and the fact that $n=6s+3$), 
we deduce that
the function $|x_0|^{s+1}\ov{f_0}$ is a linear combination of 
the constant function $\delta_{C(\OO)}$ and of $|v|^{-s}\delta_{C(\OO)}$. 
Looking at the explicit description of the Hecke algebra
representation $R^I$ one can observe that all the coefficients
of this linear combination are non-zero. 
Similar argument works in the case of $G=D_k$.  
\ed



\section{Fourier transform over finite field}\label{fouriersec}

In this section we determine the constant $a$ from
Theorem \ref{sphform}. Note that for every $l\ge 0$
the function $f_0(\pi^l,x_0,\ldots,x_n)$ is supported on the set where
$x_i\in\OO$ and depends only on $x_0,\ldots,x_n$ modulo $\pi^l$.
Therefore, we can consider it as a function $f_{\pi^l}$
on $(\OO/\pi^l\OO)^{n+1}$. Now the invariance of $f_0$ with respect
to the operator $S$ implies that $f_{\pi^l}$ is self-dual with
respect to the finite Fourier transform defined by 
$$\FF_{\pi^l}(f)(x)=q^{-\frac{(n+1)l}{2}}\sum_{x'\in(\OO/\pi^l\OO)^{n+1}}
\psi(\frac{x\cdot x'}{\pi^l})f(x').$$
This equality imposes linear equations on the constant $a$. Solving them we can find this constant.
More precisely, we will determine $a$ from the condition
$\FF_{\pi}(f_{\pi})=f_{\pi}$. 
Note that for $|x_0|=1$ we have
$$f_{\pi^l}(x_0,\ldots,x_n)=\psi(-\frac{I(x_1,\ldots,x_n)}{x_0\pi^l}).$$
In fact, for $l=1$ the function
$f_{\pi}$ is the unique Fourier self-dual
extension of $\psi(-\frac{I}{x_0\pi})$ to $(\OO/\pi\OO)^{n+1}$.
For $l>1$ the Fourier self-dual extension of $\psi(-\frac{I}{x_0\pi})$
is not necessarily unique.
However, if we already know functions $f_{\pi^k}$ for $k<l$, then
we can characterize $f_{\pi^l}$ uniquely as a Fourier self-dual
function on $(\OO/\pi^l\OO)^{n+1}$, 
such that for every $k$, $0\le k<l$, one has
$$f_{\pi^l}(u\pi^k,x_1,\ldots,x_n)=\psi(-\frac{I(x_1,\ldots,x_n)}
{u\pi^{k+l}})f_{\pi^k}(0,x_1,\ldots,x_n),$$
where $u\in\OO^*$, $x_1,\ldots,x_n\in\OO/\pi^l\OO$.

We conjecture that the functions $f_{\pi^l}$ are obtained
from the Goresky-MacPherson extensions of the perverse sheaves
over finite field corresponding to functions $\psi(-\frac{I}{x_0\pi^l})$.

\subsection{Case of $G=E_k$}

According to Theorem \ref{sphform} we have the following formula
for the function $f_{\pi}$ on $(\OO/\pi\OO)^{n+1}$:

\begin{equation}
f_{\pi}(x_0,x)=\cases \psi(-\frac{I(x)}{x_0\pi}),\ & x_0\neq 0,\\
q^{s+1},\ & x_0=0, I'_3(x)=0, x\neq 0,\\
q^{s+1}\cdot\frac{1+a q^s}{1+a},\ & x_0=x=0,\\
0,\ &\text{otherwise},\endcases
\end{equation}
where we denote $x=(x_1,\ldots,x_n)$.
Now from the equations 
$$\FF_{\pi}(f_{\pi})(1,0,\ldots,0)=f_{\pi}(1,0,\ldots,0)=1,$$
$$\FF_{\pi}(f_{\pi})(0,\ldots,0)=f_{\pi}(0,\ldots,0)=
q^{s+1}\cdot\frac{1+a q^s}{1+a}$$
we get
\begin{equation}\label{fouriereq0}
\sum_{x_0\neq 0,x}\psi(-\frac{I(x)}{x_0\pi}+\frac{x_0}{\pi})+
q^{s+1}\card\{x\neq 0: I'_3(x)=0\}+
q^{s+1}\cdot\frac{1+a q^s}{1+a}=q^{3s+2},
\end{equation}
$$\sum_{x_0\neq 0,x}\psi(-\frac{I(x)}{x_0\pi})+
q^{s+1}\card\{x\neq 0: I'_3(x)=0\}+
q^{s+1}\cdot\frac{1+a q^s}{1+a}=q^{3s+2}\cdot q^{s+1}\cdot\frac{1+a q^s}{1+a}$$
(we used the equality $n=6s+3$).
Taking the difference, we obtain the following equation for $a$:
\begin{equation}\label{fouriereq1}
q^{3s+2}(q^{s+1}\cdot\frac{1+a q^s}{1+a}-1)=
\sum_{x_0\neq 0,x}\psi(-\frac{I(x)}{x_0\pi})-
\sum_{x_0\neq 0,x}\psi(-\frac{I(x)}{x_0\pi}+\frac{x_0}{\pi}).
\end{equation}
Now we observe that for $i\neq 0$ we have 
$$\lan\ga_0,\b_i\ran=\lan\om,\b_i\ran-\lan\b_0,\b_i\ran=0,$$
hence by Lemma \ref{cubiclem2}, 
$$I(\b_i^{\vee}(x_0)x)=x_0I(x).$$
Therefore,
making the change of variables $(x_0,x)\mapsto (x_0,\b_i^{\vee}(x_0)x)$,
we get
\begin{align*}
&\sum_{x_0\neq 0,x}\psi(-\frac{I(x)}{x_0\pi}+\frac{x_0}{\pi})=
\sum_{x_0\neq 0,x}\psi(-\frac{I(x)}{\pi}+\frac{x_0}{\pi})=-\sum_x \psi(-\frac{I(x)}{\pi})=\\
&-\frac{1}{q-1}\sum_{x_0\neq 0,x}\psi(-\frac{I(x)x_0}{\pi})
\end{align*}
(in the last equality we used the same change of variables again).
On the other hand,
\begin{align*}
&\sum_{x_0\neq 0,x}\psi(-\frac{I(x)}{x_0\pi})=
\sum_{x_0\neq 0,x}\psi(-\frac{I(x)x_0}{\pi})=\\
&(q-1)\cdot\card\{x: I(x)=0\}-\card\{x: I(x)\neq 0\}=
(q-1)q^{6s+3}-q\cdot\card\{x: I(x)\neq 0\}.
\end{align*}
Hence, we can rewrite the RHS of (\ref{fouriereq1}) as
$$(1+\frac{1}{q-1})\sum_{x_0\neq 0,x}\psi(-\frac{I(x)x_0}{\pi})=
q^{6s+4}-\frac{q^2}{q-1}\cdot\card\{x: I(x)\neq 0\}.$$
It remains to compute the cardinality of the set 
$U=\{x\in(\OO/\pi\OO)^n: I(x)\neq 0\}$.
In the case $G=E_6$, this is the set of invertible $3\times 3$-matrices
with coefficients in $\OO/\pi\OO$. Hence, in this case
$$\card U=(q^3-1)(q^2-1)(q-1)q^3.$$ 
For $G=E_7$, $U$ is the set
of non-degenerate skew-symmetric $6\times 6$-matrices with coefficients
in $\OO/\pi\OO$, so 
$$\card U=\frac{\card\GL_6(\OO/\pi\OO)}{\card\Sp_6(\OO/\pi\OO)}=
\frac{(q^6-1)\ldots (q-1)q^{15}}{(q^6-1)(q^4-1)(q^2-1)q^9}=
(q^5-1)(q^3-1)(q-1)q^6.$$
Finally, for $G=E_8$, $U$ is the orbit of 
$(\G_m\times E_6)(\OO/\pi\OO)$, and the
stabilizer is $F_4(\OO/\pi\OO)$, so in this case
\begin{align*}
&\card U=\frac{\card(\G_m\times E_6)(\OO/\pi\OO)}{\card F_4(\OO/\pi\OO)}=\\
&\frac{(q-1)(q^{12}-1)(q^9-1)(q^8-1)(q^6-1)(q^5-1)(q^2-1)q^{36}} 
{(q^{12}-1)(q^8-1)(q^6-1)(q^2-1)q^{24}}=(q^9-1)(q^5-1)(q-1)q^{12},
\end{align*}
where we used the formula for the number of elements of finite Chevalley
groups from \cite{St}, \S 9.
In all cases, the equation (\ref{fouriereq1}) reduces to
$$q^{s+1}\frac{1+a q^s}{1+a}=q^{s+1}+q^{2s+1},$$
which implies that $a=-q^s$.

Note that the formula for $\card U$ can be written uniformly as
\begin{equation}\label{card1}
\card\{x\in\F_q^n: I(x)\neq 0\}=q^{3s}(q^{2s+1}-1)(q^{s+1}-1)(q-1).
\end{equation}
On the other hand, the above equations imply that
$$\sum_{x_0\neq 0,x}\psi(-\frac{I(x)}{x_0\pi}+\frac{x_0}{\pi})=
q^{3s+1}(1-q^{2s+1}-q^{s+1}).$$
Substituting this into equation (\ref{fouriereq0}) with $a=-q^s$ we
deduce the following formula:
\begin{equation}\label{card2}
\card\{x\in\F_q^n: x\neq 0, I'_3(x)=0\}=(q^{2s+1}-1)(q^{2s}+q^s+1).
\end{equation}
We will need the formulas (\ref{card1}) and (\ref{card2}) later in
some global calculations.

\subsection{Case of $G=D_k$, $k\ge 5$}

Recall that in the case $G=D_k$ we have $I(x_1,\ldots,x_n)=x_1Q(x)$, where
$Q$ is a non-degenerate quadratic form in $x=(x_2,\ldots,x_n)$
(where $n=2k-5$).
From Theorem \ref{sphform} we get the following formula for $f_{\pi}$:
\begin{equation}
f_{\pi}(x_0,x_1,x)=\cases \psi(-\frac{x_1Q(x)}{x_0\pi}),\ & x_0\neq 0,\\
q^{k-3},\ & x_0=0, x_1\neq 0, x=0,\\
q,\ & x_0=x_1=0, Q(x)=0, x\neq 0,\\ 
q\cdot\frac{1+aq^{k-4}}{1+a},\ & x_0=x_1=0, x=0\\
0,\ &\text{otherwise}.\endcases
\end{equation}
From the equation $\FF_{\pi}(f_{\pi})(1,0,\ldots,0)=1$ we obtain
\begin{equation}\label{Dkeq1}
\begin{array}{l}
\sum_{x_0\neq 0,x_1,x}\psi(-\frac{x_1Q(x)}{x_0\pi}+\frac{x_0}{\pi})+\\
q^{k-3}(q-1)+q\card\{x\neq 0: Q(x)=0\}+
q\cdot\frac{1+aq^{k-4}}{1+a}=q^{k-2}.
\end{array}
\end{equation}
Making the change of variables $x_1\mapsto x_1x_0$ we can rewrite the
first sum in the LHS as 
$$\sum_{x_0\neq 0,x_1,x}\psi(-\frac{x_1Q(x)}{\pi}+\frac{x_0}{\pi})=-
\sum_{x_1,x}\psi(-\frac{x_1Q(x)}{\pi})=
-q\cdot\card\{x: Q(x)=0\}.$$
Hence, the equation (\ref{Dkeq1}) is equivalent to
$$-q+q^{k-3}(q-1)+q\cdot\frac{1+aq^{k-4}}{1+a}=q^{k-2},
$$
or equivalently,
$$\frac{1+a q^{k-4}}{1+a}=1+q^{k-4}.$$
Therefore, $a=-q^{k-4}$.

\section{Structure of the minimal representation revisited}
\label{revisitsec}

In this section we apply our techniques to obtain a complete description
of the space of smooth vectors in the minimal representation.
In the case of series $E_k$ the main results 
of this section, Theorems \ref{coinvprop} and \ref{structurethm},
were proven by Magaard and Savin in \cite{MS}.

\subsection{$Z(H)(F)$-coinvariants and the line bundle on $C(F)$}

Recall that we denote by $R'\sub R$ the subspace of locally constant
functions with compact support contained in $\{y\neq 0\}$.
Let $R_0$ be the space of $Z(H)(F)$-coinvariants on $R$.
We have a natural surjective map $R_0\ra R/R'$.

\begin{thm}\label{coinvprop}
Assume that $G\neq D_4$. 
Then the natural map $R_0\ra R/R'$ is an isomorphism and
the map $\res:R_0\ra\Ga(C(F),\LL)$ constructed in \ref{conesec}, is an embedding.
\end{thm}

\Pf . Let us denote by $\SS(C(F),\LL)\sub\Ga(C(F),\LL)$ the subspace
of locally constant sections with compact support in $C(F)$. Then
the composition of the map $\res:R\ra\Ga(C(F),\LL)$ with the
projection $\Ga(C(F),\LL)\ra \Ga(C(F),\LL)/\SS(C(F),\LL)$ factors
through the coinvariants of $H(F)$ in $R$.
Thus, we get the following commutative
diagram of $L(F)$-representations
\begin{equation}
\begin{array}{ccc}
R &\lrar{\res} & \Ga(C(F),\LL)\\
\ldar{} & & \ldar{}\\
J(R) &\lrar{\ov{\res}} & \Ga(C(F),\LL)/\SS(C(F),\LL)
\end{array}
\end{equation}
where $J(R)$ denotes the space of coinvariants of $H(F)$ in $R$
(this is the unnormalized Jacquet functor of $R$ corresponding to 
our maximal parabolic subgroup $P$). 

Assume first that $G=E_n$. Then
we can use the fact that $J(R)$ is the direct sum of two distinct
irreducible representations of $G_1(F)$: the trivial
representation and the minimal representation (see \cite{Savin}).
We claim that the map $\ov{\res}$ is not zero on both these summands.
Indeed, as we have seen in the proof of Theorem \ref{sphform},
we can write the spherical vector $f_0$ in the form $f_0=f_1+f_2$, such
that $f_1,f_2\in R^I$ (where $I$ is the Iwahori subgroup of $G(F)$),
$T^{\b_0}f_1=q^{11s+8}f_1$ and
$T^{\b_0}f_2=q^{10s+8}f_2$. Note that the above decomposition of
$J(R)$ as $G_1(F)$-module corresponds to the decomposition of 
$R^I$ into eigenspaces of $T^{\b_0}$. Therefore, our claim follows
from the fact that the images of $f_1$ and $f_2$ 
in $\Ga(C(F),\LL)/\SS(C(F),\LL)$ are non-trivial.
For $G=D_k$ with $k\ge 5$, $J(R)$
is still the direct sum of the trivial 
and of the minimal representations of $G_1(F)$. 
However, in this case the minimal representation
of $G_1(F)$ is the direct sum of two irreducible representations, namely,
the minimal representations of the simple factors of $G_1$. 
Let us write the spherical vector $f_0$ in the form
$f_0=f_1+f_2+f_3$ where $f_i\in R^I$, $i=1,2,3$;
$f_1$ projects to the component of $J(R)$ corresponding to the trivial
representation of $G_1(F)$, $f_2$ to the minimal representation of
$D_{k-2}(F)$ and $f_3$ to the Steinberg representation of $\SL_2(F)$.
The Hecke operator $T^{\b_0}$ acts on these components as follows: 
$T^{\b_0}f_1=q^{3k-4}f_1$, $T^{\b_0}f_2=q^{4k-8}f_2$, 
$T^{\b_0}f_3=q^{3k-4}f_3$. 
Similar argument as above shows that
$\ov{\res}(f_2)\neq 0$ and $\ov{\res}(f_1+f_3)\neq 0$. 
To prove that $\ov{\res}(f_1)\neq 0$ and $\ov{\res}(f_3)\neq 0$ 
we need to consider
the action of the Hecke operator $T^{\a_1}=Iz_1I$, where
$z_1=\om_{\a_1}^{\vee}(\pi)$.
Namely, using Theorem 4.7 of \cite{L1}, it is easy to show 
that $T^{\a_1}f_2=q^{3k-5}f_2$ and that
$f_1$ and $f_3$ belong to a $2$-dimensional subspace $(v_1,v_2)$
of $R^I$ such that $T^{\a_1}v_1=q^{2k-2}v_1$, $T^{\a_1}v_2=q^{3k-5}$.
Thus, it suffices to prove linear independence of
$\ov{\res}[(T^{\a_1}-q^{3k-5})(T^{\b_0}-q^{4k-8})f_0]$ and
$\ov{\res}[(T^{\a_1}-q^{2k-2})(T^{\b_0}-q^{4k-8})f_0]$.
Set $f=(T^{\b_0}-q^{4k-8})f_0$. Using the formula for $\ov{f_0}$ we
easily get 
$$\ov{f}=q^{3k-4}|x_0|^{-1}\max(1,|\frac{x_1}{x_0}|)^{k-4}\cdot\de_{C(\OO)}.$$
We have
$$T^{\a_1}f=\sum_{h\in H'(\OO)/H'(\pi\OO),t\in\OO/\pi\OO}
h\cdot e_{\a_1}(t)z_1f,$$
where $H'=e_{\om}(\G_m)\cdot e_{\ga_0}(\G_m)\cdot e_{\b_1}(\G_m)\cdot
\prod_{i>1} e_{\ga_i}(\G_m)$.
An easy computation similar to the one in the proof of Theorem \ref{sphform}
shows that for $|x_1|<|\pi x_0|$ one has
$$\ov{T^{\a_1}f}=(q^{2k-3}-q^{2k-4}+q^k)|x_0|^{-1}\cdot\de_{C(\OO)}$$
while for $|x_1|>|x_0|$ one has
$$\ov{T^{\a_1}f}=q^{2k-3}\cdot |x_0|^{3-k}\cdot |x_1|^{k-4}\cdot\de_{C(\OO)}.
$$
Let us denote $f'=(T^{\a_1}-q^{3k-5})f$, $f''=(T^{\a_1}-q^{2k-2})f$.
Then for $|x_1|<|\pi x_0|$ we have
$$\ov{f'}=(-q^{3k-5}+q^{2k-3}-q^{2k-4}+q^k)|x_0|^{-1}\cdot\de_{C(\OO)},$$
$$\ov{f''}=(-q^{2k-2}+q^{2k-3}-q^{2k-4}+q^k)|x_0|^{-1}\cdot\de_{C(\OO)}.$$
On the other hand, for $|x_1|>|x_0|$ we have
$$\ov{f'}=(-q^{3k-5}+q^{2k-3})|x_0|^{3-k}\cdot|x_1|^{k-4}\cdot\de_{C(\OO)},$$
$$\ov{f''}=(-q^{2k-2}+q^{2k-3})|x_0|^{3-k}\cdot|x_1|^{k-4}\cdot\de_{C(\OO)}.$$
Since the asymptotics of $|x_0|^{-1}$ and 
of $|x_0|^{3-k}\cdot |x_1|^{k-4}$ as $x\to 0$ in either of the above two
regions are linearly
independent and the determinant
$$\det\left(\matrix (-q^{3k-5}+q^{2k-3}-q^{2k-4}+q^k) & (-q^{3k-5}+q^{2k-3})\\
(-q^{2k-2}+q^{2k-3}-q^{2k-4}+q^k) & (-q^{2k-2}+q^{2k-3})\endmatrix\right)=
(q^k-q^{2k-4})(q^{3k-5}-q^{2k-2})$$
is non-zero (recall that $k\ge 5$), we deduce that 
$\ov{\res}(f')$ and $\ov{\res}(f'')$ are linearly independent.

It follows that in any case the map $\ov{\res}$
is an embedding. Also, it is clear that the map
$$R_0\ra J(R)\oplus\prod_{v\in V(F)\setminus 0} R(v)$$
is an embedding. This immediately implies that the map
$\res:R_0\ra\Ga(C(F),\LL)$ is injective. Therefore, the map $R_0\ra R/R'$
is also injective. Being surjective it must be an isomorphism.
\ed

\subsection{Structure of $Z(H)(F)$-coinvariants}

\begin{lem}\label{centcoinvlem} 
Every locally constant function with compact support on $\La^0(F)$ has
form $\ov{f}$ for some $f\in R$.
\end{lem}

\Pf . Let $h$ be a function with compact support on $\La^0(F)$.
Set 
$$f(y,x_0,\ldots,x_n)=
\cases \psi(-\frac{I(x)}{x_0y})h(x_0,\ldots,x_n), & |y|\le |x_0|,\\
0 & |y|>|x_0|.\endcases$$
Then $Af$ has compact support in $F^*\times F^{n+1}$, so by Lemma
\ref{supportlem}, $Af\in R'$, hence $f\in R$. On the other hand,
clearly we have $h=\ov{f}$.
\ed

Let us denote by $\Ga_{\smooth}(C(F),\LL)$ the space of locally constant
sections of $\LL$. Recall that we denote by $\SS(C(F),\LL)\sub
\Ga_{\smooth}(C(F),\LL)$ the subspace of sections with compact support in
$\C(F)$.

\begin{thm}\label{structurethm} 
Assume that $G\neq D_4$. Then 
$R_0$ can be identified with an $L(F)$-submodule in the space
$\Ga_{\smooth}(C(F),\LL)$ containing $\SS(C(F),\LL)$.
In the case $G=E_n$ we have an exact sequence of $L(F)$-modules
$$0\ra\SS(C(F),\LL)\ra R_0\ra J(R)=\C(|\om|^{s+1})\oplus 
R_{G_1}(|\om|^{\frac{s}{2}+1})\ra 0$$
where $R_{G_1}$ is the minimal representation of $G_1$,
$J(R)$ are convariants of $H(F)$ on $R$.
In the case $G=D_k$ with $k\ge 5$ there is an exact sequence
$$0\ra\SS(C(F),\LL)\ra R_0\ra J(R)=\C(|\om|^{\frac{k}{2}-1})\oplus 
R_{\SL_2}(|\om|^{\frac{k}{2}-1})\oplus R_{D_{k-2}}(|\om|)\ra 0.$$
\end{thm}

\Pf . Lemma \ref{centcoinvlem} implies that $\SS(C(F),\LL)$
is indeed a submodule of $R_0$.
We claim that quotient $R_0/\SS(C(F),\LL)$ coincides with the space
of coinvariants of $V(F)$ on $R_0$. Indeed, by Lemma \ref{supportlem}(a)
for every $f\in R_0$ and every linear functional $v^*$ on $V$ one has
$\psi(v^*)f-f\in\SS(C(F),\LL)$. On the other hand, it is easy to see
that there are no coinvariants for the action of $V(F)$ on $\SS(C(F),\LL)$.
This proves our claim. Therefore, the quotient
$R_0/\SS(C(F),\LL)$ coincides with the Jacquet functor of $R$, $J(R)$.
Now we obtain the required exact sequences (not taking into account the action
of $Z(L)$) from the known decomposition of $J(R)$ into irreducible $G_1(F)$-modules 
(see \cite{Savin}).
The action of the center of Levi on $J(R)$ is determined using the
information about the action of $T^{\b_0}$ on $R^I$ and the explicit
form of the spherical vector $f_0$ from Theorem \ref{sphform}.
Indeed, assume for example that $G$ is of type $E_n$. Then
we know that $f_0=f_1+f_2$, where $f_1$ and $f_2$ are nonzero
vectors in $R^I$ such that $T^{\b_0}f_1=q^{11s+8}f_1$ and
$T^{\b_0}f_2=q^{10s+8}f_2$. Furthermore, the function $\ov{f}_1$ is
proportional to $|x_0|^{-s-1}\de_{C(\OO)}$, while $\ov{f}_2$
is proportional to $|x_0|^{-s-1}|v|^{-s}\de_{C(\OO)}$.
Considering asymptotics at the vertex of the cone $C$ we obtain that
$$\ov{zf}_1\equiv q^{-2s-2}\ov{f}_1\mod\SS(C(F),\LL),$$
$$\ov{zf}_2\equiv q^{-s-2}\ov{f}_2\mod\SS(C(F),\LL).$$
Since $f_1$ (resp., $f_2$) 
projects to a non-zero vector in the trivial (resp., non-trivial) 
component of $J(R)$, this proves our claim. The case of the series $D$
is similar. 
\ed

For every $a\in F$ let us denote by $R_a$ the quotient 
$$R_a=R/\{e_{\om}(t)f-\psi(at)f|\ f\in R, t\in F\}$$
(for $a=0$ this is just the space of $Z(H)(F)$-coinvariants considered above).
The above Theorem can be complemented by the following result.

\begin{prop}\label{irredprop} For $a\in F^*$
the representation of $H(F)$ on $R_a$ is irreducible.
\end{prop}

\Pf . Lemma \ref{supportlem} implies that the map
$f(y,x_0,\ldots,x_n)\mapsto f(a,x_0,\ldots,x_n)$ induces an isomorphism
of $R_a$ with the space of locally constant functions with constant
support on $\La(F)$. Furthermore, the action of $H(F)$ on $R_a$
is compatible with the standard Schr\"odinger representation of $H(F)$
associated with the character $\psi(a\cdot ?)$ on $F$, which
is known to be irreducible.
\ed

\section{Archimedian case}\label{archsec}

Some of the results described in sections \ref{constrsec}-\ref{revisitsec} 
have analogues for real and complex groups.
The minimal unitary representation 
in these cases is given by the same formulas where
$\psi$ gets replaced by the standard additive character $x\mapsto\exp(\pi ix)$ for $F=\R$ and 
by $z\mapsto\exp(\pi i(z+\ov{z}))$ for $F=\C$.
Note that for every subgroup $G'\sub G$ the ring $\SS(G'(F))$ of Schwartz
functions on $G'(F)$ (with convolution as a product) acts
naturally on the space of smooth vectors $R$. 
We denote by $\varphi*v\in R$ the result of action of a Schwartz function
$\varphi$ on $v\in R$. By irreducibility of some smooth representation
of $G'(F)$ we always mean the absence of non-trivial
proper $\SS(G'(F))$-invariant closed subspaces. 

Of course, when formulating archimedian analogues of
our results one should replace locally constant functions by 
$C^{\infty}$-functions. Another important change is that one has to impose
certain growth conditions when considering functionals on $R$.
Namely, let $G'\sub G$ be a subgroup.
Then we say that a linear functional $f^*:R\ra\C$ is $G'$-{\it tempered}
if for every $f\in R$ the functional
$\varphi\mapsto f^*(\varphi*f)$ on $\SS(G'(F))$ is continuous.
This is equivalent to the condition that the functional $f^*$ is
continuous on $R$ with respect to the topology defined by the collection
of norms $\|f\|_u=\|uf\|$, where $u\in U({\frak g'})$ (where ${\frak g'}$ is
the Lie algebra of $G'(F)$).
Dually we have to change the notion of twisted coinvariants.
The correct replacement for
the maximal quotient of $R$ on which $G'(F)$ acts by some character
$\chi:G'(F)\ra\C^*$ is 
$$J_{G',\chi}(R):=R/(\varphi*f-\int_{g'\in G'}\chi(g')\varphi(g')f)$$
where $f$ runs through $R$, $\varphi$ runs through the Schwartz space
of $G'(F)$, $(S)$ denotes the closure of the span of the set of vectors
$S$. We have to use such quotients when defining archimedian analogues
of the spaces $R(v)$ and $R_0$, or the space of coinvariants of
$H(F)$ on $R$.
Note that if a functional $f^*$ is $G'$-tempered
and satisfies $f^*(sv)=\chi(s)f^*(v)$ for some character $\chi$ of $G'(F)$
then $f^*$ descends to a functional on the space $J_{G',\chi}(R)$.

We refer to \cite{KPW} for the derivation of the formula
for the spherical vector in $R$ in the cases $F=\R$ or $F=\C$.
Below we will consider analogues of some of the results
of sections \ref{constrsec}, \ref{structuresec} and \ref{revisitsec}.

\subsection{Smooth vectors and the $L(F)$-equivariant line bundle}

Let us start by formulating analogues of parts (a) and (b)
of Lemma \ref{supportlem}.

\begin{lem}\label{archlem1}
(a) For every $f\in R$ and every differential operator $D$ in
$(y,x_0,\ldots,x_n)$ with constant coefficients, one has $D(f)\in R$.

\noindent
(b) For every differential operator $D$ as above and every $N$ one has
$D(f)|(y,x_0,\ldots,x_n)|^N\to 0$ as $|(y,x_0,\ldots,x_n)|\to\infty$.
\end{lem}

\Pf . This follows easily from
the formula for the spherical vector obtained in \cite{KPW}
and from the fact that every smooth vector can be obtained 
from the spherical vector by the action of $U({\frak p})$,
the universal enveloping algebra of the Lie algebra of $P$.
\ed

The following proposition is an analogue of part (c) of Lemma \ref{supportlem}.

\begin{prop}\label{Rinftylem}
Let $R_{\infty}$ denote the space of all $C^{\infty}$-functions
$f$ on $F^*\times F^{n+1}$ such that
$D(f)|(y^{-1},y,x_0,\ldots,x_n)|^N\to 0$ 
as $|(y^{-1},y,x_0,\ldots,x_n)|\to\infty$ for
every $N>0$ and every differential operator $D$ in $(y,x_0,\ldots,x_n)$
with constant coefficients. Then $R_{\infty}$ is contained in $R$.
\end{prop}

In the archimedian case the proof is based on the following lemma.

\begin{lem}\label{irrlem} The action of $P(F)$ on $R_{\infty}$ 
is irreducible, i.e., there are no
non-trivial proper $\SS(P(F))$-invariant closed subspaces in $R_{\infty}$.
\end{lem}

\Pf . Let $M$ be the space of Schwartz functions on $F^{n+1}$. For
every $f\in R_{\infty}$ we can consider the restriction $f|_{y=a}$
as an element of $M$. Furthermore, this defines an
isomorphism of $R_{\infty}$ with the Schwartz space of $M$-valued functions on $F^*$.
Now our assertion follows easily from the irreducibility of the
Schr\"odinger representation of $H(F)$.
\ed

\noindent
{\it Proof of Proposition \ref{Rinftylem}}.
By Lemma \ref{irrlem}, it suffices to check that $R_{\infty}\cap R$
is non-trivial.
To construct a non-zero element in $R_{\infty}\cap R$
we want to pick any non-zero $f\in R$ and
replace it by $\varphi*f$, where $\varphi$ is a Schwartz function
on $Z(H)(F)=F$. By the definition,
$\varphi*f=\hat{\varphi}(y)f$, where $\hat{\varphi}$ is the
Fourier transform of $\varphi$. Let $t\in F^*$ be such that
$f|_{t\times F^{n+1}}\neq 0$. Pick $c>0$ such that
$c<t<c^{-1}$ and then choose $\varphi$ in such
a way that $\hat{\varphi}$ is supported on the set $c\le |y|\le c^{-1}$
and $\hat{\varphi}(t)\neq 0$. Then $\varphi*f$ 
will be a non-zero element in $R_{\infty}\cap R$.
\ed

Now let us consider analogues of the results of section \ref{conesec}.
Using the formula for the spherical vector $f_0$ and the fact that
$R=U({\frak p})f_0$
one can easily see that the right-hand side of the
formula (\ref{fbareq}) is well-defined. 
Moreover, one can check that the functional
$f^*_0(f)=\lim_{y\to 0}f(y,1,0)$ is $H$-tempered.
On the other hand, for $F=\R$ or $\C$ the twisted action of $G(F)$
on $C(F)$ is still transitive as
follows from the proof of Lemma \ref{globaltrans} below.
Now it is easy to see that Proposition \ref{onedimlem} holds
in the archimedian case (the proof should be modified
similarly to the proof of Proposition \ref{archlem2} below).

The definition of the $L(F)$-equivariant line bundle $\LL$ over
$C(F)$ and of the identifications $R(v)\simeq\LL|_v$ for $v\in C(F)$
transfers to the archimedian case without any changes.
As in section \ref{conesec},
the natural projections $R\ra R(v)$ glue into a map
$\res:R\ra\Ga(C(F),\LL)$. This map factors through the quotient
$R\ra R_0:=J_{Z(H)(F),1}(R)$ and its image belongs to the subspace
$\SS_{\infty}(C(F),\LL)\sub\Ga(C(F),\LL)$ consisting of sections
$f$ of $\LL$ such that
for every $u\in U({\frak l})$ the section $uf(x)$ is
rapidly decreasing at infinity
(i.e. as $|x|\to\infty$, where $x\in C(F)\sub V$), where
${\frak l}$ is the Lie
algebra of $L$.
Indeed, this can be seen from the formula for the spherical vector in $R$.

\subsection{Twisted coinvariants with respect to $Z(H)$}

Now let us formulate an
archimedian analogue of Proposition \ref{irredprop}.

\begin{prop}\label{archlem2} For every $a\in F^*$
the map $f(y,x_0,\ldots,x_n)\mapsto f(a,x_0,\ldots,x_n)$
induces an $H(F)$-equivariant isomorphism of the space
$R_a:=J_{Z(H),\psi_a}(R)$ with the Schr\"odinger representation
of $H(F)$ on the Schwartz space of $F^{n+1}$
(here $\psi_a$ is the character of $Z(H)(F)\simeq F$ given by
$\psi_a(t)=\psi(at)$). In particular, the representation of
$H(F)$ on $R_a$ is irreducible.
\end{prop}

\Pf . In view of Proposition \ref{Rinftylem}, 
it suffices to check 
that every $f\in R$ such that $f(1,x_0,\ldots,x_n)\equiv 0$
maps to zero under the natural projection $R\ra J_{Z(H),\psi_a}(R)$.
By the definition, the latter space is the quotient of
$R$ by the closure of the span of elements of the form
$$\varphi*h-\hat{\varphi}(a)h=\hat{\varphi}(y)h-\hat{\varphi}(a)h$$
where $h\in R$,
$\varphi$ is a Schwartz function on $Z(F)=F$, $\hat{\varphi}$
is its Fourier transform. 
Therefore, it suffices to check
that every $f\in R$ such that $f(a,x_0,\ldots,x_n)\equiv 0$
is a linear combination of elements of the form
$\varphi*h-\hat{\varphi}(a)h$, where $h\in R$, $\varphi$ is a Schwartz function
on $F$. Let $\varphi_1$ be a $C^{\infty}$-function on $F$ such that
$\hat{\varphi}_1(a)=1$ and $\hat{\varphi}_1(y)\equiv 0$
for $|y-a|>\eps$ for some small $\eps>0$. Replacing $f$ by
$$f+(\varphi_1*f-\hat{\varphi}_1(a)f)=\hat{\varphi}_1(y)f$$ 
we can assume
that $f(y,x_0,\ldots,x_n)\equiv 0$ for $|y-a|>\eps$.
Now let us set
$$h(y,x_0,\ldots,x_n)=f(y,x_0,\ldots,x_n)/(y-a)\ \text{for}\ y\neq a.$$
Then $h$ extends to a $C^{\infty}$-function and clearly $h\in R_{\infty}\sub R$.
Let $\varphi_2$ be a Schwartz function on $F$ such
that $\hat{\varphi}_2(y)\equiv y$ for $|y-a|\le\eps$. Then
we have $f=\varphi_2*h-\hat{\varphi}_2(a)h$.
\ed

We believe that the following
analogues of Theorems \ref{coinvprop} and \ref{structurethm} should also hold in the
archimedian case.

\vspace{2mm}

\noindent
{\bf Conjecture 1.} If $G\neq D_4$ then the map
$R_0\ra\SS_{\infty}(C(F),\LL)$ is injective.

\noindent
{\bf Conjecture 2.} If $G\neq D_4$ then one has
exact sequences as in Theorem \ref{structurethm},
where $\SS(C(F),\LL)\sub\SS_{\infty}(C(F),\LL)$ is the space of sections $f$
such that for every $u\in U({\frak l})$ the
section $uf(x)$ is rapidly decreasing as $|x|\to\infty$ and as $|x|\to 0$.

D.~Barbasch informed us that he checked the required decomposition
of $J(R)$ for the case of series $D_k$.
In \ref{formEsec} we will assume the validity of the above conjectures for the groups
$E_6$, $E_7$ and $E_8$. However, this assumption will not be used anywhere else.

\bigskip

\centerline{\bf Part II. Global picture}

\section{Form of the automorphic functional}\label{globalsec}

Let $K$ be a global field, $G$ be a simply connected simple 
split group of type $D_k$ or $E_k$ over $K$. Let $\mA$ be the ring of adeles, 
$R_{\mA}=\otimes_v R_v$
the minimal representation of $G(\mA)$ 
which is the restricted tensor product of
the local minimal representations $R_v$ with respect
to spherical vectors. 
Below we are going to apply our local results to obtain information about
(and in some cases determine explicitly)
the unique $G(K)$-invariant tempered functional on $R_{\mA}$.
The existence of such a functional follows from the result
of \cite{GRS} (in the case of $D_4$ this is also proven in \cite{K}).
The uniqueness (up to proportionality) is an easy consequence of
the description of all $P(K)$-invariant tempered functionals on $R_{\mA}$
that we will give below. 
Note that
in the case $G=D_k$ the uniqueness follows also from the work \cite{GRS},
while in the case $G=E_7$ it was proven in \cite{Gur}.

\subsection{$P(K)$-invariant functionals}

Let $\psi$ be a non-trivial character of $\mA$, trivial on $K$,
and let $\psi_v:K_v\ra\C^*$ be its local components. 
For every finite place $v$ we can consider
the realization of the minimal representation associated 
with $\psi_v$ (see section \ref{twistedsec}). For an archimedian place $v$ such
a realization also exists. Namely, one should represent the character
$\psi_v$ in the form $t\mapsto\exp(-i \Tr_{K_v/\R}(a_v t))$
for some unique $a_v\in K_v^*$ and then conjugate the standard minimal
representation of $G(K_v)$ (considered in \cite{KPW}) by
the operator $f(y,x_0,x)\mapsto f(a_vy,a_vx_0,\ldots,a_vx_n)$.
In all these local realizations $R_v$ the centers of $H(K_v)$ act by multiplication with the
characters $t\mapsto\psi_v(ty)$. 
The restricted product $R_{\mA}=\otimes R_v$ can be viewed
as a subspace of functions of $(y,x_0,x_1,\ldots,x_n)$, where
$y$ is an idele, $x_i$'s are adeles.

It is sometimes convenient to extend the numbers $a_v$ defined above for archimedian places to
a differential idele $a=(a_v)$ attached to $\psi$ (see \cite{Weil-BNT}, 
VII-2, Def.4), so that for all finite places $v$ the additive characters 
$t\mapsto\psi_v(a_v^{-1}t)$ are of order $0$.

We are going to study linear functionals on $R_{\mA}$ that are invariant with respect
to $P(K)$. This will allow us to determine a general form of the $G(K)$-invariant
functional on $R_{\mA}$. Note that whenever 
we talk about functionals on $R_{\mA}$
we mean only $P$-tempered functionals and we
denote the space of such functionals by $R_{\mA}^*$.
The simplest $P(K)$-invariant $P$-tempered functional on $R_{\mA}$ is
$$\th_1(f)=\sum_{y\in K^*,(x_0,x)\in K^{d+1}}f(y,x_0,x).$$
Indeed, clearly, $\th_1$ is invariant with respect to the action of
$H(K)$ and of Borel subgroup $B(K)$ in $L(K)$. On the other hand, 
Poincar\'e summation formula shows its
invariance with respect to the operator $S$. It remains to note that
$L(K)$ is generated by $B(K)$ and by $S$.
The convergence of $\th_1(f)$ for $f\in R_{\mA}$ 
is clear in the functional field case since only a finite number of terms will be non-zero.
In the number field case the convergence follows from the fact that 
every smooth vector at archimedian place is obtained
from the spherical vector by the action of the universal enveloping algebra of $P$
and from the explicit formulas for archimedian spherical vectors in \cite{KPW} 
(due to the exponential decay at infinity of the relevant Bessel functions). 
Also, since the action of $P$ is given essentially by the formulas of
Weil representation, the results of \cite{Weil-un} imply
that the functional $\th_1$ is $P$-tempered.

Another $P(K)$-invariant functional is obtained by using the theory of
section \ref{conesec}. Recall that for every place $v$ we have
an $L(K_v)$-equivariant complex line bundle $\LL_v$
over the cone $C(K_v)$. Namely, we use the 
identification $L(K_v)/S_{\b_0}(K_v)\wt{\ra} C(K_v)$ given by the (twisted)
action on $e_{\b_0}\in C(K_v)$ and define $\LL_v$ to be the line bundle corresponding to
the character $|\de|_v^{-\frac{1}{2}}$ of $S_{\b_0}(K_v)$.
For every $v$ we have a unique (up to constant)
$P(K_v)$-equivariant map $\res_v:R_v\ra \Ga(C(K_v),\LL_v)$
where $H(K_v)$ acts on sections of $\LL_v$ by the rule (\ref{Vaction}).
Let us normalize these maps by the condition
$$\res_v(f_{0,v})(a_v^{-1}e_{\b_0})=\om_{\b_0}^{\vee}(a_v)(1)$$
where $f_{0,v}\in R_v$ are spherical vectors normalized as in section \ref{twistedsec}.
This is equivalent to the normalization of section \ref{twistedconesec}.
Note that for all finite places
the canonical trivializations of $\LL_v$ at $e_{\b_0}$ extend naturally to 
$L(\OO_v)$-equivariant trivializations of $\LL_v$ over $C(\OO_v)$.
Therefore, we can define the $L(\mA)$-equivariant line bundle $\LL_{\mA}$ over $C(\mA)$ 
with the fiber $\otimes_v \LL_{z_v}$ over $z\in C(\mA)$ (almost all factors
of the tensor product are trivial). 

\begin{lem}\label{globaltrans} The action of $L(K)$ on $C(K)$ is
transitive.
\end{lem}

\Pf . As in Lemma \ref{transitivelem} we see that the statement
is equivalent to the triviality of the kernel (i.e., the preimage of
the trivial class) of the map
$H^1(K,G_0)\ra H^1(K,G_1)$. In the case when $G$ is of type $E_6$ or $E_7$
this follows from the triviality of $H^1(K,G_0)$.
In the case $G=E_8$ (resp., $G=D_k$) the map in
question is $H^1(K,E_6)\ra H^1(K,E_7)$
(resp. $H^1(K,\SO_{2k-6})\ra H^1(K,\SO_{2k-4})$). In the case of series
$D$ the assertion follows easily from the theory of quadratic forms.
It remains to show that the embedding of split simply connected groups 
$E_6\ra E_7$ induces an embedding of the groups $H^1$. By the Hasse principle
this is equivalent to the similar assertion about groups over $\R$.
So let us denote by $E_6\ra E_7$ the standard embedding of split
simply connected groups over $\R$. Let
$\a$ be a $1$-cocycle of $\Z/2\Z$ with values in $E_6$
representing a non-trivial element in $H^1(\R,E_6)$ and let
$E_6^{\a}\sub E_7^{\a}$ be the twists of our groups by $\a$.
It suffices to prove that $E_7^{\a}$ is nonsplit. Note
that since the map $H^1(\R,E_6)\ra H^1(\R,E_6/A)$, where
$A$ is the center of $E_6$, is an isomorphism, the group
$E_6^{\a}$ is a nonsplit form. Let $Z\sub E_7$ be the connected
component of $1$ in the centralizer of $E_6$ in $E_7$. Then
$Z$ is a split $1$-dimensional torus and its twist $Z^{\a}\sub E_7$
is isomorphic to $Z$. Let $G\sub E_7$ be the connected component of $1$
in the centralizer of $Z$ in $E_7$. Then $G=Z\cdot E_6$ and
$G^{\a}=Z\cdot E_6^{\a}$ is the centralizer of $Z$ in $E_7^{\a}$.
If $E_7^{\a}$ were split then $G^{\a}$ would contain a split
$7$-dimensional torus which is a contradiction since $E_6^{\a}$
is nonsplit.
\ed

\begin{rem} Our proof generalizes to the case of more general
pairs of split groups 
$G_0\ra G$ over some field $k$, where $G_0$ is the semisimple part of
the Levi component
of a maximal parabolic subgroup in $G$. The statement is that
the induced map $H^1(k,G_0)\ra H^1(k,G)$ sends elements inducing non-split
forms of $G_0$ to elements inducing non-split forms of $G$.
\end{rem}

Using the above lemma we can construct a natural trivialization of the line bundle $\LL_{\mA}$
over $C(K)\sub C(\mA)$. Indeed, for every $z\in C(K)$ we can choose $g\in L(K)$ such that
$z=g*e_{\b_0}$. The corresponding isomorphism 
$$(\LL_{\mA})|_z\simeq (\LL_{\mA})|_{e_{\b_0}}\simeq\C$$
does not depend on $g$ because of the triviality of $|\de|$ on $S_{\b_0}(K)$.
Now for every $f=\otimes f_v\in R_{\mA}$ we have a well-defined section
$\res(f)$ of $\LL_{\mA}$:
$$\res(f)(z)=\prod_v \res_v(f_v)(z_v)$$
where $z\in C(\mA)$. Using the trivialization of $\LL_{\mA}$ over $C(K)$ we can define 
the functional
$$\th_2(f)=\sum_{z\in C(K)}\res(f)(z).$$
It is clear that this functional is $P(K)$-invariant.
The convergence for $f\in R_{\mA}$ is checked in the same way as for $\th_1$: if $K$ is a functional
field, then almost all terms in $\th_2(f)$ are zero; if $K$ is a number field
one has to use the explicit formula for the spherical vectors at archimedian places.
The fact that $\th_2$ is $P$-tempered is clear from the explicit formulas
for the action of $P(F)$ on $\res(f)$.

\begin{thm}\label{Pinvthm}  
The space of $P(K)$-invariant $P$-tempered functionals on $R_{\mA}$ 
is generated by $\th_1$, $\th_2$
and by the subspace of $P(K)H(\mA)$-invariant functionals.
\end{thm} 

\begin{cor}\label{uniquecor}
The space of $G(K)$-invariant $P$-tempered
functionals on $R_{\mA}$ is $1$-dimensional.
\end{cor}

\noindent
{\it Proof of the corollary}. 
Since the functional $\th_2$ is $Z(H)(\mA)$-invariant, we obtain
from the theorem that the subspace $(R_{\mA}^*)^P(K)Z(H)(\mA)\sub (R_{\mA}^*)^{P(K)}$
has codimension $1$. Since $G(K)$ and $Z(H)(\mA)$ generate $G(\mA)$, the projection
$$(R_{\mA}^*)^{G(K)}\ra (R_{\mA}^*)^{P(K)}/(R_{\mA}^*)^{P(K)Z(H)(\mA)}$$
is an embedding. On the other hand, from \cite{GRS} we know that there exists
a non-zero $G(K)$-invariant tempered functional on $R_{\mA}$.
\ed


Let $R_{\mA,0}$ be the restricted tensor product of the local spaces $R_{v,0}$
of $Z(H)$-coinvariants (recall that at infinite places we use the definition
of section \ref{archsec}). Let
$R'_{\mA}\subset R_{\mA}$ be the kernel of the map $R_{\mA}\ra
R_{\mA,0}$ and let $R_{\mA}^{\prime *}$ be the space of $H$-tempered
linear functionals on $R_{\mA}$.

\begin{lem}\label{invlem1} One has $\dim(R_{\mA}^{\prime *})^{P(K)}\leq 1$.
\end{lem}

\Pf . Consider the space $(R_{\mA}^{\prime *})^{Z(H)(K)}$. Since the
quotient group $\mA/K$ is compact, we have 
$$(R_{\mA}^{\prime *})^{Z(H)(K)}\sub\prod_\alpha J_{Z(H),\alpha}(R_{\mA})^*,$$
where $\alpha$ runs through non-trivial characters of the group  
$\mA/K$ and the twisted coinvariant spaces $J_{Z(H),\alpha}(R_{\mA})$ are defined
in the same way as in section \ref{archsec} using the ring of Schwartz functions
on $Z(H)(\mA)$, the functionals on $J_{Z(H),\alpha}(R_{\mA})^*$ are assumed to
be $H$-tempered.
It follows from Propositions \ref{irredprop} and \ref{archlem2}
that the natural action of the group $H(\mA)$ on 
$J_{Z(H),\alpha}(R_{\mA})$ is
irreducible and equivalent to the space of smooth vectors in the Schr\"odinger representation
of $H(\mA)$ corresponding to the character $\alpha$. Hence,
the spaces $J_{Z(H),\alpha}(R_{\mA})^*$ are $1$-dimensional (see \cite{Weil-un}). 
Since the natural action of the group $P(K)$ on the set of non-trivial characters 
$\alpha$ of the group $Z(H)(\mA)/Z(H)(K)\simeq\mA/K$ is {\it transitive},  
we obtain that $\dim(R_{\mA}^{\prime *})^{P(K)}\le 1$.
\ed


Let $R''_{\mA}\subset  R_{\mA,0}$ be the
kernel of the map $R_{\mA,0} \ra J(R_{\mA})$, where
$J(R_{\mA})$ is the Jacquet functor with respect to $H(\mA)$
(defined as in section \ref{archsec}). Let $R_{\mA}^{\prime\prime*}$
be the space of $P$-tempered functionals on $R''_{\mA}$.

\begin{lem}\label{invlem2}
$\dim(R_{\mA}^{\prime\prime*})^{P(K)}\le 1$.
\end{lem}
 
\Pf . The group $V(\mA)$ acts naturally on the space  
$R_{\mA,0}$ and preserves the subspace  $R''_{\mA}$. 
Consider the space $(R_{\mA}^{\prime\prime*})^{V(K)}$. Since the 
quotient group $V(\mA )/V(K)$ is compact, we have 
$$(R_{\mA}^{\prime\prime*})^{V(K)}\sub\prod_{\beta \in B}
J_{V(\mA),\beta}(R_{\mA,0})^*$$
where $\beta$ runs through the set $B$ of non-trivial characters of the
group  $V(\mA)/V(K)$.
We can identify the set $B$ with
$V(K)-\{0\}$. It follows from the adelic version of
Proposition \ref{onedimlem}
that $J_{V(\mA),\beta}(R_{\mA,0})^*=\{ 0\}$
if $\beta $ does not belong to $C(K)$ and that 
$\dim(J_{V(\mA),\beta}(R_{\mA,0})^*=1$ 
for $\beta \in C(K)$. Therefore, Lemma \ref{globaltrans} implies that
$\dim (R_{\mA}^{\prime\prime*})^{P(K)}\le 1$. 
\ed




\noindent
{\it Proof of Theorem \ref{Pinvthm}.}
According to Lemma \ref{invlem2} 
$P(K)H(\mA)$-invariant functionals form a subspace of codimension
at most $1$ in $(R_{\mA,0}^*)^{P(K)}$. Therefore, the latter space is
generated by $\th_2$ together with the subspace $(R_{\mA}^*)^{P(K)H(\mA)}$. 
On the other hand, by Lemma \ref{invlem1},
$(R_{\mA,0}^*)^{P(K)}$ is a subspace in $(R_{\mA}^*)^{P(K)}$
of codimension at most $1$. Therefore, $(R_{\mA}^*)^{P(K)}$ is generated
by this subspace and by $\th_1$. 
\ed

\subsection{Automorphic functional}

Now we can apply our study of $P(K)$-invariant functionals on $R_{\mA}$ to determine
the form of the $G(K)$-invariant (tempered) functional.
Recall that we denote by $a=(a_v)$ a differential idele corresponding to $\psi$.
It is well-known that in the number field case $|a|=|D|^{-1}$, where $D$ is the discriminant,
while in the functional field case $|a|=q^{2-2g}$, where $q$ is the number of elements in
the field of constants, $g$ is genus (see \cite{Weil-BNT}, VII-2, prop. 6).

\begin{thm}\label{mainglobalthm} 
A non-zero $G(K)$-invariant tempered functional on $R_{\mA}$ is a scalar multiple
of a functional of the form
$$\th(f)=\th_1(f)+|a|^{-\frac{n+3}{2}}\cdot\th_2(f)+\a(f),$$
where $\a$ is a $P(K)H(\mA)$-invariant functional.
\end{thm}

\Pf . 
Using Theorem \ref{Pinvthm} we obtain
that a $G(K)$-invariant tempered functional $\th$ can be written as follows
$$\th=c_1\th_1+c_2\th_2+\a$$ 
where $\a$ is $H(\mA)$-invariant.
Note that if $c_1=0$ then $\th$ is $Z(H)(\mA)$-invariant.
Since $G(K)$ and $Z(H)(\mA)$ generate $G(\mA)$, this is impossible.
So we can assume that $c_1=1$.

Let $S$ be a non-empty set of places containing all archimedian places, such
that all the characters $\psi_v$ for $v\not\in S$ have order $0$.
We define a family of function $f^{\eps}=\otimes_v f_v^{\eps}\in R_{\mA}$ 
depending on a real positive parameter $\eps$ as follows.
For places $v\not\in S$  we set
$f^{\eps}_v=f_{0,v}$ (the spherical vector). For $v\in S$ finite we set
$$f_v^{\eps}=\delta_{|y|_v=1;|(x_0,\ldots,x_n)|_v<\eps},$$
where $\delta_{U}$ denotes the characteristic function of the set $U$. 
For archimedian $v\in S$ we define $f_v^{\eps}$ to be a family of
functions in $R_v$ such that
the support of $f_v^{\eps}$ is contained
in the set $\{(y,x_0,\ldots,x_n):\ ||y|_v-1|<\eps, |x_i|_v<\eps\}$ and
such that $f_v^{\eps}(y,0,\ldots,0)=1$ if $|y|_v=1$.

From the $G(K)$-invariance of the functional $\th$ we get
\begin{equation}\label{Aequation}
\th(f^{\eps})=\th(Af^{\eps})
\end{equation}
where $A$ is the operator (\ref{operatorA}).
We are going to calculated both sides of this equality (or
rather their limits for $\eps\to 0$).

Since $\res_v(f_v^{\eps})=0$ for $v\in S$ we have $\th(f^{\eps})=\th_1(f^{\eps})$.
Furthermore, in the sum defining $\th_1(f^{\eps})$ the non-zero terms are only those
for which $y,x_i\in O_v$ ($i=0,\ldots,n$) for all $v\not\in S$ and $||y|_v-1|<\eps$, $|x_i|_v<\eps$ 
for all $v\in S$. If $\eps$ is sufficiently small this implies
that $x_0=\ldots=x_n=0$ while
$|y|_v=1$ for all $v$, hence $y$ is a root of unity in $K$. 
Note that for such $y$ we have $f^{\eps}(y,0,\ldots,0)=1$, so we obtain
$$\th(f^{\eps})=|\mu(K)|,$$
where $\mu(K)\sub K^*$ is the set of all roots of unity in $K$.

Now let us calculate the limit of the RHS of (\ref{Aequation}) as $\eps\to 0$. For $v\not\in S$
we have $Af^{\eps}_v=f_{0,v}$, while for $v\in S$ the support of $Af^{\eps}_v$ is contained in the set
$\{(y,x_0,\ldots,x_n):\ ||x_0|_v-1|<\eps, |y|_v<\eps, |x_i|_v<\eps, i\ge 1\}$ 
It follows that for sufficiently small $\eps$ we have $\th_1(Af^{\eps})=0$ (since
in this sum $y$ is required to be non-zero). 
Also, since the support of $\res_v(Af_v)$ for $v\in S$ is contained in the set
$\{(x_0,\ldots,x_n):\ ||x_0|_v-1|<\eps, |x_i|_v<\eps, i\ge 1\}$ (considered as
a subset of $\La^0(K_v)\subset C(K_v)$), we have $\a(Af^{\eps})=0$.
Arguing as before we see that in the sum defining $\th_2(Af^{\eps})$ only the terms
with $x_0\in\mu(K)$, $x_1=\ldots=x_n=0$, survive.
Now applying (\ref{resnormal}), we get
$$\res(Af^{\eps})(x_0e_{\b_0})=|a|^{\frac{n+3}{2}}\lim_{y\to 0} Af^{\eps}(y,x_0,0)=
|a|^{\frac{n+3}{2}}\lim_{y\to 0}f^{\eps}(-x_0,y,0)=|a|^{\frac{n+3}{2}}.$$
It follows that 
$$\th(Af^{\eps})=c_2\cdot|a|^{\frac{n+3}{2}}\cdot|\mu(K)|$$
Hence, from equation (\ref{Aequation}) we obtain $c_2=|a|^{-\frac{n+3}{2}}$.
\ed

\section{Automorphic functional for $G=E_k$}\label{globalEsec}

Throughout this section we assume that that $G$ is of the type 
$E_k$, $k=6,7,8$, and we set $s=1,2,4$ in these three cases as before. 

\subsection{More on the form of $G(K)$-invariant functional}\label{formEsec}
 
Recall that in the case $G=E_k$ the line bundles $\LL_v$ over $C(K_v)$ have canonical trivialization, so the
$P(K_v)$-modules $\Ga(C(K_v),\LL_v)$ are identified with the spaces of functions on $C(K_v)$ 
with the action of $P(K_v)$ twisted by $|\om|^{s+1}$ (see Corollary \ref{trivlinecor}).
This allows to rewrite the functional $\th_2$ more explicitly.

Using the trivialization of $\LL_v$ we can consider the map $\res_v$ as a $P(K_v)$-equivariant map
$$\res_v:R_v\ra |\om|^{s+1}\otimes\C(C(K_v)).$$
On the other hand, using Corollary \ref{trivlinecor} and isomorphism of $R_v$ with
standard realization, we obtain another $P(K_v)$-equivariant map between these
spaces: $f\mapsto \ov{f}^{\norm}:=|x_0|^{s+1}_v\ov{f}$, where $\ov{f}$ is given by (\ref{fbareq}) and the
function $\ov{f}^{\norm}$ is extended uniquely from the open subset 
$\La^0(K_v)\sub C(K_v)$ to a locally constant function on $C(K_v)$.

To compare these two maps it is enough to compare
$\res_v(f)(e_{\b_0})$ with $\ov{f}(e_{\b_0})$. The formula (\ref{resnormal}) immediately implies that
$$\res_v(f)(x_0,x)=|a_v|^{\frac{n+3}{2}}\cdot\ov{f}^{\norm}(x_0,x)$$
on $\La^0(K_v)\sub C(K_v)$.
Thus, defining $\ov{f}^{\norm}$ as the product of local factors for $f=\otimes_v f_v$ 
we get
\begin{equation}
\th_2(f)=|a|^{\frac{n+3}{2}}\sum_{z\in C(K)}\ov{f}^{\norm}(z).
\end{equation}

Thus, the formula of Theorem \ref{mainglobalthm}
for the $G(K)$-invariant tempered functional can be rewritten as follows:
\begin{equation}\label{invE}
\th(f)=\sum_{y\in K^*,(x_0,x)\in K^{n+1}}f(y,x_0,x)+\sum_{z\in C(K)}\ov{f}^{\norm}(z)+\a(f).
\end{equation}

In fact, we can say more about the last term $\a$.
Let us assume that the conjecture of section \ref{archsec} is true 
for the archimedian places of $K$.
Using Theorem \ref{structurethm} we can write
\begin{equation}
\a=\a_1+\a_2,
\end{equation}
where the functionals $\a_1$ and $\a_2$ are $H(\mA)$-invariant,
$\a_1$ (resp., $\a_2$) is an eigenvector for the action of $Z(L)(\mA)$
with the eigenvalue $|\om|^{s+1}$ (resp., $|\om|^{\frac{s}{2}+1}$).
Next, we claim that the functional $\a_1$ and $\a_2$ are proportional to the following
two natural functionals. For every place $v$
we have projections $J(R_v)\ra\C$ and $J(R_v)\ra R_{G_1,v}$ 
corresponding to the decomposition
of $J(R_v)$ into irreducible components (see Theorem \ref{structurethm}). 
We normalize these projections in such a way that the image of the spherical vector
under the first (resp., second) projection is the element $1\in\C$
(resp., spherical vector in $R_{G_1,v}$). 
Taking the product of these projections over all places we obtain
surjective maps $\th_{\triv}:J(R_{\mA})\ra\C$ and $J(R_{\mA})\ra R_{G_1,\mA}$.
Let us denote by $\th_{G_1}:J(R_{\mA})\ra C$
the composition of the latter map with the $G_1(K)$-invariant functional
on $R_{\G_1,\mA}$. It is easy to see that $\th_{\triv}$ and $\th_{G_1}$
are $L(K)$-invariant
functionals on $J(R_{\mA})$. Our claim follows from the following result.

\begin{prop}\label{2funprop} 
The space of $G_1(K)$-invariant tempered
functionals on $J(R_{\mA})$ is spanned
by $\th_{\triv}$ and $\th_{G_1}$.
\end{prop}

\Pf . Let us denote by $J'(R_{\mA})\sub J(R_{\mA})$ 
the kernel of the projection
$J(R_{\mA})\ra R_{G_1,\mA}$. By Corollary \ref{uniquecor} it suffices
to prove that the space of $G_1(K)$-invariant tempered
functionals on $J'(R_{\mA})$
is spanned by the restriction of $\th_{\triv}$. Let us denote by
$f_{1,v}\in J(R_{\mA})$ the element corresponding to $(1,0)$ under the
isomorphism $R_v\simeq\C\oplus R_{G_1,v}$. Then $J'(R_{\mA})$ is generated
by tensor products of vectors in $J(R_v)$ over all $v$ with $f_{1,v}$ 
in at least one place (and with the spherical vector at almost all places).
Now let $\th'$ be a $G_1(K)$-invariant functional on $J'(R_{\mA})$.
Adding to $\th'$ a scalar multiple of $\th_{\triv}$ we can assume
that $\th'(\otimes_v f_{1,v})=0$. We claim that this implies that $\th'=0$.
Indeed, let us fix a place $v_0$ and consider the subspace
$$J'_{v_0}(R_{\mA})=f_{1,v_0}\otimes J(R_{\mA^{v_0}})\sub J'(R_{\mA})$$
where $R_{\mA^{v_0}}$ is the restricted product of $R_v$ over $v\neq v_0$.
It suffices to prove that the restriction of $\th'$
to this subspace is zero. We can consider this restriction as a $G(K)$-invariant
functional on $J(R_{\mA^{v_0}})$. But $G(K)$ is dense in $G(\mA^{v_0})$, so
it is in fact $G(\mA^{v_0})$-invariant. Hence, $\th'|_{J'_{v_0}(R_{\mA})}$ is
proportional to the restriction of $\th^{\triv}$. Since we assumed that
$\th'(\otimes_v f_{1,v})=0$ this implies that $\th'|_{J'_{v_0}(R_{\mA})}=0$.
\ed

It follows that $\a_1$ is proportional to $\th_{\triv}$ while $\a_2$ is proportional
to $\th_{G_1}$. Let $f_0=\otimes_v f_{0,v}$ be the 
product of spherical vectors over all places. 
By the definition $\th_{\triv}(f_0)=1$. Therefore, $\a_1=\a_1(f_0)\th_{\triv}$.
Similarly, if we could prove that $\a_2(f_0)\neq 0$, this would imply that
$\th_{G_1}(f_0)\neq 0$ and then we would have $\a_2=\a_2(f_0)\th_{G_1}/\th_{G_1}(f_0)$.
Therefore, our study of the automorphic functional on $R_{\mA}$ essentially
reduces to finding the constants $\a_1(f_0)$ and $\a_2(f_0)$.
In the next section we will develop
a method for computing $\a_1(f_0)$ and $\a_2(f_0)$ assuming that $K$ is a functional field. 

Now let us further rewrite the first two terms of the formula (\ref{invE}) for
$f$ of the form 
$f=f_0^S\otimes f_S$, where $S$ is a finite set of places of $K$ containing all archimedian places,
$f_0^S$ is the product of spherical
vectors over all places $v$ such that $v\not\in S$, $f_S$ is an element of 
$R_{S}=\prod_{v\in S} R_v$. Let us assume also that for all $v\not\in S$ the character
$\psi_v$ has order $0$. Then
our formula for the spherical vector implies that for every $v\not\in S$ one has
$$\ov{f_{0,v}}^{\norm}(z)=\frac{q_v^{s}|z|_v^{-s}-1}{q_v^s-1}\cdot\delta_{|z|_v\le 1},$$
where $z\in C(K_v)$, $|z|_v$ as the maximum of norms of all coordinates of $z$, $q_v$
is the number of elements in the residue field of $v$.
Let us define the divisor supported at the point $v\not\in S$ by setting
$$\div_v(z)=n_v(z)\cdot v,$$
where the integer $n_v\in\Z$ is the order of $z$ at $v$.
Then we can rewrite the above formula as
$$\ov{f_{0,v}}^{\norm}(z)=\phi_v(\div_v(z))$$
where the function $\phi_v$ on divisors supported at $v$ is defined by
$$\phi_v(n\cdot v)=\cases \frac{q^{(n+1)s}_v-1}{q^{s}_v-1}, & n\ge 0,\\
0, &\text{otherwise}.\endcases$$
This implies that the vector $f_0^S=\prod_{v\not\in S} f_{0,v}$ satisfies
\begin{equation}\label{divfor0}
\ov{f_0^S}^{\norm}(z)=\phi(\div^S(z))
\end{equation}
where $\div^S(z)=\sum_{v\not\in S} \div_v(z)$ and the function $\phi$ is defined by
$$\phi(D)=\prod_v \phi_v(D_v)$$
where $D_v$ is the part of $D$ supported at $v$. 

\begin{lem}\label{philem} One has
$$\phi(D)=\sum_{0\le D'\le D} |D'|^{-s},$$
where $|\sum_v n_v\cdot v|:=\prod_v q_v^{-n_v}$.
\end{lem}

\Pf .
Indeed, it suffices to check this for divisors supported at one place in which case
this is trivial.
\ed

The formula (\ref{divfor0}) implies that
\begin{equation}\label{f0Sfor}
\begin{array}{l}
f_0^S(y,x_0,x)=\prod_{v:|y|_v\le |x_0|_v}\psi_v(-\frac{I(x)}{yx_0})|\div^S(y,x_0)|^{-s-1}\times\\
\phi(\div^S(y,x_0),x,\div^S(I'(x))-\div^S(y,x_0),\div^S(I(x))-2\div^S(y,x_0))
\end{array}
\end{equation}
where the divisor associated to a vector is set to be the minimum of the divisors of its component,
$\phi(D_1,\ldots,D_k)=\phi(\min(D_1,\ldots,D_k))$.
Therefore, we obtain
\begin{align*}
&\th(f_0^S\otimes f_S)=\sum_{y\in K^*,x_0,\ldots,x_n\in K}\ \ 
\prod_{v:|y|_v\le |x_0|_v}\psi_v(-\frac{I(x)}{yx_0})|\div^S(y,x_0)|^{-s-1}\times\\
&\phi(\div^S(y,x_0),x,\div^S(I'(x))-\div^S(y,x_0),\div^S(I(x))-2\div^S(y,x_0))f_S(y,x_0,x)+\\
&\sum_{z\in C(K)}\phi(\div^S(z))\ov{f_S}^{\norm}(z)+\a_1(f_0^S\otimes f_S)+\a_2(f_0^S\otimes f_S).
\end{align*}
Note that in the first sum the non-zero terms correspond only to $y,x_0,\ldots,x_n$
that belong to the ring of $S$-integers $\OO^S:=K\cap\prod_{v\not\in S}\OO_v$, such that
$\div^S(I'(x))\ge\div^S(y,x_0)$, $\div^S(I(x))\ge 2\div^S(y,x_0)$.  

Let us consider the particular case when $K=\Q$, $S=\infty$ is the only archimedian place of $\Q$,
$\psi$ is the product of standard additive characters $\psi_p$ of $\Q_p$ and of $\exp(-it)$. 
In this case $f_0^S=f_{0,\fin}$ is the product of spherical vectors over all finite places,
$f_S=f_{\infty}$ is an arbitrary element of $R_{\infty}$. Then  we can rewrite the above formula as
\begin{align*}
&\th(f_{0,\fin}\otimes f_{\infty})=\sum_{y\in\Z\setminus 0;x_0,\ldots,x_n\in\Z}\ \ 
\prod_{p:|y|_p\le |x_0|_p}\psi_p(-\frac{I(x)}{yx_0})\times\\
&|\gcd(y,x_0)|^{s+1}
\mu_s(\gcd(y,x_0),x,\frac{I'(x)}{\gcd(y,x_0)},\frac{I(x)}{\gcd(y,x_0)^2})f_{\infty}(y,x_0,x)+\\
&\sum_{z\in C(\Z)}\mu_s(z)\ov{f_{\infty}}^{\norm}(z)+\a_1(f_{0,\fin}\otimes f_{\infty})+
\a_2(f_{0,\fin}\otimes f_{\infty}),
\end{align*}
where $\mu_s(a)=\sum_{d|a}d^s$, $\mu_s(a_1,\ldots,a_n)=\mu_s(\gcd(a_1,\ldots,a_n))$,
$\gcd$ stands for the greatest common divisor.
Recall that the automorphic form on $G(\mA)$ corresponding to $\th$
is given by $\varphi(g)=\th(gf_0)$ where $f_0$ is the product of spherical vectors over all
places. Since $\varphi$ is $G(\Q)$-left-invariant and right-invariant with
respect to the maximal compact subgroups at all primes, it is uniquely
recovered from its restriction to $G(\R)\sub G(\mA)$. This restriction is the function
$\varphi(g)=\th(f_{0,\fin}\otimes gf_{0,\infty})$ on $G(\R)$ which is left-invariant with respect to
$G(\Z)\sub G(\R)$ and right-invariant with respect to the maximal compact
subgroup of $G(\R)$. The above formula gives an expression of $\varphi(g)$ with the first two
terms being the weighted sum of $gf_{0,\infty}$ over the integer lattice and the weighted
sum of $\ov{gf_{0,\infty}}^{\norm}$ over $C(\Z)$. Using the explicit formula for the
archimedian spherical vector $f_{0,\infty}$ and the formulas for the action of the Borel
subgroup of $G$ on $R$, one can write explicitly the first two terms of $\varphi(b)$ for
$b$ in the Borel subgroup of $G(\R)$. Also, we know that $\a_1(bf_0)=|\om(b)|^{s+1}\a_1(f_0)$ where
the character $\om:T\ra\G_m$ is extended to the Borel subgroup of $G$, so the only remaining problems
are to compute $\a_1(f_0)$ and to understand the term corresponding to $\a_2$. 
We will return to these problems in a sequel to this paper.


\subsection{Case of the functional field and $G=E_k$}

Henceforward, we assume that $K$ is a functional field
corresponding to a smooth projective curve $X$ over a finite field. 
Let $\F_q=H^0(X,\OO_X)$ be the
field of constants for $X$, $|X|$ the set of places (i.e. the set of
points of $X$ over an algebraic closure of $\F_q$ modulo the
action of the Galois group). When we talk about divisors on $X$ we always
mean divisors defined over $\F_q$.
We will prove that the constants $\a_1(f_0)$ and
$\a_2(f_0)$ are uniquely determined by the curve $X$. Moreover, we
will compute them in the case $X=\P^1$.

In the functional case it is convenient to twist our local models by the fibers of
the canonical bundle of $C$. 
Recall (see \ref{twistedsec}) that in the local
picture we can start with an arbitrary free $\OO$-module 
$M$ of rank $1$ and a character $\psi:F\otimes_{\OO} M\ra\C^*$ which is trivial
on $M$ and non-trivial on $\pi^{-1}M$. Then we can
realize the minimal representation of the local group $G(F)$ in functions of
$y,x_0,\ldots,x_n\in F\otimes_{\OO} M$ and our formula for 
the spherical vector is still valid in this realization (where
the norm on $F\otimes_{\OO}M$ is determined by the condition
that $M$ is the unit ball around $0$). 
Now for every invertible sheaf of $\OO_X$-modules $\MM$ on $X$ we
can replace adeles $\mA$ by the $\mA$-module
$$\MM_{\mA}=\prod_v (K_v\otimes_{\OO_{X,v}}\MM_v)$$
where $\MM_v$ is the stalk of $\MM$ at $v$,
the product is restricted with respect to the subgroups
$\hat{\MM}_{v}=\hat{\OO}_{X,v}\otimes_{\OO_{X,v}}\MM_{v}$.
The analogue of rational adeles in this picture is 
$\MM_K\sub\MM_{\mA}$:this is the stalk of $\MM$ at the general point of $X$
(the space of rational sections of $\MM$ on $X$).
Note that if we use local norms $|\cdot|_v$ on $K_v\otimes_{\OO_{X,v}}\MM_{v}$ 
for which $\hat{\MM}_{v}$ are unit balls, then
for $s\in\MM_K$ one has $\prod_v |s|_v=q^{-\deg\MM}$.
For $\MM=\om_X$ we denote the twisted adeles by $\om_{\mA}$ and twisted
rational ad`eles by $\om_K\sub\om_{\mA}$.
Let us choose a non-trivial character
$\psi_0:\F_q\ra\C^*$ and define the character $\psi:\om_{\mA}\ra\C^*$
by the formula 
$$\psi((\om_v))=\psi_0(\sum_{v\in |X|}\Tr_v\Res_v(\om_v))$$ 
where $\Res_v$ is the residue map with values in the residue field
$k(v)$, $\Tr_v:k(v)\ra\F_q$ is the trace map.
By the residue theorem, $\psi$ is trivial on $\om_K$.
We define the global minimal representation $R_{\mA}$ of $G(\mA)$ as the restricted
tensor product of local minimal representations $R_v$, where
for every place $v\in |X|$ we use the realization of the minimal
representation of $G(K_v)$ associated with the module
$\hat{\om}_{X,v}$.
We can think about $R_{\mA}$ as a subspace in the space of functions
$f(y,x_0,x_1,\ldots,x_n)$ where $x_i\in \om_{\mA}$, $y\in\om_{\mA}^*$.

To identify this picture with the previous way of twisting one just
has to choose a non-zero rational $1$-form $\eta\in\om_K$.
Then the map $f\mapsto f(y\cdot\eta,x_0\cdot\eta,\ldots,x_n\cdot\eta)$
will give an equivalence from the twisted realization that was defined previously
(for the character $\psi\circ\eta$ of $\mA$) with this one.

In our new twisted picture the functionals $\th_1$ and $\th_2$ assume the following form:
$$\th_1(f)=\sum_{y\in \om_K\setminus\{0\},(x_0,x)\in \om_K^{n+1}}f(y,x_0,x),$$
$$\th_2(f)=q^{(2g-2)(s+1)}\cdot\sum_{z\in C(\om_K)}\ov{f}^{\norm}(z),$$
where $\ov{f}^{\norm}$ is defined in the same way as before (but with the new meaning of
the norms $|\cdot |_v$). Here the constant $q^{(2g-2)(s+1)}$ appears as a ratio
$|x_0|^{s+1}/|x_0\eta|^{s+1}$ for $x_0\in K^*$, $\eta\in\om_K\setminus\{0\}$.
The analogue of formula (\ref{divfor0}) is
\begin{equation}\label{divfor}
\ov{f_0}^{\norm}(z)=\phi(\div(z)),
\end{equation}
where $\div(z)$ is the minimum of the divisors of the coordinates of $z\in C(\om_K)$, 
\begin{equation}\label{phifor}
\phi(D)=\sum_{0\le D'\le D}q^{s\deg D'}.
\end{equation}

Thus, the $G(K)$-invariant functional in this realization should take form
\begin{align}\label{functtwisted}
\th(f)=&\sum_{y\in \om_K\setminus\{0\},(x_0,x)\in \om_K^{n+1}}f(y,x_0,x)+
q^{(2g-2)(s+1)}\cdot\sum_{z\in C(\om_K)}\ov{f}^{\norm}(z)\nonumber\\
&+\a_1(f)+\a_2(f),
\end{align}
where $g$ is the genus of $X$.
Since the $G(K)$-invariant functional $\th$ as in (\ref{functtwisted})
is unique, the constants $\a_1(f_0)$ and $\a_2(f_0)$ 
depend only on $s\in\{1,2,4\}$ and on the curve $X$. We are going
to describe a simple general method that could be used 
to compute these constants. Then we will apply this method
in the case $X=\P^1$. The case of elliptic curve will be considered
in \ref{ellipticsec}. 

The idea is to apply the functional $\th$ to
the elements of the form $tf_0$, where $t\in Z(L)(\mA)$ and to use
the fact that $\a_1$ and $\a_2$ are eigenvectors with respect to
the action of $Z(L)(\mA)$. Using the isomorphism 
$\om_{\b_0}^{\vee}:\G_m\ra Z(L)$ we can map ideles to elements of 
$Z(L)(\mA)$. Now for every idele $a\in\mA^*$ we have
$$\om_{\b_0}^{\vee}(a)f_0(y,x_0,x)=|a|^{3s+3}f_0(a^2y,ax_0,ax),$$
$$\ov{\om_{\b_0}^{\vee}(a)f_0}^{\norm}(z)=|a|^{2s+2}\ov{f_0}^{\norm}(az).$$
Therefore, if $D$ is the divisor of $a$ then
\begin{align*}
&\th(\om_{\b_0}^{\vee}(a)f_0)=
|a|^{3s+3}\sum_{y\in H^0(X,\om_X(2D))\setminus\{0\},
(x_0,x)\in H^0(X,\om_X(D))^{n+1}}f_0^D(y,x_0,x)+\\
&q^{(2g-2)(s+1)}|a|^{2s+2}\sum_{z\in C(H^0(X,\om_X(D)))}\ov{f_0}^D(z)
+|a|^{2s+2}\a_1(f_0)+|a|^{s+2}\a_2(f_0)
\end{align*}
where $\ov{f_0}^D$ and $f_0^D$ are some modifications of the functions
$\ov{f_0}^{\norm}$ and $f_0$ defined as follows. The function 
$\ov{f_0}^D(z)=\prod_v \ov{f_0}^D_v(z)$ for $z\in C(\om_X(D)_{\mA})$
is given by the formula (\ref{divfor}), the only difference
is that the line bundle $\om_X$ is replaced by $\om_X(D)$. 
The function $f_0^D(y,x_0,x)$
for $y\in\om_X(2D)_{\mA}^*$, $(x_0,x)\in\om_X(D)_{\mA}^{n+1}$ is the product
of local factors that are given by 
$$f_0^D(y,x_0,x)_v=\cases \psi_0(-\Res_v\frac{I(x)}{yx_0})|x_0|_v^{-s-1}
\ov{f_0}^D_v(x_0,x,\frac{I'(x)}{x_0},\frac{I(x)}{x_0^2}), & |y|_v\le |x_0|_v,\\
|y|_v^{-s-1}\phi_v(\div(y),\div(x),\div(\frac{I'(x)}{y}),\div(\frac{I(x)}{y^2}))), 
& |x_0|_v<|y|_v.\endcases
$$ 
Note that this formula makes sense since the ratio $I(x)/yx_0$ can be viewed as an
element of $\om_{\mA}$. Also in this formula we extend the definition of $\phi_v$
to several divisors by setting $\phi_v(D_1,\ldots,D_n)=\phi_v(\min(D_1,\ldots,D_n))$.
To define the relevant divisors, 
we view $I'(x)/y$ as an element of $\om_{\mA}$ and $I(x)/y^2$ as an element of 
$\om_X(-D)_{\mA}$. In order for the above formula to be well-defined we
introduce the divisor of zero $\div(0)$, which should be thought of as
infinitely large, so that $\min(D,\div(0))=D$ for all $D$.

We see that the sums
$$\th_1^D(f_0):=\sum_{y\in H^0(X,\om_X(2D))\setminus\{0\},
(x_0,x)\in H^0(X,\om_X(D))^{n+1}}f_0^D(y,x_0,x),$$
$$\th_2^D(f_0):=\sum_{z\in C(H^0(X,\om_X(D)))}\ov{f_0}^D(z)$$
depend only on the linear equivalence class of $D$ and that
\begin{align*}
\th(\om_{\b_0}^{\vee}(a)f_0)=\th^D(f_0):=&
q^{-(3s+3)\deg(D)}\th_1^D(f_0)+
q^{(2s+2)(g-1-\deg(D))}\th_2^D(f_0)+\\
&q^{-(2s+2)\deg D}\a_1(f_0)+q^{-(s+2)\deg D}\a_2(f_0).
\end{align*}
On the other hand, it easy to check that the operator $(AS)^3\in G(K)$ satisfies the relation
$$(AS)^3t=t^{-1}(AS)^3$$
for all $t\in Z(L)(\mA)$. Since the functional $\th$ and the vector $f_0$
are invariant with respect to $(AS)^3$ we obtain
$$\th(tf_0)=\th((AS)^3tf_0)=\th(t^{-1}(AS)^3f_0)=\th(t^{-1}f_0)$$
for all $t\in Z(L)(\mA)$.
Hence,
\begin{equation}\label{diveq}
\th^D(f_0)=\th^{-D}(f_0).
\end{equation}
If the divisor $D$ satisfies $h^0(\om(-D))=h^0(\om(-2D))=0$
(e.g. if $\deg D>2g-2$) then
one has $\th_1^{-D}(f_0)=\th_2^{-D}(f_0)=0$.
Therefore, for such $D$ the equation (\ref{diveq}) becomes
\begin{align}\label{mainglobeq}
&q^{-(3s+3)\deg(D)}\th_1^D(f_0)+
q^{(2s+2)(g-1-\deg(D))}\th_2^D(f_0)+\nonumber\\
&q^{-(2s+2)\deg D}\a_1(f_0)+q^{-(s+2)\deg D}\a_2(f_0)=\nonumber\\
&q^{(2s+2)\deg D}\a_1(f_0)+q^{(s+2)\deg D}\a_2(f_0).
\end{align}
Clearly, these equations uniquely determine the constants $\a_1(f_0)$ and $\a_2(f_0)$
(it suffices to consider two such equations with distinct degrees of $D$).

In practice the computation of $\th_1^D(f_0)$ and $\th_2^D(f_0)$ becomes quite complicated
as degree of $D$ grows. The following proposition gives a different way of calculating
$\alpha_1(f_0)$ based on the fact that $\alpha_1$ is an eigenvector with respect to
the action of the entire torus $T(\mA)$.

\begin{prop}\label{alpha1prop} 
Let $a,a'$ be a pair of ideles such that $|aa'|=1$. Then
\begin{align*}
&(|a|-1)(|a|^s-1)(1-|a|^{-s-1})\alpha_1(f_0)=
(\th_1+\th_2)[(\b_0^{\vee}(a)^{-1}-\b_0^{\vee}(a')^{-1})\om_{\b_0}^{\vee}(aa')f_0\\
&+|a|^{s+2}(\frac{\om_{\b_0}^{\vee}(a')}{\b_0^{\vee}(a')}-\om_{\b_0}^{\vee}(a'))f_0-
|a|^{-s-2}(\frac{\om_{\b_0}^{\vee}(a)}{\b_0^{\vee}(a)}-\om_{\b_0}^{\vee}(a))f_0].
\end{align*}
\end{prop}

\Pf . For every idele $a$ let us denote $t_1(a)=\om_{\b_0}(a)$, 
$t_2(a)=\om_{\b_0}^{\vee}(a)\b_0^{\vee}(a)^{-1}$. Note that $s_{\b_0}(t_1(a))=t_2(a)$.
We consider $f\in R_{\mA}$ defined as follows:
$$f=(t_2(a')-|a'|^{s+2})(t_1(a)-|a|^{s+2})f_0.$$
Since $\alpha_2$ is the eigenvector for the action of $t_1$ with the eigenvalue $|a|^{s+2}$,
we obtain $\alpha_2(f)=0$. On the other hand, since 
the operator $A$ represents the action of $s_{\b_0}$,
we have
$$Af=(t_1(a')-|a'|^{s+2})(t_2(a)-|a|^{s+2})f_0,$$
hence $\alpha_2(Af)=0$.
Note also that since $\alpha_1$ is the eigenvector for the action of $T(\mA)$ with the
character $|\om|^{s+1}$, we have 
$$\alpha_1(f)=(|a'|^{s+1}-|a'|^{s+2})(|a|^{2s+2}-|a|^{s+2})\alpha_1(f_0)=(|a|-1)(|a|^s-1)\alpha_1(f_0),$$
$$\alpha_1(Af)=(|a'|^{2s+2}-|a'|^{s+2})(|a|^{s+1}-|a|^{s+2})\alpha_1(f_0)=
(|a|-1)(|a|^s-1)|a|^{-s-1}\alpha_1(f_0).$$
Using the equations $\th(Af)=\th(f)$ and $Af_0=f_0$ we obtain the desired formula for $\alpha_1(f_0)$.
\ed

Now we are going to apply the above methods to the case $X=\P^1$. 

\begin{thm} For $X=\P^1$ one has
$$\a_1(f_0)=\frac{1}{q^{2s+2}(q^{s+1}-1)(q^s-1)},$$
$$\a_2(f_0)=-\frac{q^{2s+1}+1}{q^{3s+2}(q^{s+2}-1)(q^s-1)}.$$
\end{thm}

\Pf .
We are going to consider the equation (\ref{mainglobeq}) for the divisor $D$ of degree $1$.
We have 
$\th_1^{-D}(f_0)=\th_2^{-D}(f_0)=0=\th_2^D(f_0)=0$ and 
$\th_1^D(f_0)=(q-1)$ since in this case $\om_{\P^1}(2D)\simeq\OO_{\P^1}$.
Thus, the equation (\ref{mainglobeq}) for $\deg(D)=1$ takes form
\begin{equation}\label{P1eq1}
q^{-3s-3}(q-1)+q^{-2s-2}\a_1(f_0)+q^{-s-2}\a_2(f_0)=
q^{2s+2}\a_1(f_0)+q^{s+2}\a_2(f_0).
\end{equation}
On the other hand, we can apply the formula of Proposition \ref{alpha1prop} to calculate
$\a_1(f_0)$. Let us take $a^{-1}=a'=\pi$, where $\pi$ is a uniformizing parameter at some
$\F_q$-point $v$. Leaving out zero terms in the right-hand side we can rewrite
this formula in our case as follows:
$$(q-1)(q^s-1)(1-q^{-s-1})\alpha_1(f_0)=
\th_2(\b_0^{\vee}(\pi)f_0)
+q^{s+2}[\th_2(\frac{\om_{\b_0}^{\vee}(\pi)}{\b_0^{\vee}(\pi)}f_0)-\th_1(\om_{\b_0}^{\vee}(\pi)f_0)],
$$
where
\begin{align*}
&\th_2(\b_0^{\vee}(\pi)f_0)=q^{-3s-3}(q-1),\\
&\th_2(\frac{\om_{\b_0}^{\vee}(\pi)}{\b_0^{\vee}(\pi)}f_0)=q^{-3s-3}(q-1),\\
&\th_1(\om_{\b_0}^{\vee}(\pi)f_0)=q^{-3s-3}(q-1).
\end{align*}
Therefore, the formula takes form
$$(q^s-1)(1-q^{-s-1})\a_1(f_0)=q^{-3s-3}.$$
It remains to use equation (\ref{P1eq1}) to find $\a_2(f_0)$.
\ed

\subsection{Elliptic curve case}
\label{ellipticsec}

Assume now that $X$ is an elliptic curve over $\F_q$ and let
us calculate the constants $\a_1(f_0)$ and $\a_2(f_0)$.

First, let us find one quantity that will often arise in the calculation.

\begin{lem}\label{conenumlem}
The number of points of $C$ over the finite field $\F_q$ is
$$\card C(\F_q)=(q^{3s+2}-1)(q^{2s+1}+1)(q^{s+1}+1).$$
\end{lem}

\Pf . The stabilizer subgroup $S_{\ga_0}\sub L$ of a point $e_{\ga_0}\in C$ is connected,
so the action of $L(\F_q)$ on $C(\F_q)$ is transitive and we have
$$\card C(\F_q)=\frac{\card L(\F_q)}{\card S_{\ga_0}(\F_q)}.$$
Using the notation of Lemma \ref{transitivelem} we have
$$\card S_{\ga_0}(\F_q)=\frac{\card Q_{\ga_0}(\F_q)}{q-1}=\card P_1(\F_q)=
q^l(q-1)\card G_0(\F_q)$$
where $l=\card\Phi_+(G_1)-\card\Phi_+(G_0)$.
On the other hand, 
$$\card L(\F_q)=(q-1)\card G_1(\F_q).$$ 
Hence,
$$\card C(\F_q)=\frac{\card G_1(\F_q)}{q^l\cdot\card G_0(\F_q)}.$$
It remains to use the table in section \ref{standsec} giving $G_0$ and $G_1$
and the formulas for the numbers of elements of
Chevalley groups over $\F_q$ (see \cite{St}, \S 9).
\ed

Let us consider equation (\ref{mainglobeq}) in the case $\deg D=1$.
We can assume that $D$ is effective, so $D=p$ for some point $p$ of degree
$1$ on $E$. Since $H^0(\OO_X(p))$ is one-dimensional, we have 
$\phi(\div(z))=\phi(p)=q^s+1$ for every $z\in C(H^0(\OO_X(p)))$.
Hence, 
$$\th_2^p(f_0)=(q^s+1)\cdot\card C(\F_q).$$

In the sum determining $\th_1^p(f_0)$ let us distinguish three cases:
(i) $x_0\neq 0$, $(y,x_0)=1$, i.e. $y$ is not divisible by $x_0$;
(ii) $x_0\neq 0$, $x_0| y$; (iii) $x_0=0$.

In case (i) we have 
$$f_0^p(y,x_0,x)=\psi_0(-\sum_{t:y(t)=0}\Res_t\frac{I(x)}{yx_0}).$$
Using the residue theorem we can rewrite this as
$$f_0^p(y,x_0,x)=\psi_0(\Res_{t:x_0(t)=0}\frac{I(x)}{yx_0})=1$$
since $x$ vanishes at the unique zero of $x_0$, hence $I(x)/yx_0$ is
regular at this point. The corresponding part of the sum determining
$\th_1^p(f_0)$ is equal to
$$\card\{(y,x_0,x):\ (y,x_0)=1\}=(q-1)q^{6s+3}(q-1)q$$ 
(we used the fact that $n=6s+3$ and the fact that the number of
$y\in H^0(\OO_X(2p))$ that are relatively prime to a given non-zero element
$x_0\in H^0(\OO_X(p))$ is equal to $q(q-1)$).

In case (ii) we have
$$f_0^p(y,x_0,x)=|x_0|^{-s-1}\phi(\div(x_0),\div(x),\div(\frac{I'(x)}{x_0}),
\div(\frac{I(x)}{x_0^2})).$$
Note that $x$ is divisible by $x_0$, hence, $I'(x)$ is divisible by $x_0^2$,
$I(x)$ is divisible by $x_0^3$. Therefore,
$$f_0^p(y,x_0,x)=|x_0|^{-s-1}\phi(\div(x_0))=q^{s+1}\phi(p)=q^{s+1}(q^s+1).$$
The corresponding part of the sum calculating $\th_1^p(f_0)$ is equal to
$$q^{s+1}(q^s+1)\cdot\card\{(y,x_0,x):\ y\neq 0, x_0\neq 0, x_0|y\}=
q^{s+1}(q^s+1)q^{6s+3}(q-1)^2.$$

In case (iii) we have
$$f_0^p(y,x_0,x)=|y|^{-s-1}\phi(\div(y),\div(x),\div(\frac{I'(x)}{y}),
\div(\frac{I(x)}{y^2})).$$
In order for this expression to be non-zero $I(x)$ should be
divisible by $y^2$. But $I(x)$ is a section of $\OO_X(3p)$ and
$y^2$ is a section of $\OO_X(4p)$, hence we should have $I(x)=0$.
Now we can distinguish the following three subcases:
(iiia) $I'(x)\neq 0$; (iiib) $I'(x)=0$, $x\neq 0$; (iiic) $x=0$.

Since $I'(x)$ should be divisible by $y$, in case (iiia) we obtain that
$I'(x)$ is a non-zero multiple of $y$. Equivalently, $y$ belongs to the
image of the map $S^2H^0(\OO_X(p))\ra H^0(\OO_X(2p))$, which is one-dimensional. Thus,
the corresponding part of the sum for $\th_1^p(f_0)$ is equal to
\begin{align*}
&q^{2s+2}(q-1)\card\{x\in\F_q^n: I(x)=0, I'(x)\neq 0\}=\\
&q^{2s+2}(q-1)[q^{6s+3}-\card\{x\in\F_q^n: I(x)\neq 0\}-\card\{x\in\F_q^n:I'(x)=0\}].
\end{align*}
The latter quantity can be calculated using formulas (\ref{card1}) and (\ref{card2}).

The part of the sum corresponding to case (iiib) is
$$q^{2s+2}\sum_{y\neq 0,x\neq 0:I'(x)=0}\phi(\div(y),\div(x)).$$
Since $\div(x)=p$, we get $\phi(\div(y),\div(x))=1$ if $y(p)\neq 0$ and
$\phi(\div(y),\div(x))=q^s+1$ if $y(p)=0$. Thus, the above expression is equal to
\begin{align*}
&q^{2s+2}[(q^s+1)\card\{y\neq 0:\ y(p)=0\}+\card\{y\neq 0:\ y(p)\neq 0\}]\card\{x\neq 0:I'(x)=0\}\\
&=q^{2s+2}[(q^s+1)(q-1)+q(q-1)]\card\{x\in\F_q^n:\ x\neq 0, I'(x)=0\}.
\end{align*}
The latter cardinality is given by the formula (\ref{card2}).

Finally, the part of the sum corresponding to case (iiic) is
$$q^{2s+2}\sum_{y\in H^0(\OO_X(2p))\setminus\{0\}}\phi(\div(y)).$$
Using formula (\ref{phifor}) we can rewrite this as follows
$$q^{2s+2}\sum_{D\ge 0} q^{s\deg D}(q^{h^0(2p-D)}-1)=q^{2s+2}[(q^2-1)+q^s(q-1)N+q^{2s}(q^2-1)],$$
where $N$ is the number of $\F_q$-points of $X$.

Combining all the terms of the sum for $\th_1^p(f_0)$ we obtain
\begin{align*}
&\th_1^p(f_0)=q^{6s+4}(q-1)^2+q^{7s+4}(q-1)^2(q^s+1)+\\
&q^{2s+2}(q-1)[q^{6s+3}-1-\card\{x\in\F_q^n: I(x)\neq 0\}+\\
&(q^s+q)\card\{x\in\F_q^n\setminus 0: I'(x)=0\}+q^s N+(q^{2s}+1)(q+1)]=\\
&q^{2s+2}(q-1)[q^{4s+2}(q-1)+q^{5s+2}(q-1)(q^s+1)\\
&+(q^{2s+1}-1)(q^{2s+1}+q+1)(q^{2s}+q^s+1)+(q^{2s}+1)(q+1)+q^s N].
\end{align*}
Therefore, equation (\ref{mainglobeq}) in this case becomes
\begin{align}\label{elleq1}
&(q^{2s+2}-q^{-2s-2})\a_1(f_0)+(q^{s+2}-q^{-s-2})\a_2(f_0)=\nonumber\\
&q^{-s-1}(q-1)[q^{4s+2}(q-1)+q^{5s+2}(q-1)(q^s+1)+\nonumber\\
&(q^{2s+1}-1)(q^{2s+1}+q+1)(q^{2s}+q^s+1)+(q^{2s}+1)(q+1)+q^s N]+\nonumber\\
&q^{-2s-2}(q^{3s+2}-1)(q^{2s+1}+1)(q^{s+1}+1)(q^s+1),
\end{align}

To get the second equation let us pick a divisor $D$ of degree $2$ on $X$ and let us compute
$\th_1^D(f_0)$ and $\th_2^D(f_0)$. One can also try to calculate $\a_1(f_0)$ using
Proposition \ref{alpha1prop} but we are not going to pursue this.

We will need the following auxiliary Lemmas.

\begin{lem}\label{ellnumlem1} 
Let $\P C\sub\P(V)$ be the projectivization of $C$. Then the number of projective lines
in $\P(V)$ that are contained in $\P C$ is equal to 
\begin{align*}
&\frac{\card\{x\in\F_q^n\setminus\{0\}: I'(x)=0\}\cdot\card C(\F_q)}{(q-1)(q^2-1)}=\\
&\frac{(q^{2s+1}-1)(q^{2s}+q^s+1)(q^{3s+2}-1)(q^{2s+1}+1)(q^{s+1}+1)}{(q-1)(q^2-1)}.
\end{align*}
\end{lem}

\Pf . Since the action of $L(\F_q)$ on $C(\F_q)$ is transitive we have
$$\card\{p\in L\sub\P C\}=\card\P C(\F_q)\card\{L\sub\P C:e_{\b_0}\in L\}.$$
On the other hand, clearly,
$$\card\{p\in L\sub\P C\}=\card\{L\sub\P C\}(q+1).$$
Therefore,
$$\card\{L\sub\P C\}=\card C(\F_q)\card\{L\sub\P C:e_{\b_0}\in L\}/(q^2-1).$$
It remains to notice that directions of lines passing through $e_{\b_0}$ that are
contained in $\P C$ correspond to directions of vectors $x\in\F_q^n$ such that
$I'(x)=0$. Hence,
$$\card\{L\sub\P C:e_{\b_0}\in L\}=\card\{x\in\F_q^n\setminus\{0\}: I'(x)=0\}/(q-1),$$
which implies our formula.

\begin{lem}\label{ellnumlem2} One has
$$\card\{x\in H^0(\P^1,\OO_X(1))^{\oplus n}\setminus\{0\}: I'(x)=0\}=
(q^{2s+1}-1)(q^{2s}+q^s+1)(q^{3s}+q^{2s+1}-q^s+1).$$
\end{lem}

\Pf . Let us denote by $Z\sub\P^{n-1}$ the projective variety given by the equations $I'(x)=0$.
Then the cardinality in question is closely related to the number of projective lines
that are contained in $Z$. More precisely, we have
\begin{align*}
&\card\{x\in H^0(\P^1,\OO_X(1))^{\oplus n}\setminus\{0\}: I'(x)=0\}=\\
&(q+1)\card\{x\in\F_q^n\setminus\{0\}:I'(x)=0\}+(q^2-1)(q-1)q\card\{\text{lines in } Z\}.
\end{align*}
Since $\card\{x\in\F_q^n\setminus\{0\}:I'(x)=0\}$ is given by formula (\ref{card2}) it
suffices to find the number of lines in $Z$ (defined over $\F_q$).

Let us denote by $\Lines(Z)$ the variety of projective lines in $Z$.
If $G=E_6$ (resp., $G=E_7$) then $Z$ is the variety of $3\times 3$-matrices of rank $1$
(resp. skew-symmetric $6\times 6$-matrices of rank $2$). So in these two cases
it is not difficult to compute the number of lines in $Z$ directly.
To compute this number in the case $G=E_8$ we will prove that 
$\Lines(Z)\simeq E_6/Q$ where $Q$ is the maximal parabolic subgroup
in $E_6$ corresponding to the vertex adjacent to the last vertex in the Dynkin diagram
(i.e., to the $5$-th vertex with respect to the standard numbering of vertices). 
First we claim that $Z$ is the projectivization of the orbit of the highest weight vector
in the $27$-dimensional representation of $E_6$, 
so $Z=E_6/P$ where $P$ is the parabolic subgroup corresponding to
the last vertex in the Dynkin diagram of $E_6$. 
Indeed, let $\la$ denote the fundamental weight corresponding
to this vertex, so that $V_{\la}$ is the
$27$-dimensional representation of $E_6$. The derivatives of the $E_6$-invariant cubic form $I$
on $V_{\la}$ correspond to the unique (up to scalar) invariant map $S^2(V_{\la})\ra V_{\la}^*$.
Therefore, the quadratic relations defining $Z$ correspond to the second factor in
the decomposition $S^2(V_{\la}^*)=V_{2\la}^*\oplus V_{\la}$. It is easy to see that the quotient
of the symmetric algebra $S^*(V_{\la}^*)$ by these relations is precisely the algebra
$\oplus_{n\ge 0}V_{n\la}^*$ of functions on the projectivization of the
orbit of the highest weight vector in $V_{\la}$.
Hence, $Z=E_6/P$. 

Let $W$ be the Weyl group of $E_6$, $W'\sub W$ be
the subgroup corresponding to $P$. One can check that the set of double cosets
$W'\backslash W/W'$ has $3$ elements. This means that $Z\times Z$ is partitioned in
three $E_6$-orbits. One of them is the diagonal $\Delta\sub Z\times Z$, another is an open orbit. 
Now the set of pairs $(z_1,z_2)$, such that $z_1\neq z_2$ and the line spanned by $z_1$ and $z_2$
is contained in $Z$, is a proper closed subset of $Z\times Z\setminus\Delta$. Therefore,
it is the third $E_6$-orbit. This implies that $E_6$ acts transitively on the set of 
lines in $Z$. It remains to find a line contained in $Z$ such that its stabilizer subgroup
in $E_6$ coincides with $Q$. To this end let us consider the subgroup $S=\SL_2\sub E_6$
corresponding to the simple root $\a_6$. Let $v_0\in V_{\la}$ be the highest weight vector.
Then $Sv_0$ is a plane in $V_{\la}$ and its projectivization $\P Sv_0$ is a line
contained in $Z$. 
For every root $\a$ let us denote by $U_{\a}$ the corresponding root subgroup in $E_6$. 
Note that $\P Sv_0$ is the closure of $U_{-\a_6}$-orbit of the line
spanned by $v_0$. For every positive root $\a$ the commutator of $U_{\a}$ and $U_{-\a_6}$
lies in the Borel subgroup, hence $U_{\a}$ preserves $\P Sv_0$. On the other hand,
if $\b$ is a root belonging to the subsystem
$A_4\sub E_6$ spanned by $\a_1,\ldots,\a_4$ then the corresponding root subgroup $U_{\b}\sub Q$
commutes with $U_{-\a_6}$, hence, it preserves $\P Sv_0$. Finally, the maximal torus
and the subgroup $S$ also preserve this line.   
It follows that the stabilizer of $\P Sv_0$
contains $Q$. Since $Q$ is a maximal parabolic, this stabilizer is equal to $Q$.

Thus, we have established that $\Lines(Z)=E_6/Q$.
Using the formula for $\card E_6(\F_q)$ we obtain
$$\card\Lines(Z)(\F_q)=\frac{(q^{12}-1)(q^9-1)(q^4+1)(q^3+1)}{(q-1)(q^2-1)}.$$
\ed

\begin{lem}\label{ellnumlem3} One has
\begin{align*}
&\card\{x\in H^0(\P^1,\OO_X(1))^{\oplus n}: \deg\div(I'(x))=2, I(x)=0\}=\\
&q^{2s+1}(q^{2s+1}-1)(q^{2s}+q^s+1)(q^{s+1}-1)(q^{s+1}+q^s+q^{s-1}+1)
\end{align*}
\end{lem}

\Pf . Let us distinguish three cases: (a) $\div(I'(x))=p_1+p_2$, where
$p_1,p_2\in\P^1(\F_q)$, $p_1\neq p_2$; (b) $\div(I'(x))=p$, where
$p\in\P^1(\F_{q^2})\setminus \P^1(\F_q)$; (c) $\div(I'(x))=2p$, where
$p\in\P^1(\F_q)$. We claim that in cases (a) and (b) the condition $I(x)=0$ is
satisfied automatically. Indeed, in these cases we have $\div(x)=0$ (otherwise $\div(I'(x))$
would contain a double point). Therefore, $x$ determines a projective line
$L\sub\P^{n-1}$. Let $Z\sub\P^{n-1}$ be the singular locus of the cubic hypersurface $I(x)=0$
(given by the equations $I'(x)=0$). Then $L$ passes through $Z$ at two distinct points
corresponding to $p_1$ and $p_2$ in case (a) (resp., to $p$ and $\Frob_q(p)$ in case (b)).
It follows that the restriction of the cubic form $I$ to $Z$ has two double zeroes.
Therefore, $L$ is contained in the hypersurface $I(x)=0$ which implies our claim.
In case (a) each $x$ is uniquely determined by $x(p_1)$ and $x(p_2)$, hence we obtain
\begin{align*}
&\card\{x\in H^0(\P^1,\OO_X(1))^{\oplus n}: \div I'(x))=p_1+p_2\}=\\
&\card\{x\in\F_q^n: I'(x)=0\}^2-\card\{x\in H^0(\P^1,\OO_X(1))^{\oplus n}: I'(x)=0\}.
\end{align*}
Similarly, in case (b) we have
\begin{align*}
&\card\{x\in H^0(\P^1,\OO_X(1))^{\oplus n}: \div I'(x))=p\}=\\
&\card\{x\in\F_{q^2}^n: I'(x)=0\}-\card\{x\in H^0(\P^1,\OO_X(1))^{\oplus n}: I'(x)=0\}.
\end{align*}
In case (c) we can have $\div(x)=0$ or $\div(x)=p$. In the former case
the line $L\sub\P^{n-1}$ corresponding to $x$ is tangent to the singular locus
$Z$ of the cubic hypersurface $I(x)=0$. This still implies that $L$ is contained in the
hypersurface $I=0$. Indeed, it is well-known that $Z$ is smooth (see \cite{EKP}).
Therefore, every tangent line to $Z$ is a limit of chords. But we have seen above that all chords of $Z$
are contained in $I=0$, hence every tangent line to $Z$ is contained in $I=0$.
Since the dimension of $Z$ is $4s$ (see \cite{EKP}), we obtain that for given $p\in\P^1(\F_q)$
\begin{align*}
&\card\{x\in H^0(\P^1,\OO_X(1))^{\oplus n}: x(p)\neq 0, I'(x)|_{2p}=0\}=\\
&q^{4s+1}\card\{x\in\F_q^n\setminus\{0\}: I'(x)=0\}.
\end{align*}
Therefore, for every $p\in\P^1(\F_q)$ one has
\begin{align*}
&\card\{x\in H^0(\P^1,\OO_X(1))^{\oplus n}: \div(I'(x))=2p, I(x)=0\}=\\
&\card\{x\in\F_q^n: I(x)=0\}+q^{4s+1}\card\{x\in\F_q^n\setminus\{0\}: I'(x)=0\}-\\
&\card\{x\in H^0(\P^1,\OO_X(1))^{\oplus n}: I'(x)=0\}.
\end{align*}
Combining the above contributions from three cases and using Lemma \ref{ellnumlem2} we get the result.
\ed

For every line bundle $M$ on a curve $Y$ let us denote by $C(M)$ the set of $z\in H^0(Y,M)\otimes V$
that satisfy the equations defining $C$ in $V$. 
Using formula (\ref{phifor}) we can write
\begin{align*}
&\th_2^D(f_0)=\sum_{z\in C(\OO_X(D))}\phi(\div(z))=\sum_{D'\ge 0}q^{s\deg D'}\card C(\OO_X(D-D'))=\\
&\card C(\OO_X(D))+ q^s\card C(\F_q)N+ q^{2s}\card C(\F_q)(q+1).
\end{align*}
Now we use the fact that the complete linear series $|D|$ defines a morphism $X\ra\P^1$
such that $D$ is the pull-back of $\OO_{\P^1}(1)$. Since the pull-back of sections defines
isomorphisms $H^0(\P^1,\OO_{\P^1}(m))\wt{\ra} \Sym^m H^0(X,\OO_X(D))\sub H^0(X,\OO_X(m D))$
for $m\ge 1$, it induced the bijection $C(\OO_{\P^1}(1))=C(\OO_X(D))$.
It is clear that
$$C(\OO_{\P^1}(1))=\card\{\text{lines in } \P C\}\card\GL_2(\F_q)+(q+1)\card C(\F_q).$$
Hence, applying Lemma \ref{ellnumlem1} we get
$$C(\OO_{\P^1}(1))=\card C(\F_q)[\card\{x\in\F_q^n\setminus\{0\}: I'(x)=0\}q+q+1].$$
Substituing this into the above expression for $\th_2^D(f_0)$ we get
\begin{align*}
&\th_2^D(f_0)=\card C(\F_q)[\card\{x\in\F_q^n\setminus\{0\}: I'(x)=0\}q+q+1+q^s N+ q^{2s}(q+1)]=\\
&(q^{3s+2}-1)(q^{2s+1}+1)(q^{s+1}+1)[q^{4s+2}+q^{3s+2}+q^{2s+2}+q^{2s}-q^{s+1}+1+q^s N].
\end{align*}

Let us distinguish four cases in the sum determining $\th_1^D(f_0)$: 
(i) $x_0\neq 0$, $(y,x_0)=1$;
(ii) $x_0\neq 0$, $(\div(y),\div(x_0))=p$ for some $\F_q$-point $p$ of $X$; 
(iii) $x_0\neq 0$, $x_0| y$; (iv) $x_0=0$.

In case (i) we have
$$f_0^D(y,x_0,x)=\psi_0(\sum_{t:x_0(t)=0}\Res_t\frac{I(x)}{yx_0}).$$
Since $I(x)$ (resp. $x_0$) is a pull-back of a section of $\OO_{\P^1}(3)$ 
(resp. $\OO_{\P^1}(1)$), either $I(x)$ is divisible by $x_0$ or $I(x)$ is relatively prime
to $x_0$. In the former case $f_0^D(y,x_0,x)=1$. In the latter case, it is convenient to
fix $x_0$ and $x$ (such that $(I(x),x_0)=1$) and to calculate the sum of $f_0^D(y,x_0,x)$ over all 
$y\in H^0(\OO_X(2D))\setminus\{0\}$ such that $y$ is relatively prime to $x_0$. 
We consider separately three subcases:
(ia) $\div(x_0)=v_1+v_2$, where $\deg(v_1)=\deg(v_2)=1$, $v_1\neq v_2$;
(ib) $\div(x_0)=v$, where $\deg(v)=2$; (ic) $\div(x_0)=2v$, where $\deg(v)=2$.
We claim that in case (ia)
$$\sum_{y:(y,x_0)=1} \psi_0(\Res_{v_1}\frac{I(x)}{yx_0}+\Res_{v_2}\frac{I(x)}{yx_0})=q^2;$$
in case (ib)
$$\sum_{y:(y,x_0)=1} \psi_0(\Res_{v}\frac{I(x)}{yx_0})=-q^2;$$
and in case (ic)
$$\sum_{y:(y,x_0)=1} \psi_0(\Res_{v}\frac{I(x)}{yx_0})=0.$$
Indeed, in case (ia) the values of $y$ at $v_1$ and $v_2$ can be arbitrary elements of 
$\F_q^*$ and each pair of values is assumed $q^2$ times. Similarly, in case (ib)
the value of $y$ at $v$ can be arbitrary element of $\F_{q^2}^*$ and each value is
assumed $q^2$ times. In case (ic), let us pick a uniformizing parameter $\pi_v$ at $v$,
and trivializations of the relevant line bundles such that the residue at $v$ can
be computed as coefficient with $\pi_v^{-1}$.
We can write $I(x)/x_0=(a_0+a_1\pi_v+\ldots)/\pi_v^2$, $y=b_0+b_1\pi_v+\ldots$,
where $a_0, b_0\in\F_q^*$, $a_1,b_1\in\F_q$. Then $\Res_v(I(x)/yx_0)=(a_1b_0-b_1a_0)/b_0^2$.
It remains to use the fact that in our sum over $y$ every pair $(b_0,b_1)$ is taken
$q^2$ times and for every fixed $b_0\in\F_q^*$ one has
$$\sum_{b_1\in\F_q}\psi_0(\frac{a_1b_0-b_1a_0}{b_0^2})=0.$$
Thus, the sum corresponding to case (ia) is equal to
\begin{align*}
&\sum_{x_0:\div(x_0)=v_1+v_2}[\card\{x\in H^0(\OO_X(D))^n: x_0|I(x)\}\card\{y:(y,x_0)=1\}+\\
&\card\{x\in H^0(\OO_X(D))^n:(I(x),x_0)=1\}q^2]
\end{align*}
Considering $x_0$ and coordinates of $x$ as elements of $H^0(\P^1,\OO(1))$ we see that
$$\card\{x\in H^0(\OO_X(D))^n: x_0|I(x)\}=q^{6s+3}\card\{x\in\F_q^n: I(x)=0\}.$$
Hence, we can rewrite this sum as
\begin{align*}
&\card\{x_0\in H^0(\OO_X(D)):\div(x_0)=v_1+v_2\}q^{6s+5}\times\\
&[\card\{x\in\F_q^n: I(x)=0\}(q-1)^2+\card\{x\in\F_q^n:I(x)\neq 0\}].
\end{align*}
Similarly, the sum corresponding to case (ib) is equal to
\begin{align*}
&\card\{x_0\in H^0(\OO_X(D)):\div(x_0)=v\}q^{6s+5}\times\\
&[\card\{x\in\F_q^n: I(x)=0\}(q^2-1)-\card\{x\in\F_q^n:I(x)\neq 0\}]
\end{align*}
and the sum corresponding to case (ic) is equal to
$$\card\{x_0\in H^0(\OO_X(D)):\div(x_0)=2v\}q^{6s+6}(q-1)
\card\{x\in\F_q^n: I(x)=0\}.$$
Let us denote by $a$, $b$ and $c$ the numbers of divisors in the linear series $|D|$ 
defined over $\F_q$ of the form $v_1+v_2$, $w$ and $2v$ respectively
(where $\deg(v_i)=\deg(v)=1$, $\deg(w)=2$).
Note that $a+b+c=q+1$ and $2a+c=N$. 
Combining cases (ia), (ib) and (ic) we get the following contribution from case (i):
\begin{align*}
&aq^{6s+5}(q-1)[\card\{x\in\F_q^n: I(x)=0\}(q-1)^2+\card\{x\in\F_q^n:I(x)\neq 0\}]+\\
&bq^{6s+5}(q-1)[\card\{x\in\F_q^n: I(x)=0\}(q^2-1)-\card\{x\in\F_q^n:I(x)\neq 0\}]+\\
&cq^{6s+6}(q-1)^2\card\{x\in\F_q^n: I(x)=0\}.
\end{align*}
Using the formula (\ref{card1}) and the relations between $a$, $b$ and $c$ this expression
can be rewritten as follows:
$$q^{9s+6}(q-1)^2[q(q+1)(q^{3s+1}+q^{2s+1}+q^{s+1}-1)+(1-q^{s+1}-q^{2s+1})N].$$

In case (ii) applying the residue theorem we can write
$$f_0^D(y,x_0,x)=\psi_0(\sum_{v:|x_0|_v<|y|_v}\Res_v\frac{I(x)}{yx_0})q^{s+1}
\phi(p,\div(x),\div(I'(x))-p,\div(I(x))-2p),$$
where $p$ is the greatest common divisor of $\div(x_0)$ and $\div(y)$.
In order for this to be non-zero $I'(x)$ should vanish at $p$. Since $I'(x)$ is a
pull-back of a section of $\OO_{\P^1}(2)$, this implies that $I'(x)$ is divisible
by $x_0$. We claim that in this case $I(x)$ is divisible by $x_0^2$.
(in particular, $\div(I(x))\ge 2p$). Indeed, this follows immediately
from the identity $I(I'(x))=I(x)^2$. It follows that the relevant residues vanish
since $\frac{I(x)}{yx_0}=\frac{I(x)}{x_0^2}\frac{x_0}{y}$ is regular at all places $v$
where $x_0/y$ is regular. Note also that we have
$$\phi(p,\div(x),\div(I'(x))-p,\div(I(x))-2p)=\cases 1, & x(p)\neq 0, I'(x)(p)=0,\\ q^s+1, & x(p)=0.
\endcases$$  
Therefore, the contribution from case (ii) is equal to
\begin{align*}
&q^{s+1}[\sum_{(y,x_0),p\in X(\F_q): (y,x_0)=p}\card\{x\in H^0(\OO_X(D))^n:x(p)\neq 0,I'(x)(p)=0\}+\\
&\sum_{(y,x_0),p\in X(\F_q): (y,x_0)=p}\card\{x\in H^0(\OO_X(D))^n:x(p)=0\}(q^s+1)]=\\
&q^{s+1}\card\{(y,x_0),p\in X(\F_q): (y,x_0)=p\}q^{6s+3}
[\card\{x\in\F_q^n\setminus\{0\}:I'(x)=0\}+q^s+1].
\end{align*}
To compute $\card\{(y,x_0),p\in X(\F_q): (y,x_0)=p\}$ we note that
for every point $p\in X(\F_q)$ there is a unique divisor $p+p'$ in $|D|$ containing $p$, hence
$x_0$ is determined up to a scalar. Next, we note that $y$ can be considered as an element of
$H^0(X,\OO_X(2D-p))$ which does not vanish at $p'$. Hence, there are $(q-1)$ choices of $x_0$
and $q^2(q-1)$ choices of $y$ for every given $p$, so 
$$\card\{(y,x_0),p\in X(\F_q): (y,x_0)=p\}=q^2(q-1)^2N.$$ 
Substituting this in the above expression and using formula (\ref{card2}) we can rewrite
the contribution from case (ii) as follows:
$$q^{9s+6}(q-1)^2(q^{2s+1}+q^{s+1}+q-1)N.$$

In case (iii) we have
$$f_0^D(y,x_0,x)=q^{2s+2}\phi(\div(x_0),\div(x),\div(I'(x)/x_0),\div(I(x)/x_0^2)).$$
In order for this to be non-zero $I'(x)$ should be divisible by $x_0$ (which automatically
implies that $I(x)$ is divisible by $x_0^2$). Note that either $x$ 
is divisible by $x_0$ or $\div(x)$ is relatively prime to $\div(x_0)$. In the latter case
$\phi(\div(x_0),\div(x),\ldots)=1$, while in the former case 
$\phi(\div(x_0),\div(x),\ldots)=\phi(\div(x_0))$. Therefore, the contribution from case (iii) is
equal to
\begin{align*}
&q^{2s+2}(q^2-1)[\card\{x_0\in H^0(\OO_X(D)), x\in H^0(\OO_X(D))^n: (x,x_0)=1, x_0|I'(x)\}+\\
&q^{6s+3}\sum_{x_0\in H^0(\OO_X(D))\setminus\{0\}}\phi(\div(x_0))].
\end{align*}
Applying (\ref{phifor}) we obtain
$$\sum_{x_0\in H^0(\OO_X(D))\setminus\{0\}}\phi(\div(x_0))=\sum_{D'\ge 0} q^{s\deg D'}(q^{h^0(D-D')}-1)=
q^2-1+q^s(q-1)N+q^{2s}(q^2-1).$$
Also we have
\begin{align*}
&\card\{x\in H^0(\OO_X(D)), x\in H^0(\OO_X(D))^n: (x,x_0)=1, x_0|I'(x)\}=\\
&(q^2-1)q^{6s+3}\card\{x\in F_q^n\setminus\{0\}: I'(x)=0\}.
\end{align*}
Therefore, using formula (\ref{card2}) we can rewrite the contribution
from case (iii) as follows:
$$q^{9s+5}(q^2-1)(q-1)[(q+1)(q^{3s+1}+q^{2s+1}+q^{s+1}-1)+N].$$

In case (iv) we have
$$f_0^D(y,x_0,x)=q^{4s+4}\phi(\div(y),\div(x),\div(I'(x)/y),\div(I(x)/y^2)).$$
In order for this to be non-zero $I(x)$ should vanish and $I'(x)$ should be divisible by $y$.
We distinguish four subcases: (iva) $I'(x)\neq 0$; (ivb) $I'(x)=0$, $\div(x)=0$;
(ivc) $I'(x)=0$, $x\neq 0$, $\div(x)>0$; (ivd) $x=0$.

The contribution from case (iva) is equal to
\begin{align*}
&q^{4s+4}(q-1)\card\{x\in H^0(\OO_X(D))^n:I(x)=0, \deg\div(I'(x))=4\}=\\
&q^{4s+4}(q-1)\card\{x\in H^0(\OO_{\P^1}(1))^n:I(x)=0, \deg\div(I'(x))=2\}.
\end{align*}
Applying Lemma \ref{ellnumlem3} we can rewrite this contribution as
$$q^{6s+5}(q-1)
(q^{2s+1}-1)(q^{2s}+q^s+1)(q^{s+1}-1)(q^{s+1}+q^s+q^{s-1}+1).$$

The contribution from case (ivb) is equal to
\begin{align*}
&q^{4s+4}(q^4-1)\card\{x\in H^0(\OO_X(D))^n:I'(x)=0, \div(x)=0\}=\\
&q^{4s+4}(q^4-1)[\card\{x\in H^0(\OO_{\P^1}(1))^n\setminus\{0\}:I'(x)=0\}-\\
&(q+1)\card\{x\in\F_q^n\setminus\{0\}:I'(x)=0\}.
\end{align*}
Using Lemma \ref{ellnumlem2} and formula (\ref{card2}) we can rewrite this
contribution as follows:
$$q^{4s+5}(q^4-1)(q^{2s+1}-1)(q^{2s}+q^s+1)(q^{2s}-1)(q^{s-1}+1).$$

The contribution from case (ivc) is equal to
\begin{align*}
&\sum_{y\in H^0(\OO_X(D))\setminus\{0\},D'\in |D|}q^{4s+4}\phi(\div(y),D')\times\\
&\card\{x\in H^0(\OO_X(D))^n\setminus\{0\}: I'(x)=0,\div(x)=D'\}=\\
&q^{4s+4}\card\{x\in F_q^n\setminus\{0\}: I'(x)=0\}\cdot
\sum_{y\in H^0(\OO_X(D))\setminus\{0\},D'\in |D|}\phi(\div(y),D').
\end{align*}
Applying formula (\ref{phifor}) we get
\begin{align*}
&\sum_{y\in H^0(\OO_X(D))\setminus\{0\},D'\in |D|}\phi(\div(y),D')=
\sum_{E\ge 0}q^{s\deg E}(q^{h^0(2D-E)}-1)\frac{q^{h^0(D-E)}-1}{q-1}=\\
&(q^4-1)(q+1)+q^s(q^3-1)N+q^{2s}(q^2-1)(q+1).
\end{align*}
Therefore, by (\ref{card2}) the contribution from case (ivc) is given by
$$q^{4s+4}(q^{2s+1}-1)(q^{2s}+q^s+1)[(q^2-1)(q+1)(q^{2s}+q^2+1)+q^s(q^3-1)N].$$

Finally, the contribution from case (ivd) is equal to
$$q^{4s+4}\sum_{y\in H^0(\OO_X(D))\setminus\{0\}}\phi(\div(y)).$$
Again applying (\ref{phifor}) we obtain
\begin{align*}
&\sum_{y\in H^0(\OO_X(D))\setminus\{0\}}\phi(\div(y))=
\sum_{E\ge 0}q^{s\deg E}(q^{h^0(2D-E)}-1)=\\
&q^4-1+q^s(q^3-1)N+q^{2s}(q^2-1)\card X^{(2)}(\F_q)+q^{3s}(q-1)\card X^{(3)}(\F_q)+q^{4s}(q^4-1),
\end{align*}
where $X^{(d)}$ denotes $d$-th symmetric power of $X$. Using the formula
$\card X^{(d)}(\F_q)=N(q^d-1)/(q-1)$ we can rewrite the contribution from case (ivd) as follows:
$$q^{4s+4}(q^4-1)(q^{4s}+1)+q^{5s+4}[(q^3-1)(q^{2s}+1)+q^s(q^2-1)(q+1)]N.$$

Combining all the above contributions we can write explicitly the equation (\ref{mainglobeq})
for the case $\deg D=2$. Together with the equation (\ref{elleq1}) this gives a system of equations 
that determines $\a_1(f_0)$ and $\a_2(f_0)$. Solving this system
we get the following answer.

\begin{thm} For an elliptic curve $X/\F_q$ one has
$$\a_1(f_0)=1+\frac{q^s N}{(q^{s+1}-1)(q^{s}-1)}$$
$$\a_2(f_0)=q^{2s+1}+1+\frac{q(q^{4s+2}-1)(q^s+1)}{q^{s+2}-1}-
\frac{(q^{2s+1}+1)N}{(q^{s+2}-1)(q^{s}-1)},$$
where $N$ is the number of $\F_q$-points of $X$.
\end{thm}

Note that in both cases $g=0$ and $g=1$ we have
$$\a_1(f_0)=q^{(2g-2)(s+1)}L(X,q^s).$$
where $L(X,t)=\frac{\det(1-t\Frob_q, H^1(X))}{(1-t)(1-qt)}$.
We conjecture that this formula holds for an arbitrary curve.
Using Proposition \ref{alpha1prop}
we can reformulate our conjecture as the following identity that
should hold for every ideles $a$ and $a'$, such that $|a|=q$, $|a'|=q^{-1}$:
\begin{align*}
&(\th_1+\th_2)[(\b_0^{\vee}(a)^{-1}-\b_0^{\vee}(a')^{-1})\om_{\b_0}^{\vee}(aa')f_0
+|a|^{s+2}(\frac{\om_{\b_0}^{\vee}(a')}{\b_0^{\vee}(a')}-\om_{\b_0}^{\vee}(a'))f_0-\\
&|a|^{-s-2}(\frac{\om_{\b_0}^{\vee}(a)}{\b_0^{\vee}(a)}-\om_{\b_0}^{\vee}(a))f_0]=
(q-1)q^{(2g-3)(s+1)}\det(1-q^s\Frob_q, H^1(X)).
\end{align*}

\end{document}